\documentclass[journal]{IEEEtran}
\usepackage[pass]{geometry}

\usepackage{scrextend}

\usepackage{footnote}
\makesavenoteenv{tabular} 
\makesavenoteenv{table}

\usepackage{tabularx,booktabs}

\newcolumntype{C}{>{\centering\arraybackslash}X} 
 \newcolumntype{b}{>{\centering\arraybackslash\hsize=2.3\hsize}X}
\newcolumntype{s}{>{\centering\arraybackslash\hsize=.45\hsize}X}
\newcolumntype{m}{>{\centering\arraybackslash\hsize=.9\hsize}X}

\usepackage{cuted}
\usepackage[linesnumbered,ruled,vlined]{algorithm2e}

\SetCommentSty{mycommfont}

\usepackage[noadjust]{cite}

\usepackage{float}
\usepackage{bbold}

\usepackage{blindtext}

\usepackage[pdfencoding=auto, psdextra]{hyperref}
\usepackage{bookmark}

\ifCLASSINFOpdf
   \usepackage[pdftex]{graphicx}
   \graphicspath{{./figure/}}
   \DeclareGraphicsExtensions{.pdf,.jpeg,.png}
\else

\fi
\usepackage[linesnumbered,ruled,vlined]{algorithm2e}
\usepackage{threeparttable}

\usepackage{rotating}
\usepackage{graphicx}
\newcommand*{\Scale}[2][4]{\scalebox{#1}{$#2$}}%
\usepackage[cmex10]{amsmath}
\usepackage{amsmath,amsfonts,amssymb}
\usepackage{eurosym}
\usepackage{epstopdf}
\usepackage{nomencl}

\usepackage{lipsum}
\usepackage{multirow}

\usepackage{booktabs,colortbl}
\usepackage{mathtools}
\DeclarePairedDelimiter\abs{\lvert}{\rvert}%

\usepackage{cuted}
\usepackage[dvipsnames]{xcolor}
\usepackage{hyperref} 
\hypersetup{colorlinks=true,                
    breaklinks=true,                
    urlcolor= black,
    urlbordercolor=blue,                
    linkcolor= blue,                
    bookmarksopen=false,
    filecolor=black,
    citecolor=green,
    linkbordercolor=blue
}
\AtBeginDocument{\hypersetup{pdfborder={0 0 1}}}
\usepackage{multirow}
\usepackage{fixltx2e}

\usepackage{tablefootnote}

\usepackage[perpage,symbol,stable]{footmisc}

\usepackage{upgreek}
\IEEEoverridecommandlockouts

\usepackage{etoolbox}
\renewcommand\nomgroup[1]{%
  \item[\bfseries
  \ifstrequal{#1}{A}{General}{%
  \ifstrequal{#1}{B}{Parameters}{%
  \ifstrequal{#1}{C}{Variables}{%
  \ifstrequal{#1}{D}{Indices}{}}}}%
]}
\DeclareSymbolFont{matha}{OML}{txmi}{m}{it}
\DeclareMathSymbol{\varv}{\mathord}{matha}{118}

\newcommand{\matr}[1]{\mathbf{#1}}

\usepackage{supertabular,capt-of}
\usepackage{multicol}
\usepackage{placeins}
\usepackage{supertabular}
\usepackage{array,supertabular,multirow}
\usepackage{float}

\usepackage{scalerel}
\usepackage{tikz}
\usetikzlibrary{svg.path}

\definecolor{orcidlogocol}{HTML}{A6CE39}
\tikzset{
  orcidlogo/.pic={
    \fill[orcidlogocol] svg{M256,128c0,70.7-57.3,128-128,128C57.3,256,0,198.7,0,128C0,57.3,57.3,0,128,0C198.7,0,256,57.3,256,128z};
    \fill[white] svg{M86.3,186.2H70.9V79.1h15.4v48.4V186.2z}
                 svg{M108.9,79.1h41.6c39.6,0,57,28.3,57,53.6c0,27.5-21.5,53.6-56.8,53.6h-41.8V79.1z M124.3,172.4h24.5c34.9,0,42.9-26.5,42.9-39.7c0-21.5-13.7-39.7-43.7-39.7h-23.7V172.4z}
                 svg{M88.7,56.8c0,5.5-4.5,10.1-10.1,10.1c-5.6,0-10.1-4.6-10.1-10.1c0-5.6,4.5-10.1,10.1-10.1C84.2,46.7,88.7,51.3,88.7,56.8z};
  }
}

\newcommand\orcidicon[1]{\href{https://orcid.org/#1}{\mbox{\scalerel*{
\begin{tikzpicture}[yscale=-1,transform shape]
\pic{orcidlogo};
\end{tikzpicture}
}{|}}}}

\begin{document}
\bstctlcite{IEEEexample:BSTcontrol}

\title{BATTPOWER Toolbox: Memory-Efficient and High-Performance Multi-Period AC Optimal Power Flow Solver }

\author{Salman Zaferanlouei  \orcidicon{0000-0001-6010-8919}\,, 
 Hossein Farahmand \orcidicon{0000-0003-1125-8296}\,,~\IEEEmembership{Senior Member,~IEEE,}
  Vijay Venu Vadlamudi \orcidicon{0000-0002-6771-1574}\,, ~\IEEEmembership{Member,~IEEE,}
   Magnus Korp\r{a}s\IEEEauthorrefmark{2} \orcidicon{0000-0001-5055-3912}\,,~\IEEEmembership{Member,~IEEE} 
\thanks{All the authors are with the Department of Electric Power Engineering, Norwegian University of Science and Technology, 7491 Trondheim, Norway (e-mails: salman.zaf@ntnu.no, hossein.farahmand@ntnu.no, vijay.vadlamudi@ntnu.no, magnus.korpas@ntnu.no).}
\thanks{\IEEEauthorrefmark{2} corresponding author}\
}

\maketitle

\begin{abstract}
With the introduction of massive renewable energy sources and storage devices, the traditional process of grid operation must be improved in order to be safe, reliable, fast responsive and cost efficient, and in this regard power flow solvers are indispensable. In this paper, we introduce an Interior Point-based (IP) Multi-Period AC Optimal Power Flow (MPOPF) solver for the integration of Stationary Energy Storage Systems (SESS) and Electric Vehicles (EV). The primary methodology is based on: 1) analytic and exact calculation of partial differential equations of the Lagrangian sub-problem, and 2) exploiting the sparse structure and pattern of the coefficient matrix of Newton-Raphson approach in the IP algorithm. Extensive results of the application of proposed methods on several benchmark test systems are presented and elaborated, where the advantages and disadvantages of different existing algorithms for the solution of MPOPF, from the standpoint of computational efficiency, are brought forward. We compare the computational performance of the proposed Schur-Complement algorithm with a direct sparse LU solver. The comparison is performed for two different applicational purposes: SESS and EV. The results suggest the substantial computational performance of Schur-Complement algorithm in comparison with that of a direct LU solver when the number of storage devices and optimisation horizon increase for both cases of SESS and EV. Also, some situations where computational performance is inferior are discussed.
\end{abstract}


\IEEEpeerreviewmaketitle
\section{Introduction}
\IEEEPARstart{L}{arge}-scale introduction of Renewable Energy Sources (RES), {Stationary} Energy Storage Systems ({SESS}) and Electric Vehicles (EV) will influence the way the electricity grid is operated. In the planning and operation of the electricity network, power flow analysis toolboxes that are reliable, computationally fast, and tractable are indispensable.\\
The optimal power flow is a non-linear and non-convex problem which was introduced in the sixties \cite{carpentier1962contribution} for the first time. Although it is considered to be a classic power systems problem among researchers, depending on the technical applications and operational dimensions, it may be adapted to various versions such as the Multi-Period AC Optimal Power Flow (MPOPF), introduced by \cite{chandy_simple_2010}, and may become intractable and computationally very demanding even after about 60 years \cite{sperstad_optimal_2016, Capitanescu_OPF_review}. \\ 
Many researchers have been trying to either simplify MPOPF by linearising the main problem \cite{carpinelli_optimal_2013, fortenbacher_optimal_2016,fortenbacher_modeling_2017}, or by making it more reliable by finding the global optimum point with different convex relation approaches such as \cite{geth_plangridev_2016}, semidefinite programming (SDP) relaxations \cite{warrington_market_2012, gopalakrishnan_global_2013} and second-order cone programming \cite{moghadasi_real-time_2014}. Moreover, MPOPF is being suggested as a potential online operational tool \cite{moghadasi_optimal_2016} for use in the near future.\\ 
Since the non-linear nature of full ACOPF problems requires non-linear solvers to be called, several Non-Linear Programming (NLP) solvers are used to solve MPOPF problems; these are primarily developed based on IP methods, such as MIPS\footnote{\label{first} { a solver name}} \cite{wang_computational_2007}, IPOPT\footref{first} \cite{wachter_implementation_2006} KNITRO\footref{first} \cite{byrd_knitro_2006}, and recently BELTISTOS\footref{first} \cite{kourounis_towards_2018}, of which the only tailored algorithm to solve MPOPF problems is BELTISTOS. An extensive review of both MPOPF problem formulations and solution methods can be found in \cite{sperstad_optimal_2016, sperstad_energy_2019}. With more and more penetration of renewables into the system, it becomes imperative to rely on SESS to overcome the variability of these renewable energy resources; MPOPF becomes a pivotal tool in this context for system operators, and as such better solvers need to be designed from the point of view of computational ease and efficiency. Though reliability is the most critical factor when it comes to the global optimum solution, computational performance is pivotal for the implementation and application of an algorithm for online operational purposes. It is well-known from the literature \cite{capitanescu_experiments_2013, castillo2013computational} that solution of linear Karush–Kuhn–Tucker (KKT) systems and calculation of gradients are the two most computationally expensive aspects in solving a MPOPF problem. Thus, we propose a solver to exploit the sparsity of a MPOPF structure (both KKT systems and gradients) and speed the solution up.\\   
Although considerable efforts have been undertaken in order to solve the MPOPF problem, no currently available work extensively elaborates the complex mathematical details of the problem, because of which it is difficult to compare the various solution algorithms in a systematic (same platform, same programming language, single-thread environment) manner. MATPOWER \cite{wang_computational_2007} does have a complete implementation of a single-period ACOPF with function evaluations and solution of KKT systems with MIPS. However, there is no package for the MPOPF function evaluations and for handling the sparsity structure of the first and second gradients, as introduced in this proposed package. Reference \cite{kourounis_towards_2018} introduced a fast multi-period KKT solver, Beltistos, as an extension of the IPOPT solver, and specifically for the MPOPF KKT structure, using the PARDISO solver. {However, reference \cite{kourounis_towards_2018} does not: 1) consider the sparsity structure of first and second gradients, and 2) discuss the dynamic behaviour of energy storage systems in two different forms of SESS and EV in MPOPF. Moreover, there are several ways to model storage devices within a time-period. Meyer-Huebner \textit{et al}.  \cite{meyer-huebner_efficient_2015} formulated three different storage device models as: (A) an inequality for the entire optimisation horizon, (B) an inequality and a variable for the entire optimisation horizon, and (C) $T$ number of equalities for each timestamp within the entire optimisation horizon. They concluded that option (C) has more efficient computational performance than options (A) and (B); therefore, we use option (C) formulation to apply the most efficient mathematical formulation with respect to the implementation of storage devices. Note that option (A) has recently been implemented in the Beltistos solver \cite{kourounis_towards_2018}. Finally, the size (number of rows and columns) of KKT systems considered to be solved in \cite{kourounis_towards_2018} is almost double the size of KKT systems we solve here.} Reference \cite{zaferanlouei_computational_2018} has developed a computationally fast solution to the MPOPF problem and applied it on a small-scale {SESS} and test system. Here, in this paper, we expand our work in \cite{zaferanlouei_computational_2018}, and propose a new solver for large-scale integration of {SESS and EV} with consideration of full AC power flow equations.\\
The main contributions of our paper could be summarised as follows: a) New input matrices are introduced in order to capture the full dynamic of a multi-period system including SESS and EV. b) The first and second analytical (hand-coded) derivatives of linear and non-linear equality constraints, inequality constraints and objective function, w.r.t. variables are extracted in the multi-period form. c) Sparsity structure of analytical derivatives is extracted in order to increase the performance in terms of memory requirements and sparsity calculations in different loops. 
d) A new re-ordering format is introduced to exploit and reveal the multi-period structure of KKT systems for both SESS and EV. e) A high-performance and memory-efficient sparse Schur-Complement algorithm is introduced in order to solve the multi-period structure of the KKT systems for both SESS and EV.
In this paper, we compare the Schur-Complement algorithm, tailored for a specific structure, with a direct sparse LU solver to shed light into a systematic comparison (same implementation platform ({MATLAB}), single-thread controlled environment, similar PC) to only compare the algorithmic complexity differences. \\
The structure of this paper is as follows: in the next section, the formulation of  MPOPF problem in the presence of SESS and EV is elaborated. The solution proposal, and the mathematical algorithms for speeding up the solution proposal, are presented in the subsequent sections. Next, the performance of the Schur-Complement algorithm is compared with that of a direct sparse LU solver for multiple numbers of storage devices and time horizons for two different cases: performance of SESS on standard mesh transmission benchmarks, and the performance of EV on radial distribution benchmarks. Finally, we summarise the paper in the last section with some concluding remarks.

\section[Formulation of Proposed Solver]{Formulation of Proposed Solver\footnote{Note that all vectors and matrices are shown with bold and non-italic notation: $\mathbf{BOLD}$}} \label{formulation}
In general, a given power system  can be represented with the following input matrices in ACOPF:  $\matr{BUS}$ matrix with $n_b \in \mathbb{N}$ number of buses. These buses are connected to each other through the total number of $n_l \in \mathbb{N}$ lines represented by the $\mathbf{BRANCH}$ matrix. The matrix consists of connecting buses, i.e.,  $\mathbf{BUS}^{\mathrm{from}}\in \mathbb{R}^{n_{l}\times 1}$ and $\mathbf{BUS}^{\mathrm{to}}\in \mathbb{R}^{n_{l}\times 1}$ and physical line parameters $\in \mathbb{R}^{n_{l}\times 1}$, such as resistance, reactance, susceptance and their apparent power capacities (MVA). The $\mathbf{GEN}$ matrix specifies the connection bus of $n_g \in \mathbb{N}$ generators  with their dispatch limits and voltage references. The $\mathbf{GENCOST}$ matrix is the associated cost functions of generators in $\mathbf{GEN}$. 
These matrices have been already defined by MATPOWER \cite{zimmerman_matpower:_2011} and are now extensively used by many researchers in the subjects of power economics and power system analysis.\\
{Here, in addition to the above input matrices, we propose new input matrices\footnote{Five non-binary and five binary matrices}, which can be seen in Appendix \ref{Appendix_A} of this paper.}   
These new input matrices are designed for the integration of large-scale SESS and EV to capture the dynamic behaviour of storage devices over the optimisation horizon. The input $\mathbf{BATT}$ represents the properties of energy storage systems and more importantly ties single-period optimal power flow equations through a positive integer parameter, called time period, $T \in \mathbb{N}$; $t$ is a time in the interval of $t \in \{1,...,T\}$. The differences between the representation of a SESS and an EV are concentrated in the binary input matrices of  $\mathbf{AVBP}$, $\mathbf{CONCH}$, $\mathbf{CONDI}$ and $\mathbf{AVBQ}$. Storage $i$  is considered stationary if the availability condition \eqref{eqn:availability} holds, otherwise it has a dynamic behaviour over time and can be considered as an EV, while charge, discharge and reactive power provision conditions could be considered as secondary or optional conditions \eqref{eqn:availabilityCH}-\eqref{eqn:availabilityQ}. 
\begin{subequations}\label{eqn:availabilityTotal}
\begin{alignat}{4}
&\Scale[0.85]{\mathbf{AVBP}\big\rvert_{i,t=1}=\mathbf{AVBP}\big\rvert_{i,t=2}=...=\mathbf{AVBP}\big\rvert_{i,t=T}=1} \label{eqn:availability}\\
&\Scale[0.85]{\mathbf{CONCH}\big\rvert_{i,t=1}=\mathbf{CONCH}\big\rvert_{i,t=2}=...=\mathbf{CONCH}\big\rvert_{i,t=T}=1}\label{eqn:availabilityCH} \\
&\Scale[0.85]{\mathbf{CONDI}\big\rvert_{i,t=1}=\mathbf{CONDI}\big\rvert_{i,t=2}=...=\mathbf{CONDI}\big\rvert_{i,t=T}=1}\label{eqn:availabilityDI}\\
&\Scale[0.85]{\mathbf{AVBQ}\big\rvert_{i,t=1}=\mathbf{AVBQ}\big\rvert_{i,t=2}=...=\mathbf{AVBQ}\big\rvert_{i,t=T}=1}\label{eqn:availabilityQ} 
\end{alignat}
\end{subequations}
In order to elaborate the mathematical formulation of the MPOPF incorporating the $\mathbf{BATT}$ matrix, we first present single-period power flow equations for a time $t$, and subsequently expand the equations to represent MPOPF.\\
Considering the above types of inputs, we have an object-oriented programming (OOP) package that constructs different types of constraints according to the input matrices and feeds the mathematical formulation of the MPOPF to the solver designed under this package.
\subsection{Single-Period Optimal Power Flow} \label{sec:single_OPF}
With the input matrices introduced above, $\matr{BUS}$, $\matr{BRANCH}$, $\matr{GEN}$ and $\matr{GENCOST}$, we show the mathematical formulation of a single-period ACOPF in this subsection. \\
Consider the vector of complex bus voltages in rectangular coordinates as illustrated  by  $\displaystyle\mathbf{\underline{V}}\in \mathbb{C}^{n_b\times 1}$, where $\mathbb{C}$ is a complex set. The voltage vector comprises complex elements as: $\mathrm{\underline{v}}_i  = \abs{\mathrm{v}_i}e^{j\mathrm{\theta}_i}$, where $\mathrm{\underline{v}}_i\in \mathbb{C}$, $\{v_i ,\theta_i \}\in  \mathbb{R}$ are the voltage magnitude and angle of the corresponding bus in polar coordinates, where $\mathbb{R}$ is a real set.  Moreover,  $\{\boldsymbol{\mathcal{V}}, \mathbf{\Theta}\} \in \mathbb{R}^{n_b \times 1}$ can be defined as vectors of real magnitude and angle of bus voltages.  In vector form, the relationship between rectangular and polar coordinates is shown as:
\begin{equation}
\label{eqn:Voltage}
\mathbf{\underline{V}}=\mathbf{diag}(\boldsymbol{\mathcal{V}})\ \exp(j\boldsymbol{\Theta})
\end{equation}
Line Connectivity matrices of $\{\mathbf{C}^{\mathrm{fr}}, \mathbf{C}^{\mathrm{to}}\} \in \mathbb{B}^{n_l \times n_b}$ can be extracted from $\mathbf{BUS}^{\mathrm{from}}$ and $\mathbf{BUS}^{\mathrm{to}}$ vectors, such that $\mathrm{c}_{ik}^{\mathrm{fr}}=1$ if bus $k$ is connected to line $i$, and otherwise $\mathrm{c}^{\mathrm{fr}}_{ik}=0$, and the same holds for $\mathbf{C}^{\mathrm{to}}$. $\{\mathbf{\underline{V}}^{\mathrm{fr}}, \ \mathbf{\underline{V}}^{\mathrm{to}}\} \in \mathbb{C}^{n_l \times 1}$ are the vectors of complex bus voltages at line terminals, including ``from" and ``to" nodes, correspondingly. These vectors can be extracted using the connectivity matrices explained above shown in Eqs. \eqref{eqn:generalProblem1} and \eqref{eqn:generalProblem2}.
\begin{equation}
\label{eqn:generalProblem1}
\mathbf{\underline{V}}^{\mathrm{fr}} = \mathbf{C}^{\mathrm{fr}}\mathbf{\underline{V}}\\
\end{equation}
\begin{equation}
\label{eqn:generalProblem2}
\mathbf{\underline{V}}^{\mathrm{to}} = \mathbf{C}^{\mathrm{to}}\mathbf{\underline{V}}\\
\end{equation}
and therefore:
\begin{align}
\label{eqn:V_line}
\mathbf{\underline{V}}^{\mathrm{Line}}=
{\begin{bmatrix}
    \mathbf{\underline{V}}^{\mathrm{fr}} \\
    \mathbf{\underline{V}}^{\mathrm{to}} \\
\end{bmatrix}\ \mkern-10mu}_{2n_l \times 1} = 
{\overbrace{\begin{bmatrix}
    \mathbf{C}^{\mathrm{fr}} \\
    \mathbf{C}^{\mathrm{to}} \\
\end{bmatrix}}^{\mathbf{C}^{\mathrm{Line}}}\ \mkern-10mu}_{2n_l \times n_b} \mathbf{\underline{V}} 
\end{align}
In order to obtain the entire network flow, the vector of complex voltages $\mathbf{\underline{V}}$ has to be determined. This can be done using the well-known Kirchhoff's current law: the sum of external current injections at a bus $\mathbf{\underline{I}}^{\mathrm{bus}} \in \mathbb{C}^{n_b\times 1}$ is equal to the sum of internal - through lines - current injections to the same bus $\mathbf{\underline{I}}^{\mathrm{bus}}= \mathbf{\underline{Y}}^{\mathrm{bus}}\mathbf{\underline{V}}$, where
$\mathbf{\underline{Y}}^{\mathrm{bus}} \in \mathbb{C}^{n_b\times n_b}$ is the bus admittance matrix. The same principle is applied to compute the complex line current using complex bus voltages of line terminals, and line admittance matrix $\mathbf{\underline{Y}}^{\mathrm{Line}} \in  \mathbb{C}^{2n_l\times n_b}$. This is shown in (\ref{eqn:generalProblem3}) 
\begin{equation}
\label{eqn:generalProblem3}
\mathbf{\underline{I}}^{\mathrm{Line}}={\begin{bmatrix}
    \mathbf{\underline{I}}^{\mathrm{fr}} \\
    \mathbf{\underline{I}}^{\mathrm{to}} \\
\end{bmatrix}\ \mkern-10mu}_{2n_l \times 1} =
{\overbrace{\begin{bmatrix}
    \mathbf{\underline{Y}}^{\mathrm{fr}} \\
    \mathbf{\underline{Y}}^{\mathrm{to}} \\
\end{bmatrix}}^{\mathbf{\underline{Y}}^{\mathrm{Line}}}\ \mkern-10mu}_{2n_l \times n_b}
\mathbf{\underline{V}} 
\end{equation}
The relation between bus admittance and line admittance matrices is defined by (\ref{eqn:generalProblem4}).
\begin{equation}
\label{eqn:generalProblem4}
\mathbf{\underline{Y}}^{\mathrm{bus}} = (\mathbf{C}^{\mathrm{fr}})^\top \mathbf{\underline{Y}}^{\mathrm{fr}}+(\mathbf{C}^{\mathrm{to}})^\top \mathbf{\underline{Y}}^{\mathrm{to}}+\mathbf{\underline{Y}}^{\mathrm{shunt}}
\end{equation}
$\{\mathbf{\underline{Y}}^{\mathrm{fr}},  \mathbf{\underline{Y}}^{\mathrm{to}}\}\in \mathbb{C}^{n_l\times n_b}$, and
$\mathbf{\underline{Y}}^{\mathrm{shunt}} \in \mathbb{C}^{n_b\times n_b}$ is the matrix of shunt admittance. 
Finally, the external complex power injections into a bus $i$ can be computed as $\mathrm{\underline{s}}^{\mathrm{bus}}_i = \underline{\mathrm{v}}_i(\underline{\mathrm{i}}^{\mathrm{bus}}_i)^*$, whereas the  complex power flow over a line at the terminal $k$ can be calculated by  $\underline{\mathrm{s}}^{\mathrm{Line}}_k = (\mathbf{C}_k^{\mathrm{Line}}\mathbf{\underline{V}})(\underline{\mathrm{i}}^{\mathrm{Line}}_k)^*$, where $\{\underline{\mathrm{s}}^{\mathrm{bus}}_i ,  \underline{\mathrm{i}}^{\mathrm{bus}}_i , \underline{\mathrm{s}}^{\mathrm{Line}}_k, \underline{\mathrm{i}}^{\mathrm{Line}}_i  \} \in \mathbb{C}$ and $\mathbf{C}_k^{\mathrm{Line}} \in \mathbb{B}^{1 \times n_b}$  is the $k$\textsuperscript{th} element of $\mathbf{C}^{\mathrm{Line}}$ matrix. In summary, power injections into a bus and into a line can be extended in the form of vectors using (\ref{eqn:generalProblem5}).   
\begin{align}
\label{eqn:generalProblem5}
&\mathbf{\underline{S}}^{\mathrm{bus}} = \mathbf{diag}(\mathbf{\underline{V}})(\mathbf{\underline{I}}^{\mathrm{bus}})^* \ \ \in \ \ \mathbb{C}^{n_b \times 1}
\\
&\mathbf{\underline{S}}^{\mathrm{Line}} = \mathbf{diag}(\mathbf{\underline{V}}^{\mathrm{Line}})(\mathbf{\underline{I}}^{\mathrm{Line}})^* \ \ \in \ \ \mathbb{C}^{2n_l \times 1}
\end{align}
Here we define generator connectivity matrix $\mathbf{C}^\mathrm{g} \in \mathbb{B}^{n_b \times n_g}$ which is a binary matrix of 0 and 1. $\mathbf{C}^\mathrm{g}_{ik} =1$ if generator $i$ is connected to the bus $k$ and otherwise 0. $\underline{\mathbf{S}}^\mathrm{g} \in \mathbb{C}^{n_g\times 1}$ and $\underline{\mathbf{S}}^\mathrm{d} \in \mathbb{C}^{n_b\times 1}$ are the complex vectors of generation units and loads. The first set of non-linear equality constraints for the ACOPF problem, can be defined as: 
\begin{equation}
\label{eqn:Gs}
\mathbf{\widetilde{g}}(\mathbf{x})=\mathbf{\underline{S}}^{\mathrm{bus}}+\underline{\mathbf{S}}^\mathrm{d}-\mathbf{C}^\mathrm{g}\underline{\mathbf{S}}^{\mathrm{g}}=0 
\end{equation} 
 In this article, $\sim$ is the sign for a non-linear equation. In the literature on power systems, \eqref{eqn:Gs} has been divided into two parts, power balance of active and reactive power as indicated in Eqs. \eqref{active_Power_flow} and \eqref{Reactive_Power_flow}, respectively. Therefore, non-linear equality constraints take the form of $\mathbf{\widetilde{g}}(x) \in \mathbb{R}^{n_{gn} \times 1}$ and $n_{gn}= 2n_b$. 
\begin{subequations}
\begin{alignat}{2}
\mathbf{C}^\mathrm{g} \boldsymbol{\mathcal{P}}^{\mathrm{g}}-\boldsymbol{\mathcal{P}}^{\mathrm{d}}=\Re{[\mathbf{C}^\mathrm{g}\underline{\mathbf{S}}^{\mathrm{g}}]} - \Re{[\underline{\mathbf{S}}^\mathrm{d}]}=\Re{[\mathbf{\underline{S}}^{\mathrm{bus}}]} \label{active_Power_flow}\\
\mathbf{C}^\mathrm{g} \boldsymbol{\mathcal{Q}}^{\mathrm{g}}-\boldsymbol{\mathcal{Q}}^{\mathrm{d}}=\Im{[\mathbf{C}^\mathrm{g}\underline{\mathbf{S}}^{\mathrm{g}}]} - \Im{[\underline{\mathbf{S}}^\mathrm{d}]}=\Im{[\mathbf{\underline{S}}^{\mathrm{bus}}]}  \label{Reactive_Power_flow} 
\end{alignat}
\end{subequations} 
where $\{\boldsymbol{\mathcal{P}}^{\mathrm{g}}, \boldsymbol{\mathcal{Q}}^{\mathrm{g}}\} \in \mathbb{R}^{n_g \times 1}$ are the vectors of active and reactive power generations, $\{\boldsymbol{\mathcal{P}}^{\mathrm{d}}, \boldsymbol{\mathcal{Q}}^{\mathrm{d}}\} \in \mathbb{R}^{n_b \times 1}$ are the vectors of active and reactive power consumption. If we take $(\abs{\mathbf{\underline{S}}^{\mathrm{Line}}_\mathrm{max}})^{2} \in \mathbb{R}^{2n_l \times 1}$ as the squared vector of apparent power flow limits, then apparent line power flow constraint can be defined as: 
\begin{align}
\begin{split}
\label{eqn:h_t}
     \mathbf{\widetilde{h}}(\mathbf{x})=
     \big[(\mathbf{\underline{S}}^{\mathrm{Line}})^*\mathbf{\underline{S}}^{\mathrm{Line}}-(\abs{\mathbf{\underline{S}}^{\mathrm{Line}}_\mathrm{max}})^{2}\big]
 \leq 0  \in \mathbb{R}^{n_{hn} \times 1}
 \end{split}
\end{align}
where $n_{hn}=2n_l$. Another type of constraint is the linear equality constraint $\mathbf{\overline{g}}(\mathbf{x}) =0$ where $ - $ is the sign for a linear vector. Voltage angle is kept equal to zero $\theta^{\mathrm{slack}}=0 \ \in \mathbb{R}$ which is a linear equality constraint for the slack bus. The last set of constraints is that of linear inequality, related to upper and lower bounds variables called box constraints \cite{kourounis_towards_2018} with vectors $\boldsymbol{\Theta}^\mathrm{min}\leq
\boldsymbol{\Theta} \leq \boldsymbol{\Theta}^{\mathrm{max}} \in \mathbb{R}^{(2n_b-1) \times 1}$,\  $\boldsymbol{\mathcal{V}}^{\mathrm{min}}\leq\boldsymbol{\mathcal{V}}\leq\boldsymbol{\mathcal{V}}^{\mathrm{max}}  \in \mathbb{R}^{2n_b \times 1},\ \{(\boldsymbol{\mathcal{P}}^{\mathrm{g}})^{\mathrm{min}}\leq\boldsymbol{\mathcal{P}}^{\mathrm{g}}\leq(\boldsymbol{\mathcal{P}}^{\mathrm{g}})^{\mathrm{max}},\ (\boldsymbol{\mathcal{Q}}^{\mathrm{g}})^{\mathrm{min}}\leq\boldsymbol{\mathcal{Q}}^{\mathrm{g}}\leq(\boldsymbol{\mathcal{Q}}^{\mathrm{g}})^{\mathrm{max}}\}\in \mathbb{R}^{2n_g \times 1}$. All the box constraints could be put together and represented with one vector $\mathbf{\overline{h}}(\mathbf{x}) \in \mathbb{R}^{n_{hl} \times 1}$ where $n_{hl}=(4n_b-1)+4n_g$.  \\
Thus, if the objective function $f(\mathbf{x})$ is an arbitrary linear or non-linear function related to the cost of power generation, the general optimisation framework could be:
\begin{equation}
\begin{multlined}
\label{eqn:generalProblem}
\min_{\mathbf{x}} f(\mathbf{x})\\
\textrm{s.t. }  \mathbf{g}(\mathbf{x})= \begin{bmatrix}
   \mathbf{\widetilde{g}}(\mathbf{x})\\
   \mathbf{\overline{g}}(\mathbf{x})
\end{bmatrix}=0 \qquad \in   \mathbb{R}^{n_{gx} \times 1}\\
\mathbf{h}(\mathbf{x})=\begin{bmatrix}
   \mathbf{\widetilde{h}}(\mathbf{x})\\
   \mathbf{\overline{h}}(\mathbf{x})
\end{bmatrix}\leq 0 \qquad \in   \mathbb{R}^{n_{hx} \times 1}
\end{multlined}
\end{equation}
where $n_{gx} = n_{gn}+n_{gl}$\footnote{{ $n_{gn}$ is the number of nonlinear equality constraints and $n_{gl}$ is the number of linear equality constraints}}, $n_{hx} = n_{hn}+n_{hl}$\footnote{{$n_{hn}$ is the number of nonlinear inequality constraints and $n_{hl}$ is the number of linear inequality constraints}}, and 
\begin{equation}
\begin{multlined}
\label{eqn:x}
\mathbf{x}= \big[\boldsymbol{\Theta}\  \boldsymbol{\mathcal{V}} \  \boldsymbol{\mathcal{P}}^{\mathrm{g}} \ \boldsymbol{\mathcal{Q}}^{\mathrm{g}} \big]^\top
 \in \mathbb{R}^{n_{x}\times 1}
\end{multlined}
\end{equation}
and $n_{x}= 2n_b+2n_g$. With the discussed input matrices and use of OOP\footnote{{object oriented programming}}, we construct {an instance of an object called $\mathbf{MP}$} which incorporates all the introduced network variables and constraints in the formulated problem. In the next subsection, we come up with a new input matrix, and make a generalised case for storage devices.  
\subsection{Multi-Period Optimal Power Flow}\label{Sec:multiAC}
``$\mathbf{BATT}$" is the main introduced matrix here in this section, to optimally operate an electricity network over a time horizon of $T$ which defines the number of time-steps in the optimisation. Therefore, multiple single-period optimal power flow problems are coupled together over a given  time horizon to determine the optimal combined operation schedule for energy storage systems. The coupling constraints over time are introduced by linear equality constraints of storage devices in  $e_{i,t}- e_{i,t-1}-\psi_{i}^\mathrm{ch}p_{i,t}^\mathrm{ch}\Delta t+\frac {p_{i,t}^\mathrm{dch}\Delta t}{\psi_{i}^\mathrm{dch}} =0$, where $soc_{i,t}= \frac {e_{i}}{e_{i}^{\mathrm{max}}}$ and $\{soc\mathrm{,} \ p^\mathrm{ch} $,  $ \ p^\mathrm{ch}, q^{s}, e, \psi \} \in \mathbb{R}$ representing storage device variables connected to bus $i$ at time $t$. The presence of: 1) cost minimisation objective function, and 2) the efficiency of charge and discharge helps to avoid the simultaneous charging and discharging. However, charging and discharging might still occur at the same time in some cases, although this is outside the scope of this paper.\\
 To extend the vector notation, for $n_y$ storage devices, we have: 
\begin{equation}\label{storages}
 \mathbf{\overline{g}}^s(\boldsymbol{\tau}_t)=\mathbf{E}_t- \mathbf{E}_{t-1}-\mathbf{\Psi}^\mathrm{ch}\boldsymbol{\mathcal{P}}_t^\mathrm{ch}\Delta t+\frac {\boldsymbol{\mathcal{P}}_t^\mathrm{dch}\Delta t}{\mathbf{\Psi}^\mathrm{dch}} =0    
\end{equation}
 where $\boldsymbol{\mathcal{SOC}}_t=\frac{\mathbf{E}_t}{\mathbf{E}^{max}}$ and $\{\boldsymbol{\mathcal{SOC}}_t$, $\boldsymbol{\mathcal{P}}_t^\mathrm{ch}, \boldsymbol{\mathcal{P}}_t^\mathrm{dch}, \boldsymbol{\mathcal{Q}}_t^\mathrm{s}\} \in \ \mathbb{R}^{n_y \times 1}$. Moreover, $\boldsymbol{\tau}_t =\{\mathbf{x}_{t-1}, \mathbf{x}_t\}$ and $\mathcal{T}=\{\boldsymbol{\tau}_1,\boldsymbol{\tau}_2,\dots,\boldsymbol{\tau}_T\}= \{\{\mathbf{x}_1\},\{\mathbf{x}_1, \mathbf{x}_2\},\dots,\{\mathbf{x}_{T-1}, \mathbf{x}_{T}\}\}=\{\mathbf{x}_1,\mathbf{x}_2,\dots,\mathbf{x}_T\}$, thus  $\overline{\mathbf{G}}^s(\mathcal{T})$\footnote{Note that $\boldsymbol{\tau}_t$ and $\mathcal{T}$ are representations of two sets such that $\boldsymbol{\tau}_t\in \mathcal{T}$}$= \overline{\mathbf{G}}^s(\mathbf{X})$. Thus, four vectors above representing $n_y$  are added to the variables shown in subsection \ref{sec:single_OPF}. Note that initial state of charge of each storage $i$ at time $t$ is defined as $e_{i,t-1}= e_{i}^{max}\mathbf{SOCi}_{i,t}$ where $\mathbf{SOCi}$ is the input matrix introduced in \ref{formulation} such that an initial value of $\mathbf{SOCi}_{i,t}$ is allocated if one of the arrival conditions are satisfied: 1) $\mathbf{AVBP}_{i,t=1}=1$. 2) $\mathbf{AVBP}_{i,t-1}=0$ and $\mathbf{AVBP}_{i,t}=1$.\\
 The procedure to initiate a MPOPF problem is as follows: First, we read information on storage devices such as bus location, charge and discharge efficiencies, capacities and initial state of charge from $\mathbf{BATT}$ matrix, and save them in the $\mathbf{MP}$ instance, previously initiated by other input matrices as introduced previously. Then, we construct  the above-mentioned storage variables and their linear equality constraints in the current $\mathbf{MP}$ instance.\\
Furthermore, four new three-dimensional connectivity matrices of $\{ \mathbf{C}^\mathrm{ch}, \mathbf{C}^\mathrm{dch}, \mathbf{C}^\mathrm{s} \} \in \mathbb{B}^{n_b \times n_y \times T}$ and $\mathbf{C}^\mathrm{g} \in \mathbb{B}^{n_b \times n_g \times T}$ are constructed. The first two are built based on: 1) bus connectivity matrix $\mathbf{BATT}$, 2) active power provision $\mathbf{AVBP}$, and 3) whether charge or discharge options are available from $\mathbf{CONCH}$ and $\mathbf{CONDI}$ matrices, such that $\mathbf{C}^\mathrm{ch}_{ikt} =1$ if both the following conditions are satisfied: (a) battery $i$ is connected to bus $k$ at time $t$, and (b) charge of battery $i$ at time $t$ is activated, otherwise 0. $\mathbf{C}^\mathrm{dch}_{jkt} =1$  if both the following conditions are satisfied (a) battery $j$ is connected to bus $k$ at time $t$, and (b) discharge of battery $j$ at time $t$ is activated, otherwise 0.\\
Connectivity matrix for reactive power provision $\mathbf{C}^\mathrm{s}$ is constructed by: 1) bus connectivity matrix $\mathbf{BATT}$, 2) reactive power provision $\mathbf{AVBQ}$ such that $\mathbf{C}^\mathrm{s}_{mkt} =1$ if both the following conditions are satisfied: (a) battery $m$ is connected to bus $k$ at time $t$, and (b) reactive power injection or absorption of battery $m$ at time $t$ is activated, otherwise 0. Similarly, a generator connectivity matrix $\mathbf{C}^\mathrm{g}$ is built with: 1) bus connectivity matrix of generators $\mathbf{GEN}$, 2)  active and reactive power provision $\mathbf{AVG}$ such that $\mathbf{C}^\mathrm{g}_{nkt} =1$ if the generator $n$ is connected to bus $k$ is running at time $t$, otherwise 0. Thus,
Eqs. \eqref{active_Power_flow} and \eqref{Reactive_Power_flow} can be re-formulated now to Eqs. \eqref{active_Power_flow2} and \eqref{Reactive_Power_flow2} respectively.
\begin{subequations}
\begin{alignat}{2}
&\mathbf{C}_t^\mathrm{g} \boldsymbol{\mathcal{P}}_t^{\mathrm{g}}-\boldsymbol{\mathcal{P}}_t^{\mathrm{d}}-\mathbf{C}_t^\mathrm{ch} \boldsymbol{\mathcal{P}}_t^{\mathrm{ch}}+\mathbf{C}_t^\mathrm{dch} \boldsymbol{\mathcal{P}}_t^{\mathrm{dch}}-\Re{[\mathbf{\underline{S}}_t^{\mathrm{bus}}]}=0 \label{active_Power_flow2}\\
&\mathbf{C}_t^\mathrm{g} \boldsymbol{\mathcal{Q}}_t^{\mathrm{g}}-\boldsymbol{\mathcal{Q}}_t^{\mathrm{d}}+\mathbf{C}_t^\mathrm{s} \boldsymbol{\mathcal{Q}}_t^{\mathrm{s}}-\Im{[\mathbf{\underline{S}}_t^{\mathrm{bus}}]}=0 \label{Reactive_Power_flow2}
\end{alignat}
\end{subequations} 
where, $\boldsymbol{\mathcal{P}}_t^{\mathrm{d}}=\mathbf{PD}_t$. and $\boldsymbol{\mathcal{Q}}_t^{\mathrm{d}}= \mathbf{QD}_t$ In summary, the total number of variables for each time-step becomes:
\begin{equation}
\begin{multlined}
\label{eqn:x_t}
{\mathbf{x}_{t}= \big[\boldsymbol{\Theta}_{t} \  \boldsymbol{\mathcal{V}}_{t} \ \boldsymbol{\mathcal{P}}^{\mathrm{g}}_{t} \ \boldsymbol{\mathcal{Q}}^{\mathrm{g}}_{t} \   \boldsymbol{\mathcal{SOC}}_{t}\  \boldsymbol{\mathcal{P}}_{t}^\mathrm{ch} \ \boldsymbol{\mathcal{P}}_{t}^\mathrm{dch} \ \boldsymbol{\mathcal{Q}}_{t}^\mathrm{s}\big]^\top
 \ \mkern-10mu}_{1 \times N_{x_t}}
\end{multlined}
\end{equation}
$N_{x_t} = n_{x}+4n_y$. Subscript $t$ stands for a specific time-step in this paper. In addition to box constraints defined in \ref{sec:single_OPF}, we define more box constraints corresponding to the new defined storage variables.  $\{\mathbf{SOCMi}_t\leq \boldsymbol{\mathcal{SOC}}_t \leq \boldsymbol{\mathcal{SOC}}^\mathrm{max} $, $(\boldsymbol{\mathcal{P}}^\mathrm{ch})^\mathrm{min} \leq \boldsymbol{\mathcal{P}}_t^{ch} \leq(\boldsymbol{\mathcal{P}}^\mathrm{ch})^\mathrm{max}$, $(\boldsymbol{\mathcal{P}}^\mathrm{dch})^\mathrm{min} \leq \boldsymbol{\mathcal{P}}_t^\mathrm{dch}\leq (\boldsymbol{\mathcal{P}}^{\mathrm{dch}})^\mathrm{max}$ and $(\boldsymbol{\mathcal{Q}}^\mathrm{s})^\mathrm{min} \leq
\boldsymbol{\mathcal{Q}}_t^\mathrm{s} \leq(\boldsymbol{\mathcal{Q}}^\mathrm{s})^\mathrm{max}\} \in \mathbb{R}^{n_y \times 1}$
The vector of total variables in the MPOPF problem $\mathbf{X} \in \ \mathbb{R}^{N_{x} \times 1}$  where $N_{x} = TN_{x_t}$, is shown in (\ref{X}):
\begin{equation}
\begin{multlined}
\label{X}
\mathbf{X}= \big[\mathbf{x}_{1} \quad  \mathbf{x}_{2}\quad ...\quad \mathbf{x}_{t}\quad ...\quad  \mathbf{x}_{T}\big]^\top
\end{multlined}
\end{equation}
Finally, a general MPOPF problem can be formulated as:
\begin{subequations}
\begin{alignat}{4}
&\min_{\mathbf{X}} F(\mathbf{X})\\ \label{multiperiodG.a}
\textrm{s.t. } G(\mathbf{X})&= \begin{bmatrix}
   \widetilde{G}(\mathbf{X}) \
   \overline{G}(\mathbf{X}) \
   \overline{G}^s(\mathbf{X})
\end{bmatrix}^\top&&=0  \in   \mathbb{R}^{N_{g} \times 1}\\\label{multiperiodH.b}
H(\mathbf{X})&=\begin{bmatrix}
   \widetilde{H}(\mathbf{X}) \quad
   \overline{H}(\mathbf{X})
\end{bmatrix}^\top&&\leq 0  \in   \mathbb{R}^{N_{h} \times 1}
\end{alignat}
\end{subequations}
where $F(\mathbf{X})=f_{t=1}(\mathbf{x}_{1})+f_{t=2}(\mathbf{x}_{2})+\dots+f_{t=T}(\mathbf{x}_{T})$, $\widetilde{\mathbf{G}}(\mathbf{X})\in \mathbb{R}^{N_{gn} \times 1}$, $\overline{\mathbf{G}}(\mathbf{X})\in \mathbb{R}^{N_{gl} \times 1}$, $\overline{\mathbf{G}}^s(\mathbf{X})\in \mathbb{R}^{N_{gs} \times 1}$, $\widetilde{\mathbf{H}}(\mathbf{X})\in \mathbb{R}^{N_{hn} \times 1}$ and $\overline{\mathbf{H}}(\mathbf{X})\in \mathbb{R}^{N_{hl} \times 1}$ are as shown:
\begin{subequations}\label{19}
\begin{alignat}{4}
   \widetilde{\mathbf{G}}(\mathbf{X})=&
\begin{bmatrix}
    \widetilde{\mathbf{g}}(\mathbf{x}_{1}) \
    \widetilde{\mathbf{g}}(\mathbf{x}_{2}) \
    \dots \
    \widetilde{\mathbf{g}}(\mathbf{x}_{T}) 
\end{bmatrix}^\top \label{19.a}\\
   \overline{\mathbf{G}}(\mathbf{X})=&
\begin{bmatrix}
    \overline{\mathbf{g}}(\mathbf{x}_{1}) \
    \overline{\mathbf{g}}(\mathbf{x}_{2}) \
    \dots \
    \overline{\mathbf{g}}(\mathbf{x}_{T}) \
\end{bmatrix}^\top \label{19.b}\\
  \overline{\mathbf{G}}^s(\mathbf{X})=&
\begin{bmatrix}
    \overline{\mathbf{g}}^s(\boldsymbol{\tau}_1) \
    \overline{\mathbf{g}}^s(\boldsymbol{\tau}_2) \
    \dots \
    \overline{\mathbf{g}}^s(\boldsymbol{\tau}_T) 
\end{bmatrix}^\top \label{19.c} \\
   \widetilde{\mathbf{H}}(\mathbf{X})=&
\begin{bmatrix}
    \widetilde{\mathbf{h}}(\mathbf{x}_{1}) \
    \widetilde{\mathbf{h}}(\mathbf{x}_{2}) \
    \dots \
    \widetilde{\mathbf{h}}(\mathbf{x}_{T}) 
\end{bmatrix}^\top \label{19.d}\\
  \overline{\mathbf{H}}(\mathbf{X})=& 
\begin{bmatrix}
    \overline{\mathbf{h}}(\mathbf{x}_{1}) \
    \overline{\mathbf{h}}(\mathbf{x}_{2}) \
    \dots \
    \overline{\mathbf{h}}(\mathbf{x}_{T}) 
\end{bmatrix}^\top \label{19.e}
\end{alignat}
\end{subequations}
where $N_{g} =N_{gn}+N_{gl}+N_{gs}$,\  $N_{gn} = Tn_{gn},\ N_{gl}= n_{{gl}_{t=1}}+n_{{gl}_{t=2}}+...+n_{{gl}_{t=T}},\ N_{gs}=Tn_y,\ N_{h} =N_{hn}+N_{hl},\ N_{hn}= Tn_{hn}, \ N_{hl}=n_{{hl}_{t=1}}+n_{{hl}_{t=2}}+...+n_{{hl}_{t=T}}+T(8n_{y})$, $\boldsymbol{\tau}_1 =\{\mathbf{x}_1\}$. Furthermore, $\widetilde{\mathbf{g}}(\mathbf{x}_t)$ contains the two new defined constraints of \eqref{active_Power_flow2} and \eqref{Reactive_Power_flow2}.
 $\overline{\mathbf{g}}(\mathbf{x}_t)$ includes \eqref{eqn:equalitylineargrida}-\eqref{eqn:equalitylineargena} plus any other upper and lower bounds of variable $\mathbf{x}_t$ such that $x_t^\mathrm{min}=x_t^\mathrm{max}$, which can be user defined, and as such can be removed from the list of box constraints in \eqref{19.e} and is introduced here as a new linear equality \eqref{eqn:equalitylineargridd}.
\begin{subequations} \label{eqn:equalitylineargrid}
\begin{alignat}{6}
&\theta_t^{\mathrm{slack}}=0&&&\ \ \ \ &\label{eqn:equalitylineargrida}\\
 &p_{i,t}^\mathrm{ch}=0,\ \ \text{if } \quad \{\mathbf{AVBP}_{i,t}	\lor \mathbf{CONCH}_{i,t} \}=0\\
 & p_{i,t}^\mathrm{dch}=0,\ \text{if }\quad \{\mathbf{AVBP}_{i,t}	\lor \mathbf{CONDI}_{i,t} \}=0\\
&q_{i,t}^\mathrm{s}=0,\ \ \text{if }\quad \{\mathbf{AVBP}_{i,t}	\lor \mathbf{AVBQ}_{i,t}\}=0 \label{eqn:equalitylineargridc}\\
&p_{i,t}^\mathrm{g}=0,\ \ \text{if }\quad \mathbf{AVG}_{i,t}=0 \label{eqn:equalitylineargena}\\
&q_{i,t}^\mathrm{g}=0,\ \ \text{if }\quad \mathbf{AVG}_{i,t}=0 \label{eqn:equalitylineargenb}\\
 & x_t=x_t^\mathrm{min}= x_t^\mathrm{max} \  \ \text{if }  x_t^\mathrm{min}= x_t^\mathrm{max} \label{eqn:equalitylineargridd}
 \end{alignat}
 \end{subequations}
 $\overline{\mathbf{g}}^\mathrm{s}(\mathbf{x}_t)$ is defined from \eqref{storages}, $\widetilde{\mathbf{h}}(\mathbf{x}_{t})$ is the non-linear inequality constraints for time $t$ and is similar to \eqref{eqn:h_t}. Finally $\overline{\mathbf{h}}(\mathbf{x}_{t})$ is the set of box constraints of all variables except the slack bus. Our proposed formulation through the instance of $\mathbf{MP}$ is fed to the solver that is introduced in the next section.\\
Note that the number of linear equality constraints of $n_{{gl}_t}$ and the number of linear inequality constraints of $n_{{hl}_t}$, in each time $t = \{1,...,T\}$, are dependent on the availability matrices $\mathbf{AVBP}$, $\mathbf{CONCH}$, $\mathbf{CONDI}$, $\mathbf{AVBQ}$ and $\mathbf{AVG}$ introduced in \ref{formulation}. These numbers play an important role in the Jacobian structure of solution proposal and the follow-up re-ordering section, which will be explained in \ref{reordering}. In brief, they are constant numbers over time $t = \{1,...,T\}$ in the optimisation, if all the storage devices and generators have similar input availability matrices over time $t$, as shown in \eqref{eqn:availabilityTotal} (if all storage devices are {SESS}). 
\section{Solution proposal}
\subsection{Primal-Dual Interior Point}
The problem formulated in Section \ref{Sec:multiAC} can be solved using primal-dual interior method \cite{wang_computational_2007}. 
This can be implemented by converting the inequality equations to equality in \eqref{eqn:generalProblem} using slack variable of $z_i \in \mathbb{R}$, where $i$ denotes the number of inequality equation $\{i|i \in \mathbb{N}, 1\leq i \leq N_h\}$ and applying barrier function for slack variables:
\begin{subequations}
\begin{alignat}{4}
\min_{\mathbf{X}} &\bigg[F(\mathbf{X})-\gamma\sum_{i=1}^{N_h}{ln(z_{i})}\bigg] \label{eqn:generalbarrier}\\
&\textrm{s.t. }  \mathbf{G}(\mathbf{X})=0,\\
&\mathbf{H}(\mathbf{X})+\mathbf{Z}=0\\
&\mathbf{Z}\geq 0 \label{z}
\end{alignat}
\end{subequations}
where $\mathbf{Z} \in \mathbb{R}^{N_g\times 1}$ is the vector of slack variables and $\gamma$ is the perturbation parameter which reduces to zero when the problem approaches to optimal point. Lagrangian function of the  sub-problems \eqref{eqn:generalbarrier}-\eqref{z} becomes:

\begin{equation}
\begin{multlined}
\label{eqn:lagrangian}
\boldsymbol{\mathcal{L}}^{\gamma}(\mathbf{X},\mathbf{Z},\boldsymbol{\lambda},\boldsymbol{\mu})=f(\mathbf{X})+\boldsymbol{\lambda}^\top \mathbf{G}(\mathbf{X})\\
+\boldsymbol{\mu}^\top(\mathbf{H}(\mathbf{X})+\mathbf{Z})-\gamma \sum_{i=1}^{N_g}{ln(z_{i})}
\end{multlined}
\end{equation}
where $\boldsymbol{\lambda} \in \mathbb{R}^{N_g\times 1}, \boldsymbol{\mu} \in \mathbb{R}^{N_h\times 1}$ are the vectors of Lagrange multipliers for equality and inequality constraints. To write Karush-Kuhn-Tucker (KKT) conditions, partial differentials of \eqref{eqn:lagrangian} can be extracted with respect to the all variables:
\begin{subequations}
\begin{alignat}{4}
&\boldsymbol{\mathcal{L}}_{\mathbf{X}}^{\gamma}(\mathbf{X},\mathbf{Z},\boldsymbol{\lambda},\boldsymbol{\mu})=f_{\mathbf{X}}+\boldsymbol{\lambda}^\top \mathbf{G}_{\mathbf{X}}+\boldsymbol{\mu}^\top \mathbf{H}_{\mathbf{X}}=0 \label{kkt1}\\
&\boldsymbol{\mathcal{L}}_{\mathbf{Z}}^{\gamma}(\mathbf{X},\mathbf{Z},\boldsymbol{\lambda},\boldsymbol{\mu})=\boldsymbol{\mu}^\top - \gamma \mathbf{e}^\top\mathbf{diag}(\mathbf{Z})^{-1}=0\label{2}\\
&\boldsymbol{\mathcal{L}}_{\boldsymbol{\lambda}}^{\gamma}(\mathbf{X},\mathbf{Z},\boldsymbol{\lambda},\boldsymbol{\mu})=\mathbf{G}^\top (\mathbf{X})=0\label{kkt3}\\
&\boldsymbol{\mathcal{L}}_{\boldsymbol{\mu}}^{\gamma}(\mathbf{X},\mathbf{Z},\boldsymbol{\lambda},\boldsymbol{\mu})=\mathbf{H}^\top(\mathbf{X})+\mathbf{Z}^\top=0\label{kkt4}
\end{alignat}
\end{subequations}
where $f_\mathbf{X} \in \mathbb{R}^{N_x\times 1}$, $\mathbf{G}_{\mathbf{X}}  \in \mathbb{R}^{N_g\times N_x}$ and $\mathbf{H}_\mathbf{X} \in \mathbb{R}^{N_h\times N_x}$ are partial differentials of objective function, equality constraints and inequality constraints with respect to $\mathbf{X}$, and $\mathbf{e} \in \{\mathbb{1}\}^{N_h\times 1}$. Eqs. \eqref{kkt1}-\eqref{kkt4} can be written as \eqref{kkm1}  in a matrix form.
\begin{subequations}
\begin{alignat}{4}
\boldsymbol{\Omega}(\mathbf{X},\mathbf{Z},\boldsymbol{\lambda},\boldsymbol{\mu})
=&\begin{bmatrix}
    f_{\mathbf{X}}+\boldsymbol{\lambda}^\top \mathbf{G}_{\mathbf{X}}+\boldsymbol{\mu}^\top \mathbf{H}_{\mathbf{X}} \\
   \mathbf{diag}(\mathbf{Z})\boldsymbol{\mu}^\top - \gamma \mathbf{e}^\top\\
    \mathbf{G}^\top (\mathbf{X})\\
    \mathbf{H}^\top(\mathbf{X})+\mathbf{Z}^\top
\end{bmatrix}=0 \label{kkm1}\\
& \quad \mathbf{Z} > 0\\
& \quad \boldsymbol{\mu} > 0
\end{alignat}
\end{subequations}
 We applied  Newton-Raphson method \cite{zimmerman2010ac} to solve sets of equations in \eqref{kkm1}, and hence we have:
\begin{equation}\label{newton1}
    \Scale[0.88]{[\boldsymbol{\Omega}_\mathbf{X} \ \boldsymbol{\Omega}_\mathbf{Z} \ \boldsymbol{\Omega}_{\boldsymbol{\lambda}} \ \boldsymbol{\Omega}_{\boldsymbol{\mu}}]^k{[\Delta \mathbf{X} \ \Delta \mathbf{Z} \ \Delta \boldsymbol{\lambda} \ \Delta \boldsymbol{\mu}]^\top}^k=-\boldsymbol{\Omega}(\mathbf{X},\mathbf{Z},\boldsymbol{\lambda},\boldsymbol{\mu})^k}
\end{equation}
where $k$ is the iteration number in each step. In order to use Newton-Raphson's method to solve equation, partial differential equations of $\boldsymbol{\Omega}_X$, $\boldsymbol{\Omega}_Z$, $\boldsymbol{\Omega}_{\boldsymbol{\lambda}}$ and $\boldsymbol{\Omega}_{\boldsymbol{\mu}}$ must be calculated as shown in \eqref{Hessian of Newton}.
\begin{equation}\label{Hessian of Newton}
\Scale[0.89]{
    \begin{bmatrix}
       \boldsymbol{\mathcal{L}}_{\mathbf{X}\mathbf{X}}^{\gamma} & 0 & \mathbf{G}_{\mathbf{X}}^\top & \mathbf{H}_{\mathbf{X}}^\top \\
       0 & \mathbf{diag}(\boldsymbol{\mu}) &0 & \mathbf{diag}(\mathbf{Z})\\
       \mathbf{G}_{\mathbf{X}} & 0 &0 &0\\
       \mathbf{H}_{\mathbf{X}} & I & 0&0
    \end{bmatrix}^k} \Scale[0.9]{\begin{bmatrix}
       \Delta \mathbf{X}\\
       \Delta \mathbf{Z}\\
       \Delta \boldsymbol{\lambda}\\
       \Delta \boldsymbol{\mu}
    \end{bmatrix}^k}\Scale[0.89]{=-\begin{bmatrix}
       {\boldsymbol{\mathcal{L}}_{\mathbf{X}}^{\gamma}}^\top\\
       {\boldsymbol{\mathcal{L}}_{\mathbf{Z}}^{\gamma}}^\top\\
       {\boldsymbol{\mathcal{L}}_{\boldsymbol{\lambda}}^{\gamma}}^\top\\
      {\boldsymbol{\mathcal{L}}_{\boldsymbol{\mu}}^{\gamma}}^\top
    \end{bmatrix}^k}
\end{equation}
where $\boldsymbol{\mathcal{L}}_{\mathbf{X}\mathbf{X}}^{\gamma}(\mathbf{X},\mathbf{Z},\boldsymbol{\lambda},\boldsymbol{\mu})=f_{\mathbf{X}\mathbf{X}}+ \mathbf{G}_{\mathbf{X}\mathbf{X}}(\boldsymbol{\lambda})+\mathbf{H}_{\mathbf{X}\mathbf{X}}(\boldsymbol{\mu})$. Looking at the structure of the coefficient matrix in \eqref{Hessian of Newton}, we are able to eliminate two of the rows. In the second row of \eqref{Hessian of Newton}, we have the equation as $\mathbf{diag}(\boldsymbol{\mu})\Delta \mathbf{X}+\mathbf{diag}(\mathbf{Z})\Delta\boldsymbol{\mu}=-\mathbf{diag}(\boldsymbol{\mu})\mathbf{Z}+\gamma \mathbf{e}$, where $\Delta \boldsymbol{\mu}$ can re-written as a function of $\Delta \mathbf{Z}$ as in \eqref{deltamu}. 
\begin{equation}\label{deltamu}
    \Delta \boldsymbol{\mu} = -\boldsymbol{\mu}+\mathbf{diag}(\mathbf{Z})^{-1}(\gamma \mathbf{e} -\mathbf{diag}(\boldsymbol{\mu})\Delta \mathbf{Z})
\end{equation}
The same holds for the fourth row of \eqref{Hessian of Newton} which can be replaced by $\Delta \mathbf{Z}$ as a function of $\Delta \mathbf{X}$ as in \eqref{deltaZ}. 
\begin{equation} \label{deltaZ}
    \Delta \mathbf{Z} = -\mathbf{H}(\mathbf{X})-\mathbf{Z}-\mathbf{H}_\mathbf{X} \Delta \mathbf{X}
\end{equation}
Thus, two rows of matrix \eqref{Hessian of Newton} are taken out of the sets of equation and can be calculated by substituting $\Delta X$ in \eqref{deltaZ} and then $\Delta Z$ in \eqref{deltamu}. After the eliminations of above mentioned rows and several stages of simplifications, \eqref{Hessian of Newton} can finally be written as :
\begin{align}
\begin{split}
\label{eqn:kktn}
{\begin{bmatrix}
    \mathbf{M}&  \mathbf{G}_{\mathbf{X}}^\top\\
    \mathbf{G}_{\mathbf{X}} & 0\\
\end{bmatrix}}^k
{\begin{bmatrix}
    \Delta \mathbf{X}\\
   \Delta\boldsymbol{\lambda}
\end{bmatrix}}^k=
{\begin{bmatrix}
    -\mathbf{N}\\
   -\mathbf{G}(\mathbf{X})
\end{bmatrix}}^k
\end{split}
\end{align}
$\mathbf{M} \in \mathbb{R}^{N_x \times N_x}$ and $\mathbf{N} \in \mathbb{R}^{N_x \times 1}$ are defined as:
\begin{subequations}
\begin{align}
& \mathbf{M}=\boldsymbol{\mathcal{L}}_{\mathbf{X}\mathbf{X}}^{\gamma}+\mathbf{H}_{\mathbf{X}}^\top\mathbf{diag}(\mathbf{Z})^{-1}\mathbf{diag}(\boldsymbol{\mu}) \mathbf{H}_{\mathbf{X}} \label{M}\\
  \begin{split}
&\mathbf{N}=f_{\mathbf{X}}^\top+\mathbf{G}_{\mathbf{X}}^\top\boldsymbol{\lambda}+\mathbf{H}_{\mathbf{X}}^\top\boldsymbol{\mu}\\ &\qquad \qquad +\mathbf{H}_{\mathbf{X}}^\top\mathbf{diag}(\mathbf{Z})^{-1}(\gamma \mathbf{e} +\mathbf{diag}(\boldsymbol{\mu})\mathbf{H}(\mathbf{X})) \label{N}  \end{split}\\
& \boldsymbol{\mathcal{L}}_{\mathbf{X}\mathbf{X}}^{\gamma} = f_{\mathbf{XX}}+\mathbf{G}_{\mathbf{XX}}(\boldsymbol{\lambda})+\mathbf{H}_{\mathbf{XX}}(\boldsymbol{\mu}) \label{LXX}
\end{align}
\end{subequations}
Solving \eqref{eqn:kktn} numerically results to the locally optimum point $\mathbf{X}^{*}$. By replacing $ \Delta \mathbf{X}^k$ into Eq. \eqref{deltaZ}, $\Delta \mathbf{Z}^k$ is obtained and subsequently $\Delta \mathbf{Z}^k$ into Eq. \eqref{deltamu} finally $\Delta \boldsymbol{\mu}^k$ is computed. Therefore, $\mathbf{X}^{k+1}=\mathbf{X}^k+\alpha \Delta \mathbf{X}^k$, $\boldsymbol{\lambda}^{k+1}=\boldsymbol{\lambda}^k+\alpha \Delta \boldsymbol{\lambda}^k$, $\mathbf{Z}^{k+1}=\mathbf{Z}^k+\alpha \Delta \mathbf{Z}^k$ and $\boldsymbol{\mu}^{k+1}=\boldsymbol{\mu}^k+\alpha \Delta \boldsymbol{\mu}^k$, where $\alpha$ is the step-control parameter which can be chosen arbitrary or based on a optimal multiplier method, which is beyond the scope of this paper\cite{Euzebe_2005}. The detail of numerical solution in Newton's method can be found in \cite{zimmerman_MIPS_2016}.  

It is common to construct Jacobian matrix of \eqref{eqn:kktn} using numerical derivatives and solve the KKT equations using $LU$ factorization\cite{boyd_convex_2004}. However, note that $\mathbf{M}$ is an asymmetric matrix even though it is structurally symmetric, The reason is the term $\mathbf{H}_{\mathbf{X}}^\top\mathbf{diag}(\mathbf{Z})^{-1}\mathbf{diag}(\boldsymbol{\mu}) \mathbf{H}_{\mathbf{X}}$ in Eq. \eqref{M} which makes it asymmetrical. Therefore it is not possible to apply the $LDL^\top$ factorization technique \cite{boyd_convex_2004}. 

\section{Speed up of the solution proposal} 
\label{sec:Speed-up}
{In terms of both memory allocation and computational operations, a MPOPF problem can be attributed to two parts: 1) Input data\footnote{14 input matrices for BATTPOWER solver, see Appendix \ref{Appendix_A}} preparation, and 2) Core optimisation solver which itself consists of: (a) function evaluation, which is the calculus of partial differentials of constraints and objective function w.r.t. all variables in the solution proposal section, (b) computing the inverse of Newton-Raphson Jacobian \eqref{eqn:kktn}, and (c) computational efforts regarding the update of step-control parameter in each iteration.} It is well known that step (b) is the most computationally expensive step of an interior point (IP) algorithm with a large-scale number of variables and constraints \cite{capitanescu_interior-point_2007, capitanescu_experiments_2013}. In this section, first, we mathematically derive the analytical derivative of partial differentiation of constraints w.r.t. all existing variables which, in turn, is the fastest way to solve step (a) \cite{jiang_efficient_2010}, and then, we tailor an algorithm to exploit the structure of KKT systems, specifically using a Schur-Complement approach to accelerate step (b) of the IP method. In the following subsections we will elaborate these two steps. 
\subsection{Analytical Derivatives}
In this subsection, the first and second analytical derivatives of $\mathbf{H}(\mathbf{X})$, $\mathbf{G}(\mathbf{X})$ and $F(\mathbf{X})$ will be extracted. These derivations will be used to construct and consequently solve Eqs. \eqref{M}, \eqref{N}, \eqref{LXX} and finally \eqref{eqn:kktn}. Details of equations regarding the extraction of analytical derivatives can be found in Appendix \ref{Appendix_B} .\\
Furthermore, to exploit the sparsity structure of each block of $\mathbf{G}_{\mathbf{X}}=\frac{\partial \mathbf{G}}{\partial \mathbf{X}}$, $\mathbf{H}_{\mathbf{X}}=\frac{\partial \mathbf{H}}{\partial \mathbf{X}}$, $F_{\mathbf{X}}=\frac{\partial F}{\partial \mathbf{X}}$, $\mathbf{G}_{\mathbf{X}\mathbf{X}}=\frac{\partial}{\partial \mathbf{X}}(\mathbf{G}_\mathbf{X}^\top \boldsymbol{\lambda})$, $\mathbf{H}_{\mathbf{X}\mathbf{X}}=\frac{\partial}{\partial \mathbf{X}}(\mathbf{H}_\mathbf{X}^\top \boldsymbol{\lambda})$, $F_{\mathbf{X}\mathbf{X}}=\frac{\partial}{\partial \mathbf{X}}({F}_\mathbf{X}^\top)$  and their subsequent sub-blocks according to \eqref{19.a}-\eqref{19.e}, a robust and simplified form of structure is developed in order to use highly efficient mathematical operations and matrix substitutions to form \eqref{eqn:kktn}. The general format of these structures can be found in Appendix \ref{Appendix_C} of this paper. 
\subsection{Structure Exploitation and Following Re-Ordering}
\label{reordering}
\begin{figure}[!htbp]
\centering
\includegraphics[width=3.5 in , height=1.8 in]{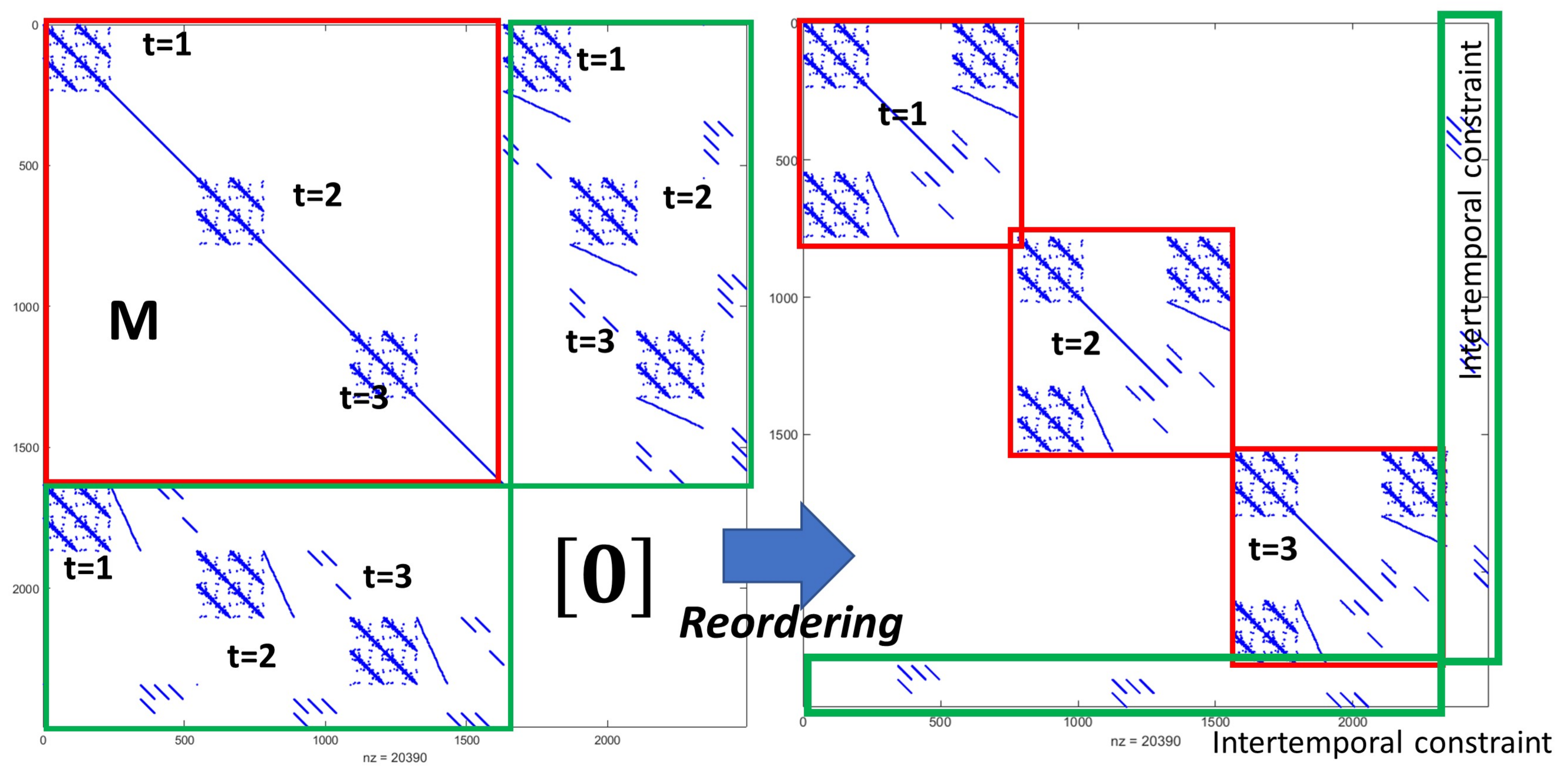}
\caption{Structure of Jacobian of the Newton-Raphson's algorithm before and after reordering.}
\label{fig:1}
\end{figure}

Theoretically, the concept of MPOPF is several snapshots of OPF coupled in time. Each snapshot is a dispatch operational problem with its own variables introduced with \eqref{eqn:x_t}. If we exploit the structure of the Hessian matrix in \eqref{eqn:kktn} with ordering of constraints and variables as shown in (\ref{19}), then, the detailed structure of $\mathbf{M}$ and $\mathbf{G}_\mathbf{X}$ can be seen as:
\begin{align}
\begin{split}
\label{eqn:kktnn}
\begin{bmatrix}
\begin{bmatrix}
 \mathbf{M}_1   &  \\
  & \ddots  \quad   \mathbf{M}_T \\
\end{bmatrix}& 
   \mathbf{G}_\mathbf{X}^\top\\
{\begin{bmatrix}
   \widetilde{\mathbf{G}}_{\mathbf{x}_1}   &  \\
                & \Scale[1]{\ddots} \quad \widetilde{\mathbf{G}}_{\mathbf{x}_T}       \\
 \overline{\mathbf{G}}_{\mathbf{x}_1}    &  \\
 &  \Scale[1]{\ddots} \quad \overline{\mathbf{G}}_{\mathbf{x}_T}\\
  \overline{\mathbf{G}}^s_{\boldsymbol{\tau}_1}    &  \\
 &  \Scale[1]{\ddots} \quad \overline{\mathbf{G}}^s_{\boldsymbol{\tau}_T}
 \end{bmatrix}}& 
    \mathbf{O}
\end{bmatrix}
\Scale[0.8]{\begin{bmatrix}
    \Delta \mathbf{x}_1\\
    \Scale[0.6]{\vdots}\\
    \Delta \mathbf{x}_T\\
   \Delta\widetilde{\boldsymbol{\lambda}}_1\\
  \Scale[0.6]{\vdots}\\
   \Delta\widetilde{\boldsymbol{\lambda}}_T\\
     \Delta\overline{\boldsymbol{\lambda}}_1\\
   \Scale[0.6]{\vdots}\\
   \Delta\overline{\boldsymbol{\lambda}}_T\\
     \Delta\overline{\boldsymbol{\lambda}}^s_1\\
   \Scale[0.6]{\vdots}\\
   \Delta\overline{\boldsymbol{\lambda}}^s_T\\
\end{bmatrix}
=\begin{bmatrix}
    -\mathbf{N}_1\\
    \Scale[0.6]{\vdots}\\
    -\mathbf{N}_T\\
    -\widetilde{\mathbf{g}}({\mathbf{x}_1})\\
    \Scale[0.6]{\vdots}\\
    -\widetilde{\mathbf{g}}({\mathbf{x}_T})\\
    -\overline{\mathbf{g}}({\mathbf{x}_1})\\
    \Scale[0.6]{\vdots}\\
    -\overline{\mathbf{g}}({\mathbf{x}_T})\\
    -\overline{\mathbf{g}}^s({\boldsymbol{\tau}}_1)\\
    \Scale[0.6]{\vdots}\\
    -\overline{\mathbf{g}}^s({\mathbf{\boldsymbol{\tau}}_{T}})
\end{bmatrix}}
\end{split}
\end{align}
where $\mathbf{M}_t \in  \mathbb{R}^{N_{x_t}\times N_{x_t}} $, $ \mathbf{G}_\mathbf{X}^\top \in \mathbb{R}^{N_x\times N_g}$ is the transpose of the left bottom block in the coefficient matrix of \eqref{eqn:kktnn} , $\mathbf{O}\in \mathbb{O}^{N_g\times N_g}$, (set of zeros: $\mathbb{O}$), $\widetilde{\mathbf{G}}_{\mathbf{x}_t}= \frac{\partial \widetilde{\mathbf{g}}}{\partial \mathbf{x_t}} \in  \mathbb{R}^{N_x\times n_{gn}}$, $\overline{\mathbf{G}}_{\mathbf{x}_t}= \frac{\partial \overline{\mathbf{g}}}{\partial \mathbf{x_t}} \in  \mathbb{R}^{N_x\times n_{{gl}_t}}$. As defined in section \ref{Sec:multiAC}, for the sake of simplicity of notation we take  $\overline{\mathbf{G}}^s_{\boldsymbol{\tau}_t}=\frac{\partial \overline{\mathbf{g}}^s}{\partial \boldsymbol{\tau}_t} = \frac{\partial \overline{\mathbf{g}}^s}{\partial (\mathbf{x}_{t-1},\mathbf{x}_{t}) } \in  \mathbb{R}^{N_x\times n_{gs}}$. If we reorder the current variables and consequently re-construct the coefficient and righthand side matrix such that all variables corresponding to time $t$ are assembled together except the variables of inter-temporal constraints, then the vector of variables with its corresponding righthand side and coefficient matrix could be written as \eqref{var}, \eqref{reighthand} and \eqref{coefficient}, respectively.
\begin{subequations}\label{reordered}
\begin{alignat}{4}
    &\Scale[0.95]{\begin{bmatrix}
    \Delta \mathbf{x}_1 \ \Delta\widetilde{\boldsymbol{\lambda}}_1 \ \Delta\overline{\boldsymbol{\lambda}}_1 \Scale[1]{\dots}\Delta \mathbf{x}_T \ 
   \Delta\widetilde{\boldsymbol{\lambda}}_T \ \Delta\overline{\boldsymbol{\lambda}}_T
      ,\quad\Delta\overline{\boldsymbol{\lambda}}^s_1    \Scale[0.6]{\dots}   \Delta\overline{\boldsymbol{\lambda}}^s_T
\end{bmatrix}^\top} \label{var}\\
   \Scale[0.5]{\pmb{-}}&\Scale[0.92]{\begin{bmatrix}
    \mathbf{N}_1 \ \widetilde{\mathbf{g}}({\mathbf{x}_1})\ \overline{\mathbf{g}}({\mathbf{x}_1}) \Scale[0.5]{\dots} \mathbf{N}_T \ 
   \widetilde{\mathbf{g}}({\mathbf{x}_T}) \ \overline{\mathbf{g}}({\mathbf{x}_T})
      , \overline{\mathbf{g}}^s({\boldsymbol{\tau}_1})   \Scale[0.5]{\dots}   \overline{\mathbf{g}}^s({\boldsymbol{\tau}_T})
\end{bmatrix}}^\top \label{reighthand}\\
&\begin{bmatrix}
\begin{bmatrix}
      \mathbf{\Upsilon}_1           &  \\
&\Scale[0.6]{\ddots}  \quad \mathbf{\Upsilon}_T    \\
\end{bmatrix}&
   {\overline{\mathbf{G}}_{\mathbf{X}}^{sr}}^\top\\
\begin{bmatrix}
  \overline{\mathbf{G}}^{sr}_{\boldsymbol{\tau}_1}    & \\
 &  \Scale[0.6]{\ddots} \quad \overline{\mathbf{G}}^{sr}_{\boldsymbol{\tau}_T} 
\end{bmatrix}&\mathbf{O}^r
\end{bmatrix}\label{coefficient}
\end{alignat}
\end{subequations}
where $\Scale[0.88]{
    \mathbf{\Upsilon}_t = \begin{bmatrix}
 \mathbf{M}_t     &\begin{bmatrix}
    \widetilde{\mathbf{G}}_{\mathbf{x}_t}^\top &\overline{\mathbf{G}}_{\mathbf{x}_t}^\top 
 \end{bmatrix} \\
\begin{bmatrix}
   \widetilde{\mathbf{G}}_{\mathbf{x}_t}\\
   \overline{\mathbf{G}}_{\mathbf{x}_t}
\end{bmatrix}  &{\mathbf{O}^{r}}'\\
    \end{bmatrix}} \in \mathbb{R}^{N_{\mathbf{\Upsilon}_t} \times N_{\mathbf{\Upsilon}_t}}$, $N_{\mathbf{\Upsilon}_t}=N_{x_t}+n_{gn}+n_{{gl}_t}$, ${\mathbf{O}^{r}}' \in \mathbb{O}^{[n_{gn}+n_{{gl}_t}]\times [n_{gn}+n_{{gl}_t}]}$ and  $\Scale[0.88]{\overline{\mathbf{G}}_\mathbf{X}^{sr} = \begin{bmatrix}
  \overline{\mathbf{G}}^{sr}_{\boldsymbol{\tau}_1}    & \\
 &  \Scale[0.6]{\ddots} \quad \overline{\mathbf{G}}^{sr}_{\boldsymbol{\tau}_T} 
\end{bmatrix}} \in \mathbb{R}^{N_{gs} \times N_{gsr}}$, $N_{gsr}=N_{\mathbf{\Upsilon}_{t=1}}+N_{\mathbf{\Upsilon}_{t=2}}+\dots+N_{\mathbf{\Upsilon}_{t=T}}$, $\mathbf{O}^r \in \mathbb{O}^{N_{gs}\times N_{gs}}$.
In order to illustrate the re-ordering more clearly, we define: ${\delta \pmb{\omega}}_t=\begin{bmatrix} \Delta \mathbf{x}_t \Delta\widetilde{\boldsymbol{\lambda}}_t\Delta\overline{\boldsymbol{\lambda}}_t\end{bmatrix}^\top \in \mathbb{R}^{N_{\mathbf{\Upsilon}_t} \times1}$,  ${\delta \boldsymbol{\lambda}}=\begin{bmatrix} \Delta\overline{\boldsymbol{\lambda}}^s_1    \Scale[0.6]{\dots}   \Delta\overline{\boldsymbol{\lambda}}^s_T\end{bmatrix}^\top \in \mathbb{R}^{N_{gs} \times1}$, ${\pmb{\zeta} }_t=-\begin{bmatrix} \mathbf{N}_t \widetilde{\mathbf{G}}_{\mathbf{x}_t}\overline{\mathbf{G}}_{\mathbf{x}_t}\end{bmatrix}^\top \in \mathbb{R}^{N_{\mathbf{\Upsilon}_t} \times1}$ and ${ \mathbf{\Gamma}}=-\begin{bmatrix} \overline{\mathbf{g}}^s(\boldsymbol{\tau}_1)   \Scale[0.8]{\dots}   \overline{\mathbf{g}}^s(\boldsymbol{\tau}_T)\end{bmatrix}^\top \in \mathbb{R}^{N_{gs} \times1}$. Therefore we can convert \eqref{reordered} to \eqref{commonStruct} which is a well-known ``arrowhead" structure that can be found in the literature \cite{petra_augmented_2014,petra_real-time_2014}.
\begin{align}
\begin{split}
\label{commonStruct}
\begin{bmatrix}
 \mathbf{\Upsilon}_1   &    & & &\pmb{\rho}_1^\top\\
  & \mathbf{\Upsilon}_2 & && \pmb{\rho}_2^\top\\
  & & \ddots&&\vdots\\
  &&&\mathbf{\Upsilon}_T&\pmb{\rho}_T^\top\\
  \pmb{\rho}_1&\pmb{\rho}_2&\dots&\pmb{\rho}_T&0
\end{bmatrix}
\Scale[1]{\begin{bmatrix}
    \delta \mathbf{\pmb{\omega}}_1\\
    \delta \pmb{\omega}_2\\
   \vdots\\
   \delta\pmb{\omega}_T\\
     \delta \boldsymbol{\lambda}\\
\end{bmatrix}
=\begin{bmatrix}
      \pmb{\zeta}_1\\
   \pmb{\zeta}_2\\
   \vdots\\
  \pmb{\zeta}_T\\
     \pmb{\Gamma}\\
\end{bmatrix}}
\end{split}
\end{align}
where the coupling matrices of $\pmb{\rho}_t \in \mathbb{R}^{N_{gs} \times N_{\mathbf{\Upsilon}_t}}$ are:
\begin{equation}\label{rho}
    \pmb{\rho}_1=\begin{bmatrix}
       \overline{\mathbf{G}}^s_{\boldsymbol{\tau}_1}\\
       0\\
       0\\
       0\\
       0\\
    \end{bmatrix}, \pmb{\rho}_2=\begin{bmatrix}
       0\\
       \overline{\mathbf{G}}^s_{\boldsymbol{\tau}_2}\\
       0\\
       0\\
       0\\
    \end{bmatrix}, \pmb{\rho}_T=\begin{bmatrix}
       0\\
       0\\
       0\\
       0\\
       \overline{\mathbf{G}}^s_{\boldsymbol{\tau}_T}
    \end{bmatrix}
\end{equation}
Fig. \ref{fig:1} illustrates the reordering of Hessian matrix of Eqs. \eqref{eqn:kktnn} to \eqref{commonStruct}. In the next subsection, we propose a Schur-Complement technique tailored for the reordered structure of Eq. \eqref{commonStruct} to save computational time. 
Considering Eqs. \eqref{reordered} and \eqref{commonStruct}, it should be kept in mind that in the introduced re-ordering structure, there are some vectors which could have a different length in each time $t$ as discussed in Section \ref{Sec:multiAC}  when the availability of either storage devices, charge, discharge, reactive power provision, and generators would alternate over the optimisation horizon ($\mathbf{AVBP}\big\rvert_{t=1}\neq \mathbf{AVBP}\big\rvert_{t=2}\neq ...\neq \mathbf{AVBP}\big\rvert_{t=T}$), and consequently, ${n_{gl}}_t$ and ${n_{hl}}_t$ are not constant through time $t = \{1,...,T\}$. Therefore, specific indices are introduced to keep track of them at each time and over each iteration $k$, in Eq. \eqref{eqn:kktnn}. 
\subsection{Schur-Complement Technique}\label{ShcurComplement}
The sparse arrowhead structure of the coefficient matrix of \eqref{commonStruct}, is suited for block elimination using the Schur-Complement technique \cite{boyd_convex_2004}. The algorithm is tailored for solving any problem featured with repeatable matrices and coupling constraints  between them, such as: (a) multi-period systems as this paper presents, and in the literature  \cite{kourounis_towards_2018}, (b) stochastic problems with a large number of scenarios \cite{petra_augmented_2014, petra_real-time_2014}, and (c) security constrained problems with a large number of contingencies \cite{kardos_structure_2019}. Algs. \ref{Alg1} and \ref{Alg2} are proposed to solve for a generic multi-period KKT system with a structurally symmetric structure. Alg. \ref{Alg1} is for Schur-Complement factorisation and Alg. \ref{Alg2} is for forward and backward substitution. The substructures of  \eqref{commonStruct} are inputs to both Algs. \ref{Alg1} and \ref{Alg2}.
\subsubsection{Alg. \ref{Alg1}: Schur-Complement Factorisation}
\begin{figure}[!htbp]
\centering
\includegraphics[width=2.5 in , height=2.1 in]{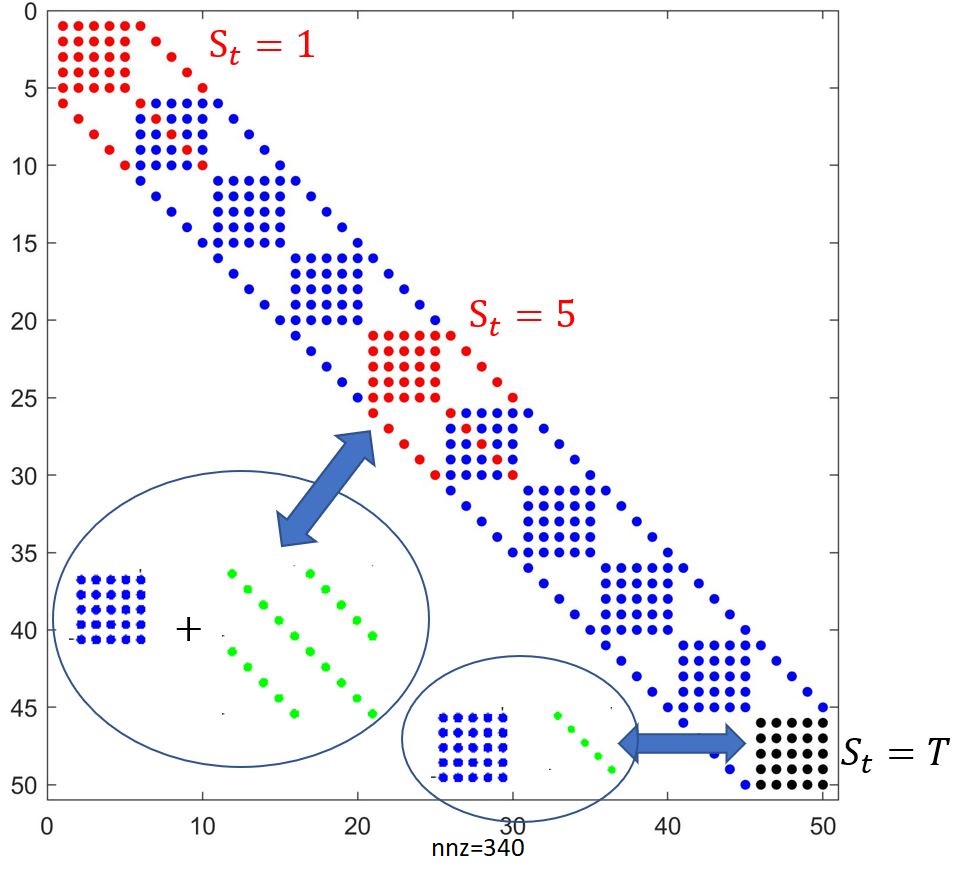}
\caption{Overall structure of the main Schur-Complement: $\pmb{\sigma}^c$. The structure of $\mathbf{S}_{t}$ is shown with red dots, and $\mathbf{S}_{t=T}$ with black dots, for a network with $n_y = 5$ and $T=10$. \footnotemark} 
\label{fig:SchurComp1}
\end{figure}

\begin{figure}[!htbp]
\centering
\includegraphics[width=3.5 in , height=2.5 in]{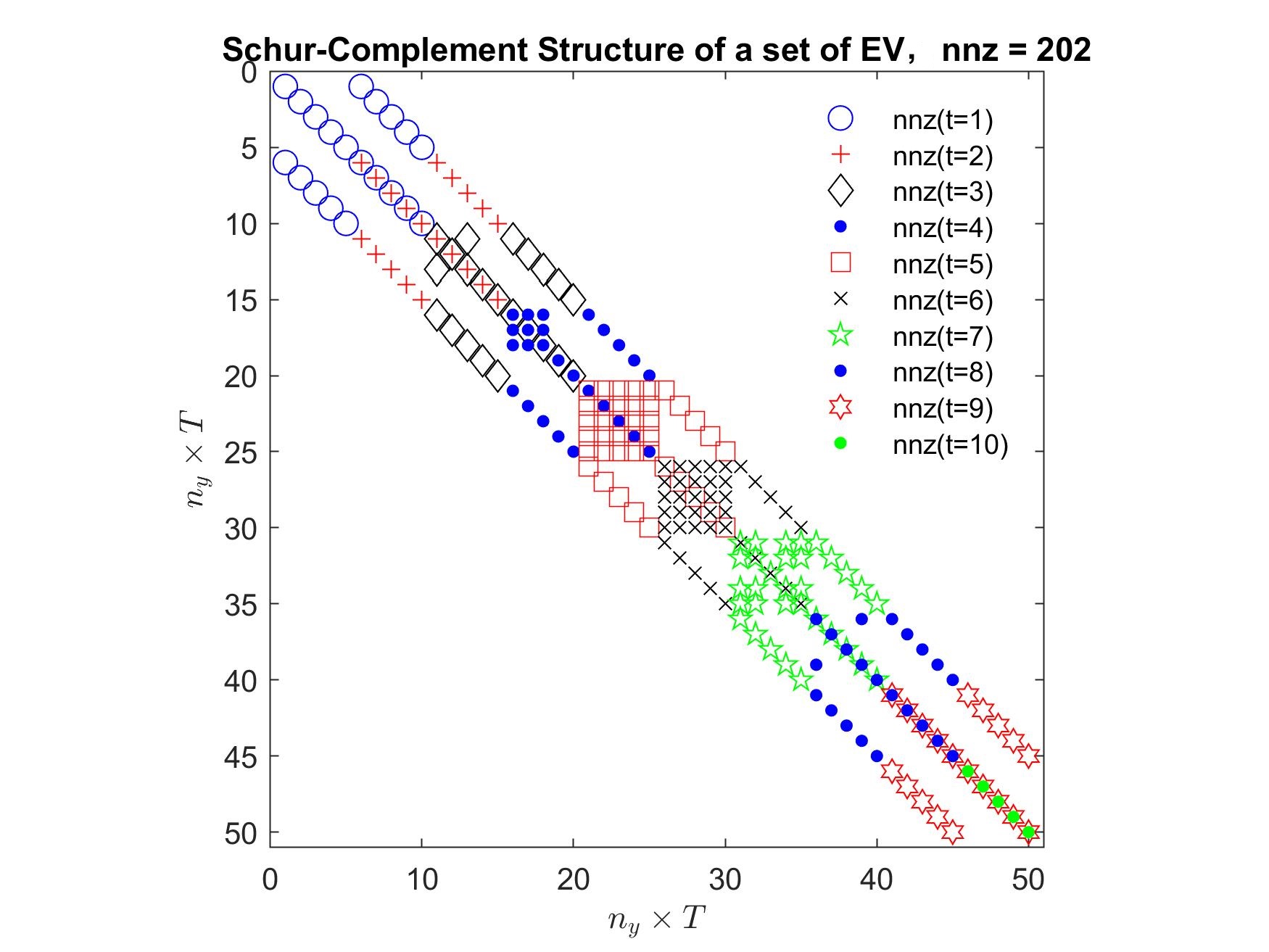}
\caption{Overal structure of the main Schur-Complement: $\pmb{\sigma}^c$, for a network with $n_y = 5$ and $T=10$, corresponding to input matrices of \eqref{eqn:structureEV1}-\eqref{eqn:structureEV2} \footnotemark. Number of non-zero elements of $\mathbf{S}_{t}$ changes through time $t = \{1,...,T\}$.} 
\label{fig:SchurComp}
\end{figure}

Alg. 1 starts with finding a permutation matrix of $\mathbf{Q}_t^{am} \in \mathbb{B}^{N_{\mathbf{\Upsilon}_t} \times N_{\mathbf{\Upsilon}_t}}$, generated in order to capture the sparse structure of $\mathbf{\Upsilon}_t$ to reduce the number of non-zeroes in LU factorisation in L.\ref{factorize}. Permutation matrix $\mathbf{Q}_t^{am}$ is produced based on an approximate minimum degree permutation method \cite{noauthor_approximate_nodate}. It should be kept in mind that the structure of $\mathbf{Q}_t^{am}$ is dependent on $\mathbf{\Upsilon}_t$ which, in turn, is dependent on the structure of input matrices and the conditions of \eqref{eqn:availability}-\eqref{eqn:availabilityDI}. If these conditions hold (all storage devices are SESS), then $N_{\mathbf{\Upsilon}_t}$ is constant through time $t$. Therefore,   $\mathbf{\Upsilon}_t$ and $\mathbf{Q}_t^{am}$ have constant structures over time and L.\ref{permute2} will not be executed; in other words $\mathbf{Q}_1^{am}=\mathbf{Q}_2^{am}=\mathbf{Q}_t^{am}$. If these conditions do not hold, which in turn mean dynamic storage devices (EV), then L.\ref{permute1} is not executed and instead, L.\ref{permute2} will be executed. The factorisation here is based on an \textit{incomplete augmented factorisation} technique \cite{kourounis_towards_2018} in order to compute the Schur-Complement of each augmented matrix of $\mathbf{A}^a_t=\begin{bmatrix}\mathbf{\Upsilon}_t &\pmb{\rho}_t^\top\\ \pmb{\rho}_t &0 \end{bmatrix}$. Permuted $\mathbf{\Upsilon}_t$ factorises with LU factorisation technique in L.\ref{factorize} and afterwards, the Schur-Complement of each block of $\mathbf{A}^a_t$ is computed as $\mathbf{S}_t=-\pmb{\rho}_t\mathbf{\Upsilon}_{t}^{-1}\pmb{\rho}_t^\top \in \mathbb{R}^{N_{gs} \times N_{gs}} $ in L.\ref{Schur_Comp_block} in each iteration and summed together in $\pmb{\sigma}^c \in \mathbb{R}^{N_{gs} \times N_{gs}}$ L.\ref{sigmac} to shape the main Schur-Complement of arrowhead structure of Eq. \eqref{commonStruct}. \\ 
\begin{algorithm}\caption{Schur-Complement Factorization}
  \DontPrintSemicolon \label{Alg1}
  \SetKwFunction{FMain}{SchurComI}
  \SetKwProg{Fn}{Function}{:}{}
  \Fn{\FMain{$\{\mathbf{\Upsilon}_t$, $\pmb{\rho}_t$, $\pmb{\zeta}_t$, \ $\forall \ t=1,\dots, T\}$ }}{
  [$\mathbf{Q}_t^{am}$]=\footnotesize{ApproxMinDegrPermut}($\mathbf{\Upsilon}_1$)\tcp*[h]{\Scale[0.9]{\mathrm{If}\ \eqref{eqn:availability}-\eqref{eqn:availabilityDI} \ \mathrm{hold}}}\normalsize\; \label{permute1}
  $\pmb{\sigma}^c=0$\;
  $\pmb{\sigma}^l=0$\;
  \For{t$\mathrm{=1:}$T}{
        [$\mathbf{Q}_t^{am}$]=ApproxMinDegrPermut($\Scale[0.9]{\mathbf{\Upsilon}_t}$)\;\label{permute2} \qquad \tcp*[h]{\Scale[0.8]{\mathrm{If}\ \eqref{eqn:availability}-\eqref{eqn:availabilityDI} \ \mathrm{does\ NOT\ hold}}}\; 
        $\mathbf{\Pi}$ = \normalsize {SparsePermute}(${\big[[\mathbf{Q}_t^{am}}]^\top \mathbf{\Upsilon}_t\mathbf{Q}_t^{am}\big]$)\;
        [$\mathbf{L}^{lu}_t$, $\mathbf{U}^{lu}_t$, $\mathbf{P}_t^{lu}$, $\mathbf{Q}_t^{lu}$, $\mathbf{R}^{lu}_t$] \;\qquad=  SparseLUFactorize($\mathbf{\Pi}$)\;\label{factorize}
        $\mathbf{inf}^\mathbf{\Upsilon}_t$ =\small
 struct\normalsize ($\mathbf{L}^{lu}_t$, $\mathbf{U}^{lu}_t$, $\mathbf{P}_t^{lu}$, $\mathbf{Q}_t^{lu}$, $\mathbf{R}^{lu}_t$, $\mathbf{Q}_t^{am}$)\; \label{save}
        $\mathbf{S}_t=-\pmb{\rho}_t\mathbf{\Upsilon}_{t}^{-1}\pmb{\rho}_t^\top$ \tcp*[h]{SchurCompAuxiliary $\mathbf{A}^a_t$}\; \label{Schur_Comp_block}
        $\pmb{\sigma}^c = \pmb{\sigma}^c+\mathbf{S}_t$ \tcp*[h]{\Scale[0.9]{Main SchurCompArrowhead-\ref{commonStruct}}}\;\label{sigmac}
        $\mathbf{\Xi}_t= -\pmb{\rho}_t\mathbf{\Upsilon}_{t}^{-1}\pmb{\zeta}_t $ \tcp*[h]{SchurCompAuxiliary $\mathbf{A}^b_t$}\;\label{12}
        $\pmb{\sigma}^l= \pmb{\sigma}^l+\mathbf{\Xi}_t$\tcp*[h]{righthand side Alg 2,L.\ref{Lambda}} \; \label{sigmal}
        \textbf{End for}\;
  }
  [$\mathbf{L}^{ldl}$, $\mathbf{D}^{ldl}$, $\mathbf{P}^{ldl}$, $\mathbf{S}^{ldl}$] =SparseLDLFactorize($\pmb{\sigma}^c$)\;\label{ldlfactorize}
  $\mathbf{inf}^c$ = struct ($\mathbf{L}^{ldl}$, $\mathbf{D}^{ldl}$, $\mathbf{P}^{ldl}$, $\mathbf{S}^{ldl}$)\;
  \KwRet $\{\pmb{\sigma}^l$, $\mathbf{inf}^c\}$ and $\{\mathbf{inf}^\mathbf{\Upsilon}_t$, \ $\forall \ t=1,\dots, T\}$   \;
  \textbf{End Function}\; 
}\end{algorithm} 
With almost the same procedure explained above, in order to compute $\pmb{\xi}$ which is the righthand side of the main Schur-Complement equation in $[\pmb{\sigma}^c][\delta\boldsymbol{\lambda}]= [\pmb{\xi}]$, (refer to Alg 2, L.\ref{Lambda}), we define another auxiliary block matrix of $\mathbf{A}^b_t=\begin{bmatrix}\mathbf{\Upsilon}_t &\pmb{\zeta}_t\\ \pmb{\rho}_t &0 \end{bmatrix}$ and consequently compute its Schur-Complement as $\mathbf{\Xi}_t= -\pmb{\rho}_t\mathbf{\Upsilon}_{t}^{-1}\pmb{\zeta}_t$ where $\mathbf{\Upsilon}_{t}^{-1}$ is factorised in previous step L.\ref{factorize}. Thus, we only recall the stored ``struct" of $\mathbf{inf}^\mathbf{\Upsilon}_t$ L.\ref{save} from memory. Value of $\mathbf{\Xi}_t$ is aggregated in each iteration with $\pmb{\sigma}^l$ in L.\ref{sigmal}. $\pmb{\sigma}^c$ is the Schur-Complement of arrowhead structure of Eq. \eqref{commonStruct} and has an interesting pattern that can be exploited further, and it is dependent on the input matrices of $\mathbf{AVBP}$, $\mathbf{CONCH}$ and $\mathbf{CONDI}$.
\paragraph{Static Schur-Complement Structure:  {SESS}}\label{ESS} If $\{\mathbf{AVBP}$,$\mathbf{CONCH}$,$\mathbf{CONDI}\}\in \{\mathbb{1}\}^{n_y \times T}$ holds, then all the storage devices are considered as {SESS}. As elaborated above, $\pmb{\sigma}^c$ is the aggregation of Schur-Complement of each auxiliary block of $\mathbf{A}^a_t$, therefore, in each iteration $\mathbf{S}_t=-\pmb{\rho}_t\mathbf{\Upsilon}_{t}^{-1}\pmb{\rho}_t^\top $. $\pmb{\sigma}^c$ is a sparse bandwidth matrix such that $\{\forall i,j \ \pmb{\sigma}_{i,j}^c=0$ if $\abs{i-j}>n_y\}$  where the number of non-zero matrix elements is only dependent on the number of storage devices $n_y$ and simulation horizon $T$, and not on network properties. More precisely, each element of ${s}_{ij,t}$ is computed through  ${s}_{ij,t}=-\pmb{\rho}_{i,t}\mathbf{\Upsilon}_{t}^{-1}\pmb{\rho}_{j,t}^\top$ considering the only non-zero part of $\pmb{\rho}_{t}$ is $\overline{\mathbf{G}}^s_{\boldsymbol{\tau}_t}$ which moves from top to bottom while $t$ moves from $t=1$ to $t=T$ as illustrated in Eq. \eqref{rho} and in the for loop of Alg. \ref{Alg1}, with the essential assumption that the condition $\mathbf{AVBP}= \mathbf{CONDI} =\mathbf{CONCH}= \begin{bmatrix}
\mathbb{1}
\end{bmatrix}_{n_y\times T}$ holds, which intuitively means that $\overline{\mathbf{G}}^s_{\boldsymbol{\tau}_t}$ has a constant structure over optimisation horizon. $\pmb{\rho}_{t}^\top$ also has the same pattern through a loop from $t=1$ to $t=T$. {Therefore, $\mathbf{S}_{t}$ moves in a bandwidth structure as shown in Fig. \ref{fig:SchurComp1} from the top left corner to the bottom right. Each block of $\mathbf{S}_{t\neq T}$ is constructed by two sub-blocks of $\overline{\mathbf{G}}^s_{\boldsymbol{\tau}_t}$, the blue rectangular sub-block with $\boldsymbol{\mathcal{P}}^{ch}$ and $\boldsymbol{\mathcal{P}}^{dch}$ variables and the light green sub-block with $\boldsymbol{\mathcal{SOC}}$ variables at each time. The structure of $\mathbf{S}_{t=T}$ is different since $\overline{\mathbf{G}}^s_{\boldsymbol{\tau}_T}$ is the last linking equation over time, and therefore,  $\boldsymbol{\mathcal{SOC}}$ is no longer linked to a next time, but  $\boldsymbol{\mathcal{P}}^{ch}$ and $\boldsymbol{\mathcal{P}}^{dch}$ still exist.} 
\paragraph{Dynamic Schur-Complement Structure---EV} \label{EV}
\footnotetext{{$nnz=340$ stands for number of non-zero elements}}If conditions of \eqref{eqn:availability}-\eqref{eqn:availabilityDI} {do not} hold, the structure shown in Fig. \ref{fig:SchurComp1} would change depending on the dynamic behaviour of the input matrices of $\mathbf{AVBP}$, $\mathbf{CONCH}$ and $\mathbf{CONDI}$. Fig. \ref{fig:SchurComp} illustrates the sparse Schur-Complement pattern of $\pmb{\sigma}^c$, made by input matrices shown in Eqs.\eqref{eqn:structureEV1}-\eqref{eqn:structureEV2}. Each colour/legend shows the non-zero elements of $\mathbf{S}_{t}$ for each time $t={1,...,T}$, which show up in the main Schur-Complement structure $\pmb{\sigma}^c$. As presented in Alg. \ref{Alg1}, L.\ref{sigmac}, the structure of $\pmb{\sigma}^c$, and also as shown in Fig. \ref{fig:SchurComp}, is constructed inside a for loop from $t=1$ to $t=T$. The non-zero elements of $\mathbf{S}_{t}$  have overlaps with non-zero elements of $\mathbf{S}_{t-1}$ and $\mathbf{S}_{t+1}$, see Fig. \ref{fig:SchurComp}. Therefore, they are added together in the for loop while $t$ progresses from $t=1$ to $t=T$. In this algorithm, it is important to allocate a size of memory proportional to the number of non-zero elements of $\mathbf{S}_{t}$ and $\pmb{\sigma}^c$ to achieve high performance. It should be noted that the number of non-zero elements in $\mathbf{S}_{t}$ and $\pmb{\sigma}^c$ (and the size of allocated memory) for each time $t={1,...,T}$ are in turn function of input matrices and they can be predetermined before the starting of operation of Alg. \ref{Alg1}.
\begin{subequations}\label{eqn:structureEV}
\begin{alignat}{4} 
& \Scale[0.85]{\mathbf{AVBP}= \mathbf{CONCH}=
\begin{bmatrix}
0&	0	&1	&1	&1	&1	&1	&1&	0	&0 \\
0&	0&	0	&1	&1&	1	&1&	0	&0&	0\\
0&	0	&1	&1	&1	&1	&0	&0&	0&	0\\
0	&0	&0&	0	&1	&1&	1	&1	&1&	0\\
0	&0	&0	&0	&1&	1&	1	&0	&0	&0\\
\end{bmatrix}}\label{eqn:structureEV1}\\
&\Scale[0.85]{\mathbf{CONDI}=\begin{bmatrix}
\mathbb{0}
\end{bmatrix}_{(n_y=5)\times (T=10)}\label{eqn:structureEV2}}
\end{alignat}
\end{subequations}
{ There is no EV in the system at time $t=1,2$; two EV appear at time $t=3$, full EV connected system at $t=5,6$; similarly, corresponding structures alter, as can be seen in Fig. \ref{fig:SchurComp1}.}
\footnotetext{{$nnz=202$ stands for number of non-zero elements, different value for each time}} { The reason for the different structure shown Fig. \ref{fig:SchurComp1} is the same reason elaborated in Section \ref{ESS}. The only non-zero part of $\pmb{\rho}_{t}$ is $\overline{\mathbf{G}}^s_{\boldsymbol{\tau}_t}$ which moves from top to bottom while $t$ moves from $t=1$ to $t=T$ illustrated in  Eq. \eqref{rho} and in the for loop of Alg. \ref{Alg1}. $\overline{\mathbf{G}}^s_{\boldsymbol{\tau}_t}$ is the first derivative of energy storage systems in the time $t$ and its structure alternates over time, due to the dynamic behaviour of input matrices of $\mathbf{AVBP}$, $\mathbf{CONCH}$ and $\mathbf{CONDI}$. Therefore, the Schur-Complement structure of $\pmb{\sigma}^c$ gets more sparse than that of SESS and could be factorised even faster than that of SESS.} \\
In the rest of Alg. \ref{Alg1}, $\pmb{\sigma}^c$ is factorised by sparse LDL factorisation technique, since it is a sparse symmetric matrix. Finally Alg. \ref{Alg1} returns $\pmb{\sigma}^l$ matrix, and $\mathbf{inf}^c$ and $\mathbf{inf}^\mathbf{\Upsilon}_t$ structs containing factorisation information to be called in Alg. \ref{Alg2}.
\subsubsection{Alg. \ref{Alg2} Forward and Backward Substitution}
Inputs are matrices of $\mathbf{\Gamma},\pmb{\sigma}^l, \pmb{\rho}_t,\pmb{\zeta}_t$,  and  structs of $\mathbf{inf}^c$, $\mathbf{inf}^\mathbf{\Upsilon}_t$. Alg. \ref{Alg2} is to solve $\delta{\boldsymbol{\lambda}}$ through sparse LDL forward and backward substitution, and further, to compute $\delta \pmb{\omega}_t$ through a for loop with the help of sparse LU forward and backward substitution. \\
First the righthand side of the main Schur-Complement equation of $\pmb{\sigma}^c \delta\boldsymbol{\lambda}= \pmb{\xi}$, which is $\pmb{\xi}$,  is computed in Alg 2, L.\ref{Lambda}. Using the ``struct" of $\mathbf{inf}^c$, the sparse LDL forward and backward substitution is performed in L.\ref{solve1} and $\delta{\boldsymbol{\lambda}}$ is cleared. Second, in a for loop, a slack variable called $\pmb{\kappa}_t$ is constructed with $\pmb{\kappa}_t=\pmb{\zeta}_t-\pmb{\rho}_t^\top\delta{\boldsymbol{\lambda}}$ in L.\ref{slack} and recalled to solve for the sparse LU forward and backward substitution using the ``struct" of $\mathbf{inf}^\mathbf{\Upsilon}_t$ in each $t$ and thus, $\delta \pmb{\omega}_t$ is cleared. Finally, Alg. \ref{Alg2} returns vectors $\delta{\boldsymbol{\lambda}}$ and $\delta \pmb{\omega}_t$. 
\begin{algorithm}\caption{Forward and Backward Substitution}
  \DontPrintSemicolon \label{Alg2}
  \SetKwFunction{FMain}{SchurComII}
  \SetKwProg{Fn}{Function}{:}{}
  \Fn{\FMain{$\{\mathbf{\Gamma},\pmb{\sigma}^l$, $\mathbf{inf}^c\}$, $\{\pmb{\rho}_t,\pmb{\zeta}_t, \mathbf{inf}^\mathbf{\Upsilon}_t$, \ $\forall \ t=1,\dots, T\}$ }}{
  $\pmb{\xi}=\mathbf{\Gamma}-\pmb{\sigma}^l$\;\label{Lambda}
  $\delta{\boldsymbol{\lambda}}$=SparseLDLForBackSolve($\mathbf{inf}^c$, $\pmb{\xi}$)\;\label{solve1}
  \For{t$\mathrm{=1:}$T}{
        $\pmb{\kappa}_t=\pmb{\zeta}_t-\pmb{\rho}_t^\top\delta{\boldsymbol{\lambda}}$\; \label{slack}
        $\delta \pmb{\omega}_t$=  SparseLUForBackSolve($\mathbf{inf}^\mathbf{\Upsilon}_t, \pmb{\kappa}_t$)\;\label{solve2}
        \textbf{End for}\;
  }
  \KwRet $\{\delta{\boldsymbol{\lambda}}\}$ and $\delta \pmb{\omega}_t$, \ $\forall \ t=1,\dots, T\}$   \;
  \textbf{End Function}\;
}\end{algorithm}
\subsection{Computational Performance and Memory Efficiency}
\subsubsection{Sparse Matrix Operations}
Since most of the operations here are done in sparse format efficiently, sparse indexing and libraries are developed to construct and handle re-ordering explained above. Sparsity is a method to save significant memory in large-scale simulations.  
\subsubsection{Function Evaluation}
First and foremost, analytical derivative functions, known as hand-coded functions, are implemented here, and are the fastest possible way to compute all partial derivatives of the first and second order of objective function and constraints w.r.t. all variables, especially when it comes to large-scale optimisation systems \cite{jiang_efficient_2010}. Efficiently exploiting and handling the operations of sparse matrices and subsequently updating their indices over the optimisation time horizon accelerates the computational performance.  
\subsubsection{KKT Systems} \label{KKTSystems}
The computational complexity of sparse matrix operations is proportional to the number of non-zero elements in each sparse matrix and independent of the size of matrices \cite{gilbert_sparse_1992}. Note that the number of non-zero elements of Jacobian matrix for each case and per row is almost constant and is similar to the number of non-zero elements of a coefficient matrix obtained from discretisation of finite element methods in three-dimensional meshes.\\
Therefore, the only possible way to assess the performance of a solution proposal for KKT systems is through the estimation of sparse matrix operations. Factorisation L.\ref{factorize} in Alg. \ref{Alg1} is the most computationally expensive step in the loop when $n_y \leq \mathcal{K}$,  where $\mathcal{K} \propto N_{\mathbf{\Upsilon}_t}=N_{x_t}+n_{gn}+n_{{gl}_t}= 2n_b+2n_g+4n_y+n_{gn}+n_{{gl}_t}$ \footnote{In other words, Factorisation L.\ref{factorize} is the most computationally expensive step, when the number of storage devices are smaller than a certain number, where this number is proportional to size of blocks of $\mathbf{\Upsilon}_t$ in Eq. \eqref{commonStruct}}. Note that the rest of Alg. \ref{Alg1} and Alg. \ref{Alg2} are dominated by this step and can be ignored in this case. If the number of non-zero elements per row of coefficient matrix of \eqref{rho} is approximately similar to Laplacian matrices discretised by finite elements, then LU factorisation of L.\ref{factorize} has the complexity of $O(N_{\mathbf{\Upsilon}_t}^2)$. Consequently, the complexity of factorisation for all blocks of   $\mathbf{\Upsilon}_t$ will be $O(TN_{\mathbf{\Upsilon}_t}^2)$. \\
However, for $n_y > \mathcal{K}$ the complexity of Alg. \ref{Alg1} is dominated by (a) L.\ref{Schur_Comp_block}, (b) L.\ref{12}, since the size of matrix $\pmb{\rho}_t \in \mathbb{R}^{N_{gs} \times N_{\mathbf{\Upsilon}_t}}$ gets larger where $N_{gs}=Tn_y$, and (c) $LDL^\top$ factorisation of $\pmb{\sigma}^c$, in L.\ref{ldlfactorize}; this is because the main Schur-Complement structure of $\pmb{\sigma}^c$  becomes dense and therefore, the complexity of $LDL^\top$ will be $O(\frac{1}{3}(Tn_y)^3)$. Note that overhead of factorisation of $\pmb{\sigma}^c$ is more dependent on $n_y$ than on $T$ since non-zero elements of each sub-structure of $\mathbf{S}_t$ are $\{nnz(\mathbf{S}_t)=n_y(n_y+3) \ | \ t \neq T,\ nnz(\mathbf{S}_T)=n_y^2\}$. Note that these arguments are only valid when $\mathbf{AVBP}= \mathbf{CONDI} =\mathbf{CONCH}= \begin{bmatrix}
\mathbb{1}
\end{bmatrix}_{n_y\times T}$ holds, which in turn would lead to a static Schur-Complement structure discussed in Section \ref{ESS}. On the contrary, the complexity of dynamic Schur-Complement structure shown in Section \ref{EV} has even less overhead than the static one. 
{\subsubsection{Memory Efficiency}\label{Memory}
\paragraph{{Input Matrices}}\label{InputMatirice}
BATTPOWER has 14 input matrices in total, as described in the Appendix \ref{Appendix_A}: Four are similar to MATPOWER ($\matr{BUS}$, $\matr{BRANCH}$, $\matr{GEN}$, $\matr{GENCOST}$) and ten new matrices to capture the mulitperiod formulation and energy storage. Except the mentioned MATPOWER matrices and the new matrises $\mathbf{PD}$ and $\mathbf{QD}$, which are all dense matrices, the rest is stored in a memory-efficient manner, either binary or spare format matrices. $\mathbf{AVBP}$, $\mathbf{CONCH}$, $\mathbf{CONDI}$, $\mathbf{AVBQ}$, and $\mathbf{AVG}$ matrices are the binary ones. Finally, $\mathbf{SOCi}$ and $\mathbf{SOCMi}$ are neither dense nor binary, so these are stored with sparse format.
\paragraph{{Core Optimisation Solver}}\label{core}  
As noted in Section \ref{sec:Speed-up}, the solution of the linear KKT system \eqref{eqn:kktn} is the most computationally expensive step in an IP algorithm. This step attributes also to the highest peak memory footprint, since the Newton-Raphson Jacobian \eqref{eqn:kktn} is the largest structure built\footnote{In general, both sparse matrix size ($m\times n$) and density have a direct relationship with the size of allocated memory in a computational program.}. The solution of the multi-period KKT structure of \eqref{eqn:kktn} through Schur-Complement breaks it into smaller $\mathbf{\Upsilon}_t$ blocks in \eqref{commonStruct}, and thus, significantly less memory allocations. In fact, the line with peak memory allocation is located in Alg. \ref{Alg1},L.\ref{save} where struct $\mathbf{inf}^\mathbf{\Upsilon}_t$ stores info to be called in Alg. \ref{Alg2} L.\ref{solve2}.} 
\section{Case Study and Results\footnote{Please note that, the efficiency of calculating analytical derivatives and their structures is illustrated in Appendix \ref{Appendix_D} of this paper}}
\begin{figure}[!htbp]
\centering
\includegraphics[width=3.5 in , height=1.5 in]{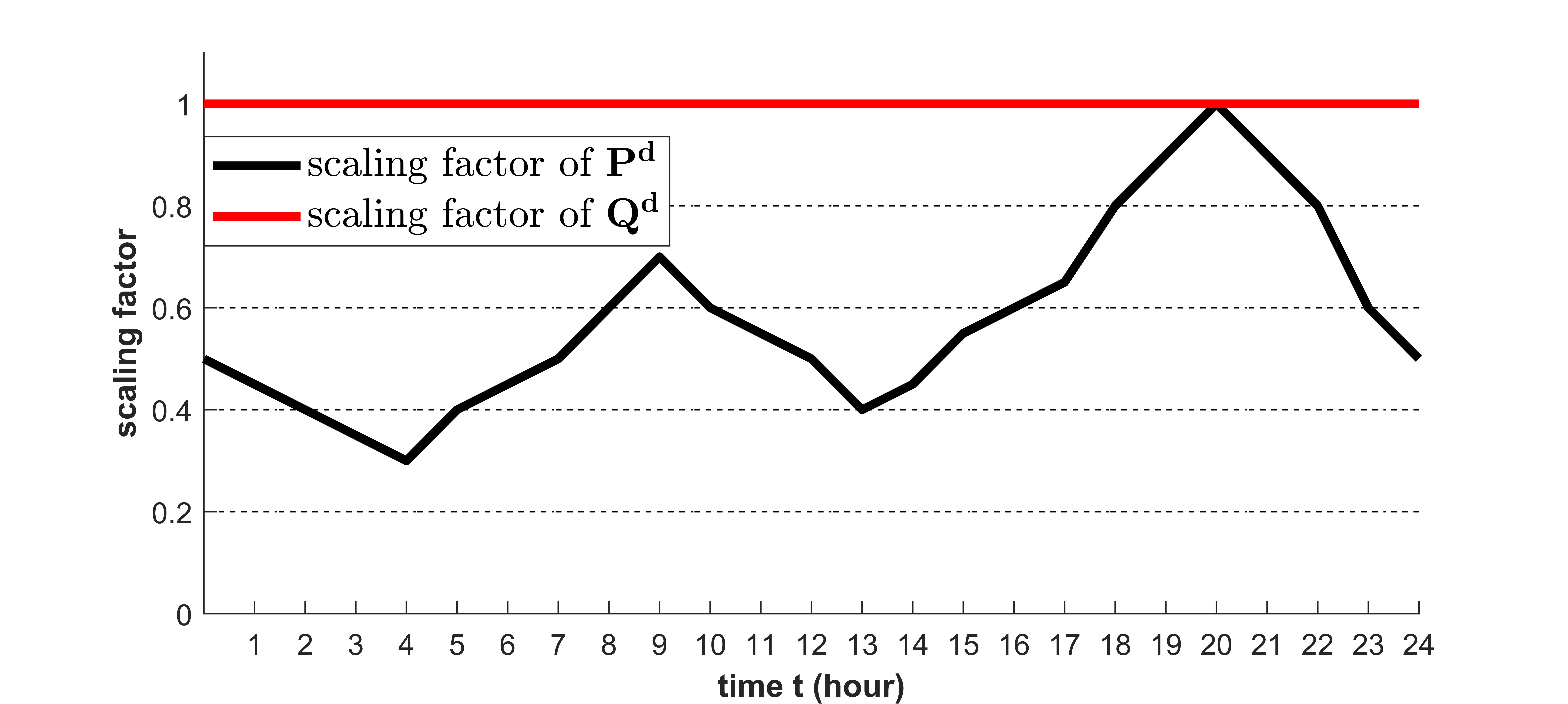}
\caption{Scaling factor multiplied by vector of consumption of active load $\boldsymbol{\mathcal{P}}^\mathrm{d}$ in order to simulate one period of 24 hours} 
\label{fig:scalingfactor}
\end{figure}

In this section, we present the results of benchmarking, obtained from implementations of different algorithms on similar platforms and workstations. The aim is to show the computational differences among mathematical algorithms when implemented on similar platforms.\\
 In this respect, standard case-files are adopted, Case9 \cite{schulz_long_1974}, IEEE30 \cite{alsac_optimal_1974}, IEEE118 \cite{noauthor_pg_tca118bus_nodate} and PEGASE1354 \cite{josz_ac_2016,fliscounakis_contingency_2013} are chosen for the study of SESS. Moreover, three distribution networks are considered for the simulation of EV: Case85 \cite{case85}, Case141 \cite{case141}, and a case study based on a real distribution grid in Mid-Norway \cite{zaferanlouei}; except for the last case study, the other cases can be found in the MATPOWER data folder \cite{zimmerman_matpower:_2011}. Details of these benchmarks are shown in Table \ref{tab:grid}. The transmission networks are considered in order to simulate the SESS and consequently the performance of static Schur-Complement structure that was discussed in Section \ref{sec:Speed-up}. The distribution networks are considered to simulate and compare the performance of dynamic Schur-Complement algorithm on arrival and departure of EV in one 24-hour period.\\ 
\begin{table}[htbp!]
\caption{Power grid Benchmark Models}
\label{tab:grid}
\begin{center}
\begin{tabular}{c c c c c} 
\toprule
 Type&case study&$n_b$ & $n_l$&$n_g$   \\  
\midrule \midrule
\multirow{4}{5em}{Transmission Grid}& Case9 & 9 & 9 & 4 \\
\cline{2-5}
&IEEE30 &30 & 41 & 6   \\
\cline{2-5}
&IEEE118 & 118 & 186 & 54 \\
\cline{2-5}
&PEGASE1354 & 1354 & 1991 & 260 \\
\midrule \midrule
\multirow{3}{5em}{Distribution Grid}&Case85 & 85 &84 & 1\\
\cline{2-5}
&Case141 & 141 &140&1\\
\cline{2-5}
&Mid-Norway &974 & 1023&2\\
\bottomrule 
\end{tabular}
\end{center}
\end{table}
For all simulation results, both transmission and distribution networks in the next two subsections, $\mathbf{BUS}$, $\mathbf{BRANCH}$,  $\mathbf{GEN}$, and  $\mathbf{GENCOST}$ matrices are taken from the original test-cases and are not modified. {Moreover, the flat initialisation strategy ($\frac{\mathbf{X}^{\mathrm{max}}-\mathbf{X}^{\mathrm{min}}}{2}$) is taken for all results presented in this paper.} Vector of bus active load is obtained through $\boldsymbol{\mathcal{P}}^\mathrm{d}=c^p(t).\boldsymbol{\mathcal{P}}$, and $c^p(t)$ is illustrated in Fig. \ref{fig:scalingfactor}, which fluctuates similar to a base load of households. Vector of bus reactive load is simulated as $\boldsymbol{\mathcal{Q}}^{\mathrm{d}} = c^q(t).\boldsymbol{\mathcal{Q}}$ with a constant scaling factor shown in Fig. \ref{fig:scalingfactor} by a red line. Note that the objective function of the optimisation problem in this paper includes only active power minimisation. $\boldsymbol{\mathcal{P}}$ and $\boldsymbol{\mathcal{Q}}$ are taken from the original values in $\mathbf{BUS}$ matrix.  {All simulations codes are developed in MATLAB environment.} They are performed on a computer with Intel(R) Xeon(R) CPU E5-2690 v4 @ 2.60 GHz and 384 GB RAM, and controlled with the single-thread environment to compare the computational differences in only the single-thread mode. 
\subsection{Transmission Network with Stationary Storage} \label{Trans}
\begin{figure}[!htbp]
\centering
\includegraphics[width=4 in , height=2 in]{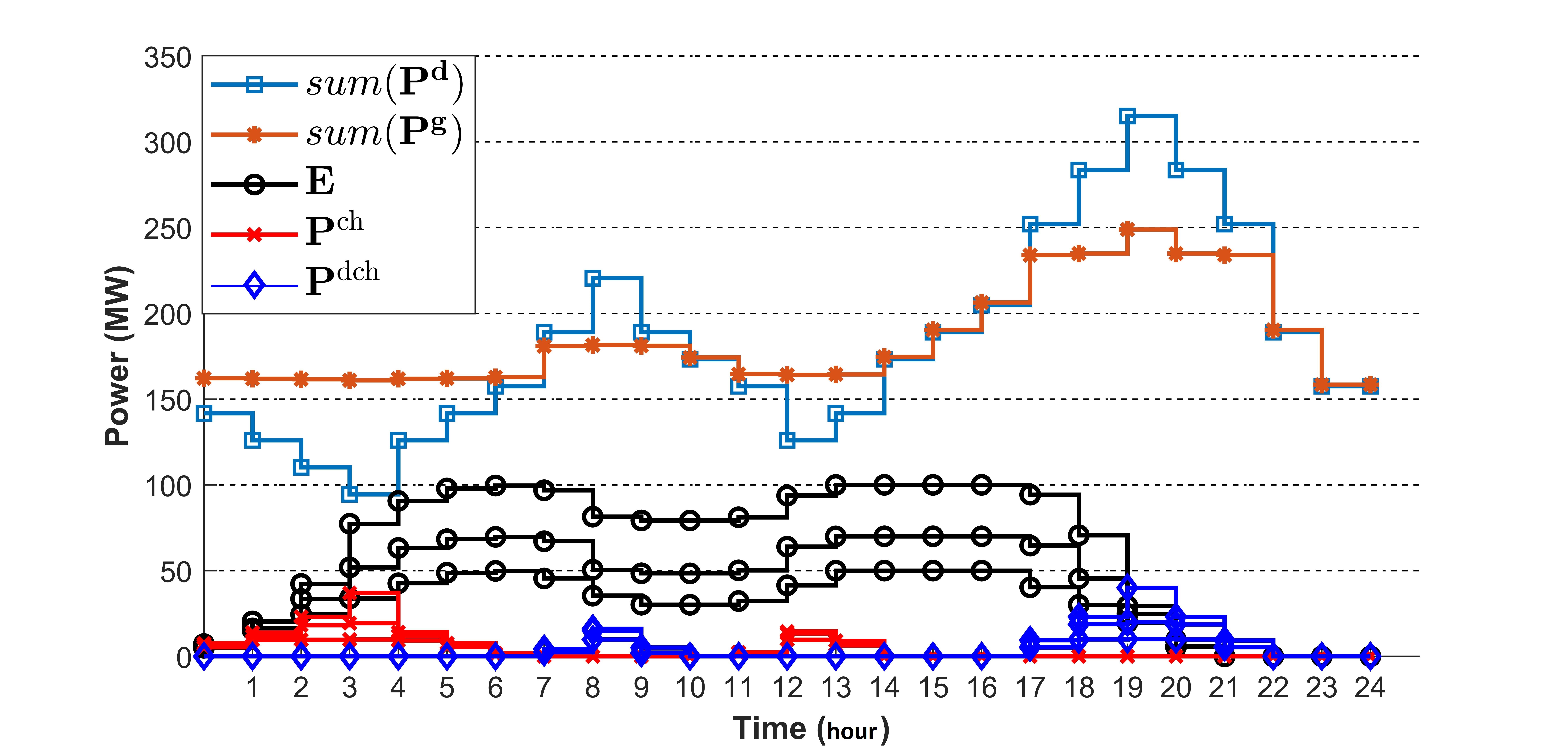}
\caption{IEEE Case9, $n_y=3$, $T=24$ and $n_g=3$. Total loads and generation vs operational variables of batteries---$\boldsymbol{\mathcal{SOC}}, \boldsymbol{\mathcal{P}}^{\mathrm{ch}}$ and $\boldsymbol{\mathcal{P}}^{\mathrm{dch}}$ } 
\label{fig:generaloverview}
\end{figure}

For transmission networks, the time-step is taken to be $\Delta t=1$ {hour}. The capacity of storage devices is $e_{i}^{\mathrm{max}}= 100 \ \mathrm{MWh}$, and consequently $e_{i}^{\mathrm{min}}= 0 \  \mathrm{MWh}$.  Charge and discharge limits are considered as $(\boldsymbol{\mathcal{P}}_t^\mathrm{ch})^\mathrm{min} = (\boldsymbol{\mathcal{P}}_t^\mathrm{dch})^\mathrm{min} = 0 \ \mathrm{MW}$ and $(\boldsymbol{\mathcal{P}}_t^\mathrm{ch})^\mathrm{max} = (\boldsymbol{\mathcal{P}}_t^\mathrm{dch})^\mathrm{max} = 10 \ \mathrm{MW}$. All charging and discharging efficiencies are taken as $\psi_{i}^\mathrm{ch}=0.95$ and $\psi_{i}^\mathrm{dch}=0.97$, respectively. Moreover, the initial status of storage devices is taken to be zero for all cases: $\mathbf{SOCi}= \begin{bmatrix}
\mathbb{0}
\end{bmatrix}_{n_y\times T}$.
Fig. \ref{fig:generaloverview} shows a typical optimisation outcome where the case study is Case9. Summation of loads and generations in each hour is depicted here with the operational strategy of energy storage systems with $n_y= 3$. Storage devices are located at buses 1, 2 and 3, where generators are located originally in the case-file. Since the objective function is a quadratic cost function, the storage devices are charging when the sum of loads is at the minimum and discharging when at the maximum. Fig. \ref{fig:generaloverview} is only for  illustration of the optimisation outcome and does not provide more insight.\\    
\subsubsection{Independency of Distribution of Storage Devices and Computational Performance}
Optimisation problems are solved for different distribution of storage devices at buses. Various types of scenarios are tested in order to verify that the distribution of storage devices does not have an impact on the overall computational time of each benchmark for convergence. In this respect, Table \ref{tab:dist} illustrates the iterations and overall time spent in solving each case study with specific strategies for the distribution of storage devices. In Table \ref{tab:dist}, First-Last strategy is for the one with distribution from first bus number to the $n_y^{\mathrm{th}}$ bus number when $n_y\leq n_b$. If $n_y > n_b$, then First-Last strategy repeats until all storage devices are located at the buses. For the Last-First strategy, the starting point is from the last bus gradually to the first bus. The Load-Bus strategy is for locating the storage devices at the buses that only have non-zero load in their original case-files: bus $i$ adopts one storage if $p_i^\mathrm{d}\neq 0$; this process repeats until all storage devices have been located. The last strategy is Fair-Dist, which uniformly distributes the storage devices among the buses. For instance if we have $n_b = 100$ and $n_y=10$ then every tenth bus adopts a storage device. Table \ref{tab:dist} shows that the time spent for each iteration, and for each scenario explained above, is very similar.\\ 
\begin{table}[htbp!]
\caption{Strategies on distribution of Storage devices}
\label{tab:dist}
\begin{threeparttable}
\begin{tabular}{c c C c c c c c c} 
\toprule
 Case&Distribution & \begin{tabular}[c]{@{}c@{}}$T$\\(period)\end{tabular}  &$n_y$ &\begin{tabular}[c]{@{}c@{}}No. of\\Iter.\end{tabular}   & Time (s)& $\frac{\mathrm{Time}}{\mathrm{iter}}$  \\  
\midrule \midrule
\tnote{A\textsuperscript{a}} & First-Last & 240 & 10 & 69 & 145 &  2.10
\\
\hline
\tnote{A\textsuperscript{a}} &Last-First & 240 & 10 & 90 & 190 &2.11  \\
\hline
\tnote{A\textsuperscript{a}} & Load-Bus & 240 & 10 & 67 & 141.9 & 2.11 \\
\hline
\tnote{A\textsuperscript{a}} & Fair-Dist & 240 & 10 & 81 & 170  &2.10\\
\hline
\tnote{B\textsuperscript{b}} & First-Last & 24 & 50 & 43 & 338 &7.86 \\
\hline
\tnote{B\textsuperscript{b}} & Last-First & 24 & 50 & 43 & 338 &7.86  \\
\hline
\tnote{B\textsuperscript{b}} & Load-Bus & 24 & 50 & 44 & 348&7.91  \\
\hline
\tnote{B\textsuperscript{b}} & Fair-Dist & 24 & 50 & 43 & 342 &7.95 \\
\bottomrule
\end{tabular}
\begin{tablenotes}
\item[a] A:{IEEE118}
\item[b] B:{PEGASE1354}
\end{tablenotes}
\end{threeparttable}
\end{table}
\subsubsection{Reported Time in the Benchmark}
\begin{figure}[!htbp]
\centering
\includegraphics[width=3.2 in , height=2.5 in]{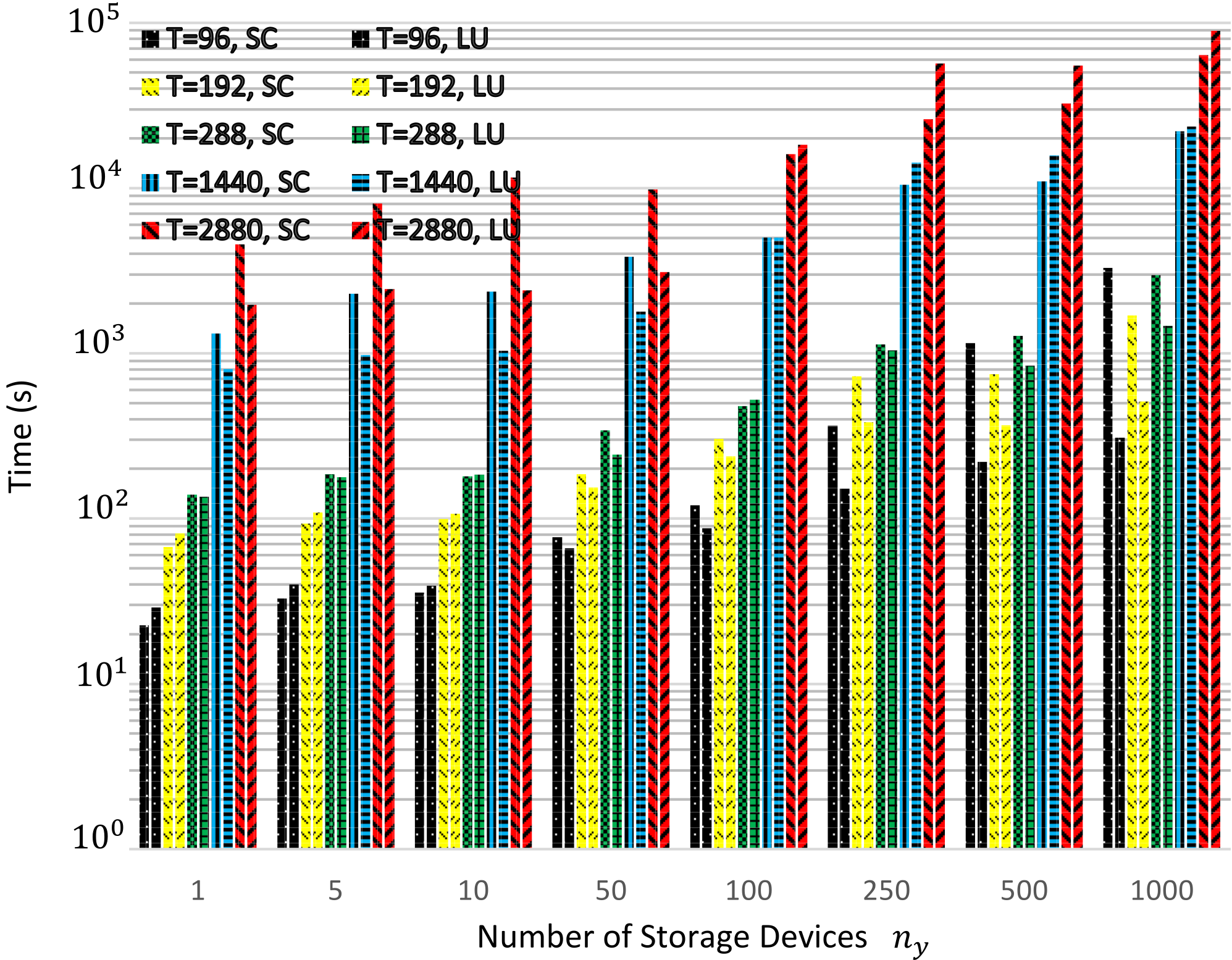}
\caption{Total time ($\mathrm{TotalTime= No_{\cdot} of \ Iter_{\cdot} \times TimePerIter}$) for solution of the linear KKT systems of \eqref{commonStruct} solved by Schur-Complement algorithm vs direct sparse LU solver, applied on Mid-Norway local Norwegian distribution grid} 
\label{Norwegian} 
\end{figure}

All the reported time values in this paper are shown in Appendix \ref{Appendix_D} (Table \ref{tab:numericVSanalytic}) and Figs. \ref{9bus}, \ref{30bus}, \ref{118bus}, \ref{1356bus}, \ref{85bus}, \ref{141bus} and \ref{Norwegian} show the total time taken for different benchmarks, which is the elapsed time to execute an algorithm until the optimum point is found. In other words, total time is $\mathrm{TotalTime= No_{\cdot} of \ Iter_{\cdot} \times TimePerIter}$. It should be kept in mind that the same number of iterations is considered for benchmarks reported and compared in Appendix \ref{Appendix_D} (Table \ref{tab:numericVSanalytic}) and Figs. \ref{9bus}, \ref{30bus}, \ref{118bus}, \ref{1356bus}, \ref{85bus}, \ref{141bus} and \ref{Norwegian}. Moreover, the selected strategy is First-Last distribution for all the above benchmarks. \\
Lastly, it should be noted that different distribution strategies have a different impact on the outcome of optimum operational values, such as objective function, generator scheduling, total losses, and voltage fluctuations in the grid.  
\subsubsection{Linear Algebra Overhead of {SESS} \texorpdfstring{$\mathbf{A}\mathbf{X}=\mathbf{B}$}{}}\label{AXB1}
\begin{figure}[!htbp]
\centering
\includegraphics[width=3.2 in , height=2.5 in]{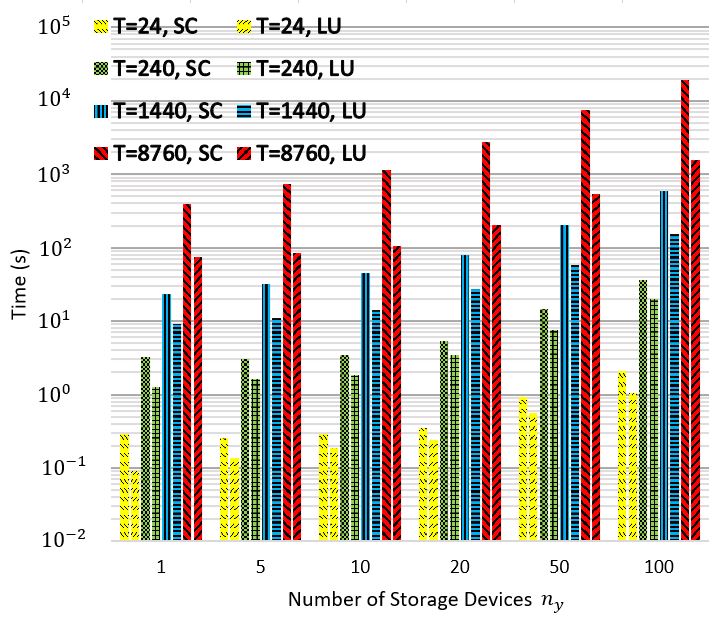}
\caption{Total time ($\mathrm{TotalTime= No_{\cdot} of \ Iter_{\cdot} \times TimePerIter}$) for solution of the linear KKT systems of \eqref{commonStruct} solved by Schur-Complement algorithm vs direct sparse LU solver, applied on Case9}
\label{9bus}
\end{figure}

\begin{figure}[!htbp]
\centering
\includegraphics[width=3.2 in , height=2.5 in]{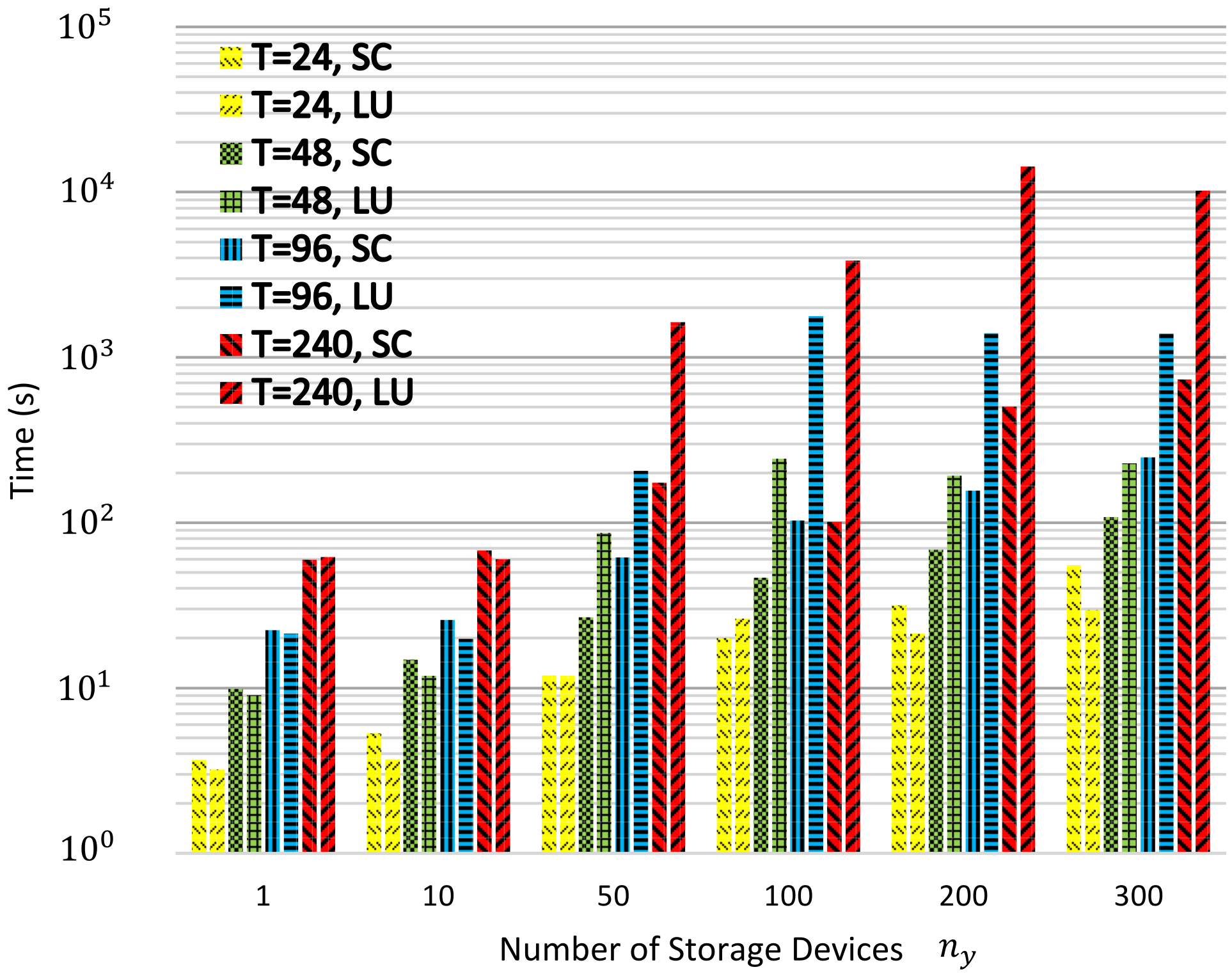}
\caption{Total time ($\mathrm{TotalTime= No_{\cdot} of \ Iter_{\cdot} \times TimePerIter}$) for solution of the linear KKT systems of \eqref{commonStruct} solved by Schur-Complement algorithm vs direct sparse LU solver, applied on IEEE118}
\label{118bus}
\end{figure}

\begin{figure}[!htbp]
\centering
\includegraphics[width=3.2 in , height=2.5 in]{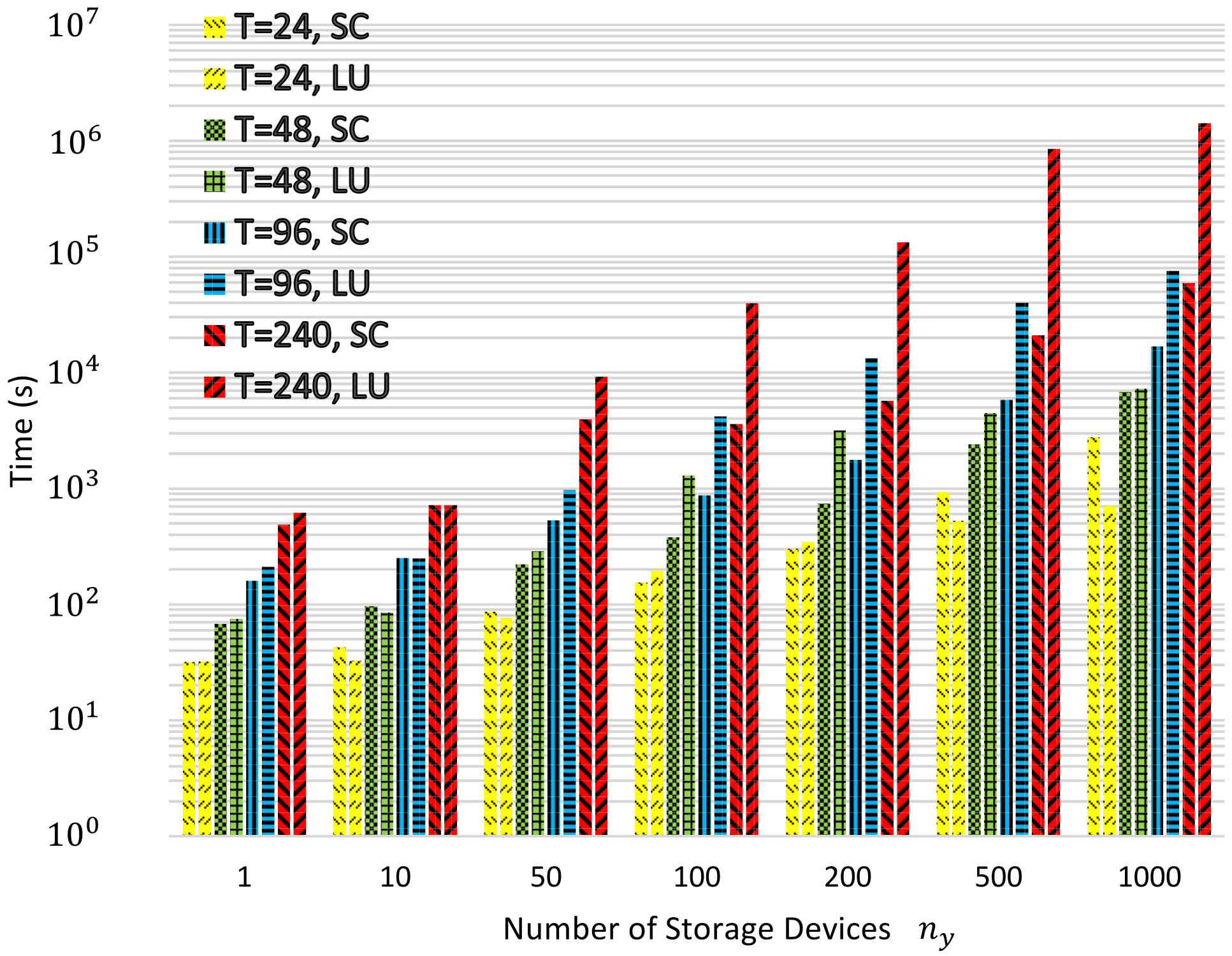}
\caption{Total time ($\mathrm{TotalTime= No_{\cdot} of \ Iter_{\cdot} \times TimePerIter}$) for solution of the linear KKT systems of \eqref{commonStruct} solved by Schur-Complement algorithm vs direct sparse LU solver, applied on PEGASE1354} 
\label{1356bus}
\end{figure}

As we stated in Section \ref{sec:Speed-up}, the most computationally expensive part of the IP algorithm is the solution of linear algebraic equations of the KKT system, which is similar to $\mathbf{A}\mathbf{X}=\mathbf{B}$, where $\mathbf{A}$ is the coefficient matrix, and $\mathbf{B}$ is the right-hand side vector. In order to assess the computational performance of the Schur-Complement technique, here we compare the performance of the tailored algorithm with direct sparse LU factorisation of complete structure of Eq. \eqref{commonStruct} (well-known arrowhead structure) and consequently the forward and backward solution, to present the computational time for the two algorithms implemented in the same platform. Note that most of the current IP-based solvers such as IPOPT, MIPS and KNITRO embed a direct solution method (LU/LDL). The main focus of the numerical results in this section is short-term horizon $T<240$ for IEEE30, IEEE118 and PEGASE1354 since full ACOPF would be more application oriented within short time-periods.\\
Fig. \ref{9bus} depicts the computational time needed to solve Eq. \eqref{commonStruct} of Case9 with Schur-Complement and a direct sparse-LU solver for number of storage devices $n_y= 1,5,10,20,50$ and $100$; each case takes into account time horizons of $T=24, 240, 1440$ and $8760$. Direct sparse-LU solver outperforms in all cases as $n_y$ and $T$ increase. This informs us that the Schur-Complement method is not efficient in comparison with a direct sparse-LU solver when the case study is a comparatively small network.\\ 
However, results for IEEE30, IEEE118 and PEGASE1354  suggest that the Schur-Complement method outperforms the direct sparse-LU solver when $T$ significantly increases. Since very large time-periods $T>240$ would not be considered as applied cases, we do not include them here. Moreover, for relatively small number of time periods $T \leq 24$ direct sparse-LU solver outperforms the Schur-Complement method in almost all cases \cite{kourounis_towards_2018}; therefore, we focus our results {for IEEE30, IEEE118 and PEGASE1354} when the $24<T<240$. It can be seen that when number of storage devices increases $n_y>10$, then the Schur-Complement method provides a more computationally efficient outcome.\\   
\begin{figure}[!htbp]
\centering
\includegraphics[width=3.2 in , height=2.5 in]{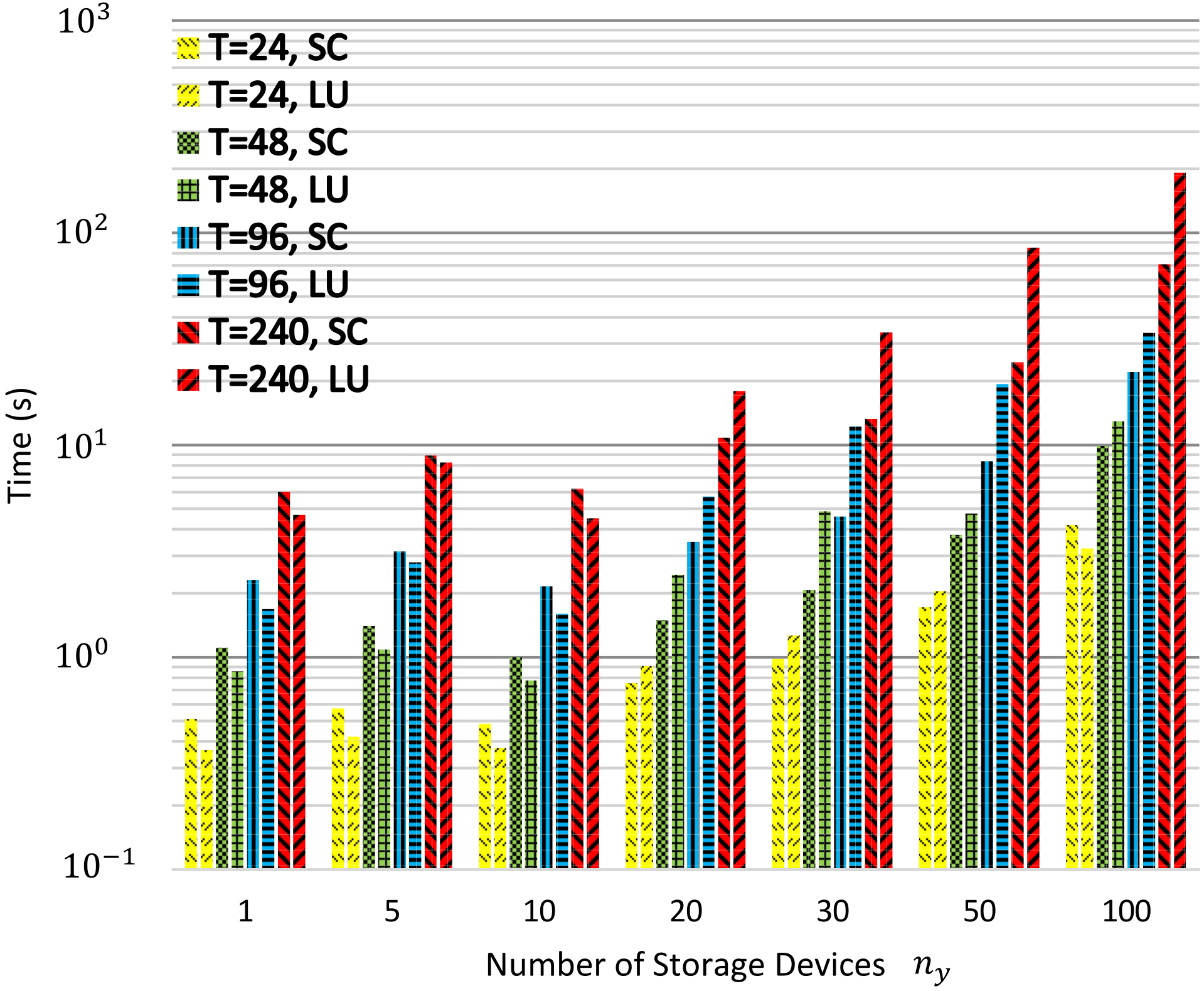}
\caption{Total time ($\mathrm{TotalTime= No_{\cdot} of \ Iter_{\cdot} \times TimePerIter}$) for solution of the linear KKT systems of \eqref{commonStruct} solved by Schur-Complement algorithm vs direct sparse LU solver, applied on IEEE30 } 
\label{30bus}
\end{figure}  
Fig. \ref{30bus} illustrates the comparative computational performance of the Schur-Complement solver and direct sparse-LU solver in order to solve the IEEE30 case study. As $n_y>10$ the Schur-Complement solver has higher performance which increases considerably when $T>24$. Note that the direct sparse-LU solver is dominant again when $n_y>300$, due to the reason that $n_y > \mathcal{K}$  where $\mathcal{K} \propto N_{\mathbf{\Upsilon}_t}=2n_b+2n_g+4n_y+n_{gn}+n_{{gl}_t}$, where $n_b$, $n_g$, $n_y$, $n_{gn}$ and $n_{{gl}_t}$ are respectively the number of buses, generators, storage devices, grid non-linear equalities and grid linear equalities at time $t$. Put simply, when $n_y>300$, the number of storage devices is larger than a certain number which is proportional to the size of blocks of $\mathbf{\Upsilon}_t$ in \eqref{commonStruct}. Therefore, the computationally demanding terms would be the calculation of Schur-Complement auxiliary blocks of $\mathbf{A}^a_t$ and $\mathbf{A}^b_t$ as, respectively, in terms of $\mathbf{S}_t=-\pmb{\rho}_t\mathbf{\Upsilon}_{t}^{-1}\pmb{\rho}_t^\top$ and $\mathbf{\Xi}_t= -\pmb{\rho}_t\mathbf{\Upsilon}_{t}^{-1}\pmb{\zeta}_t $ in Alg. \ref{Alg1}.\\ 
Numerical results of IEEE118 follow approximately a similar pattern as for IEEE30 when $n_y$ and $T$ increase, cf. Fig. \ref{118bus}. It can be observed that the difference between the direct sparse LU method and Schur-Complement method gets larger in IEEE118, which in turn, proves that the latter outperforms when size of $\mathbf{\Upsilon}_t$ blocks in \eqref{commonStruct} becomes larger.\\  
Fig. \ref{1356bus} shows the performance of the Schur-Complement method in comparison with the direct sparse LU solver where the case study is a large-scale optimisation of PEGASE1354. As expected, the Schur-Complement method outperforms the sparse-LU solver considerably when $n_y>10$ and $T>24$. \\ 
The simulation shows that computational time is highly dependent on the number of battery energy storage devices in the grid and the time-period. The Schur-Complement method is an efficient method to accelerate the MPOPF solution time when $n_y$ and $T$ are large numbers in optimisation. 
\subsection{{Memory Efficiency}}
\begin{figure}[!htbp]
\centering
\includegraphics[width=3.5 in , height=2.5 in]{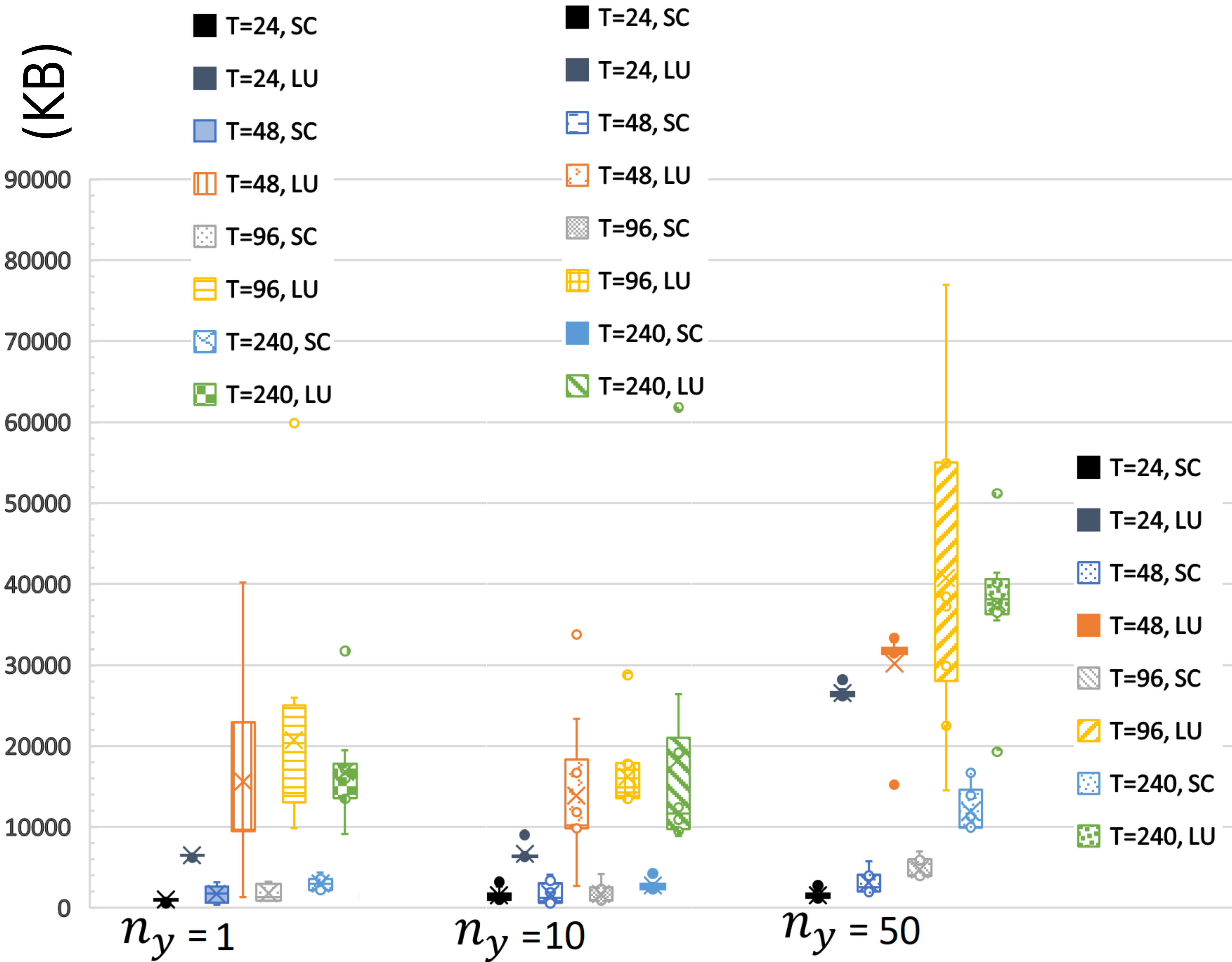}
\caption{Peak Memory ($KB$) for solution of the linear KKT systems of \eqref{commonStruct} solved by Schur-Complement algorithm vs direct sparse LU solver, IEEE 118, where $n_y=1$, $n_y=10$, and $n_y=50$.} 
\label{fig:Memory} 
\end{figure}

{In order to back the statement in subsection \ref{core}, the maximum memory consumption is tested. Therefore, the peak memory usage of Schur-Complement method is compared with that of sparse LU solver. Fig. \ref{fig:Memory} presents the peak memory consumption to solve IEEE118, where $n_y=1$, $n_y=10$, and $n_y=50$. In each $n_y$, four different time horizons $T=24$, $T=48$, $T=96$, and $T=240$ are tested. Moreover, each simulation test $\in \{1,..,24\}$ is repeated 10 times, and the results are shown as box plots in Fig. \ref{fig:Memory}. The results show that the average peak memory consumption of Schur-Complement method is more than 7 times less than sparse LU solver. It should be noted that increase in both $T$ and $n_y$ would result in higher maximum memory usage. Lastly, a larger variation of the peak memory usage is observed in the case of direct sparse-LU solver. It should be noted that LU solver is an internal MATLAB library.}\\  
\subsection{Distribution Network with EV Storage} \label{DIST}
In this section, three distribution networks are considered for the benchmarking study:  Case85 \cite{case85}, Case141 \cite{case141}, and a real Mid-Norway distribution grid.\\ References \cite{report_learning_2016,thingvad_economic_2019} report that the average driving distance is 52 km and as per reference \cite{bretteville-jensen_norwegian_2016}, the standard deviation is 22 km. EV fleet's arrival and departure are derived from the work hour lifestyle survey results presented in \cite{sterud_working_nodate}. A summary of the data for EV charge profile generation is provided in Table \ref{tab:EV}. Note that the departure times of the EV owners are calculated as 9.5 hours after their arrival time by considering 8 hours of working and 1.5 hours for total commuting time.\\ 
\begin{table}[htbp!]
\caption{Data for EV charge profile generation}
\label{tab:EV}
\begin{threeparttable}
\begin{tabular}{l l} 
\toprule
 Mean daily drive distance &52 km  \\  
\hline
Standard deviation of daily drive distance& 22 km\\
\hline
Standard deviation of daily drive distance distribution&10\% \tnote{1}  \\
\hline
Percentage of EV population that consume $\leq$ 18 kWh/100km &80\%\\
\hline
Percentage of EV population that consume $\geq$ 18 kWh/100km &20\%\\
\hline
Mean arrival time for the EV population& \begin{tabular}[c]{@{}c@{}}17:00\\hours\end{tabular}  \\
\hline
Standard deviation of arrival time for the EV population& 90 min\\
\hline
Standard deviation of daily arrival time for individual EV & 15 min\\
\hline
Percentage of 230V, 10A chargers & 70\%\\
\hline
Percentage of 230V, 16A chargers & 20\%\\
\hline
Percentage of 230V, 48A chargers & 10\%\\
\bottomrule
\end{tabular}
\begin{tablenotes}
\item[1] {of daily drive distance} 
\end{tablenotes}
\end{threeparttable}
\end{table}
The number of EV, departure time, and arrival time are selected as an input and solved for one period of 24 hours. {One full EV optimisation period is applied for the entire simulation from 12:00 PM until 12:00 PM the next day.} Time resolution in each profile can be seen in Table \ref{tab:EV_time}. Three types of charger power capacity are chosen and shown in Table \ref{tab:EV}. Discharge is considered to be inactive for all optimisation senarios: $\mathbf{CONDI}=\begin{bmatrix}
\mathbb{0}
\end{bmatrix}_{n_y\times T}$, input matrices of $\mathbf{AVBP}$ and $\mathbf{CONCH}$ are considered to be similar, i.e. $\mathbf{AVBP}=\mathbf{CONCH}$, which are derived from arrival and departure distribution functions. {The initial state of charge input $\mathbf{SOCi}$ is calculated according the distance each EV has traveled such that $[\mathbf{SOCi}_i]^{\mathrm{day_n}}={[e_{i}^{\mathrm{max}}-\mathrm{Energy(x)}]}^{\mathrm{day_{n-1}}}$, where $\mathrm{Energy(x)}$ is a function that calculates the energy roughly consumed.} Lastly, the departure state of charge is controlled through the last input matrix $\mathbf{SOCMi}$, which is the minimum state of charge, cf. the box constraint introduced as $\mathbf{SOCMi}_t\leq \boldsymbol{\mathcal{SOC}}_t \leq \boldsymbol{\mathcal{SOC}}^\mathrm{max} $, such that $\mathbf{SOCi}_{i,t}=\boldsymbol{\mathcal{SOC}}_i^\mathrm{max} \ \ \text{if }  \ \{\mathbf{AVBP}_{i,t}=1	\land \mathbf{AVBP}_{i,t+1}=0\} \lor \{\mathbf{AVBP}_{i,t=T}=1\}$.\\ 
\begin{table}[htbp!]
\caption{Time Resolution in different Optimisation sessions}
\label{tab:EV_time}
\begin{center}
\begin{tabular}{l l} 
\toprule
$T$ & $\Delta t $\\
\midrule \midrule
96 &15 min  \\  
\hline
192 & 7.5 min \\
\hline
288 &5 min  \\
\hline
1440 &1 min\\
\hline
2880 &30 sec\\
\bottomrule
\end{tabular}
\end{center}
\end{table} 
\subsubsection{Standard Distribution Cases}
 Case85 and Case141 are open-source radial distribution cases in the MATPOWER data folder; they are 11 kV and 12.5 kV medium voltage distribution grids, respectively. They are each fed by only one feeder, cf. Table \ref{tab:grid}. Linear cost functions are chosen for both of them such that $f(\boldsymbol{\mathcal{P}}^\mathrm{g})=\boldsymbol{\mu}^{\top}\boldsymbol{\mathcal{P}}^\mathrm{g}$ where $\boldsymbol{\mu} \in \mathbb{R}^{T\times 1}$ is a vector of marginal price from 12:00 PM until 12:00 PM the next day, randomly selected from the Nordpool \cite{nordpool} spot price, $\boldsymbol{\mathcal{P}}^\mathrm{g} \in \mathbb{R}^{T\times 1}$ is the power bought from the upstream network, and $T$ is the optimisation horizon. $\boldsymbol{\mathcal{P}}^\mathrm{d}$ and $\boldsymbol{\mathcal{Q}}^\mathrm{d}$ of Case85 and Case141 are calculated similar to the procedure in Section \ref{Trans}. Arrival and departure timetables of EV are calculated according to the description in Section \ref{DIST}. 
\subsubsection{Local Mid-Norway Distribution Grid}
\begin{figure}[!htbp]
\centering
\includegraphics[width=3.5 in , height=2 in]{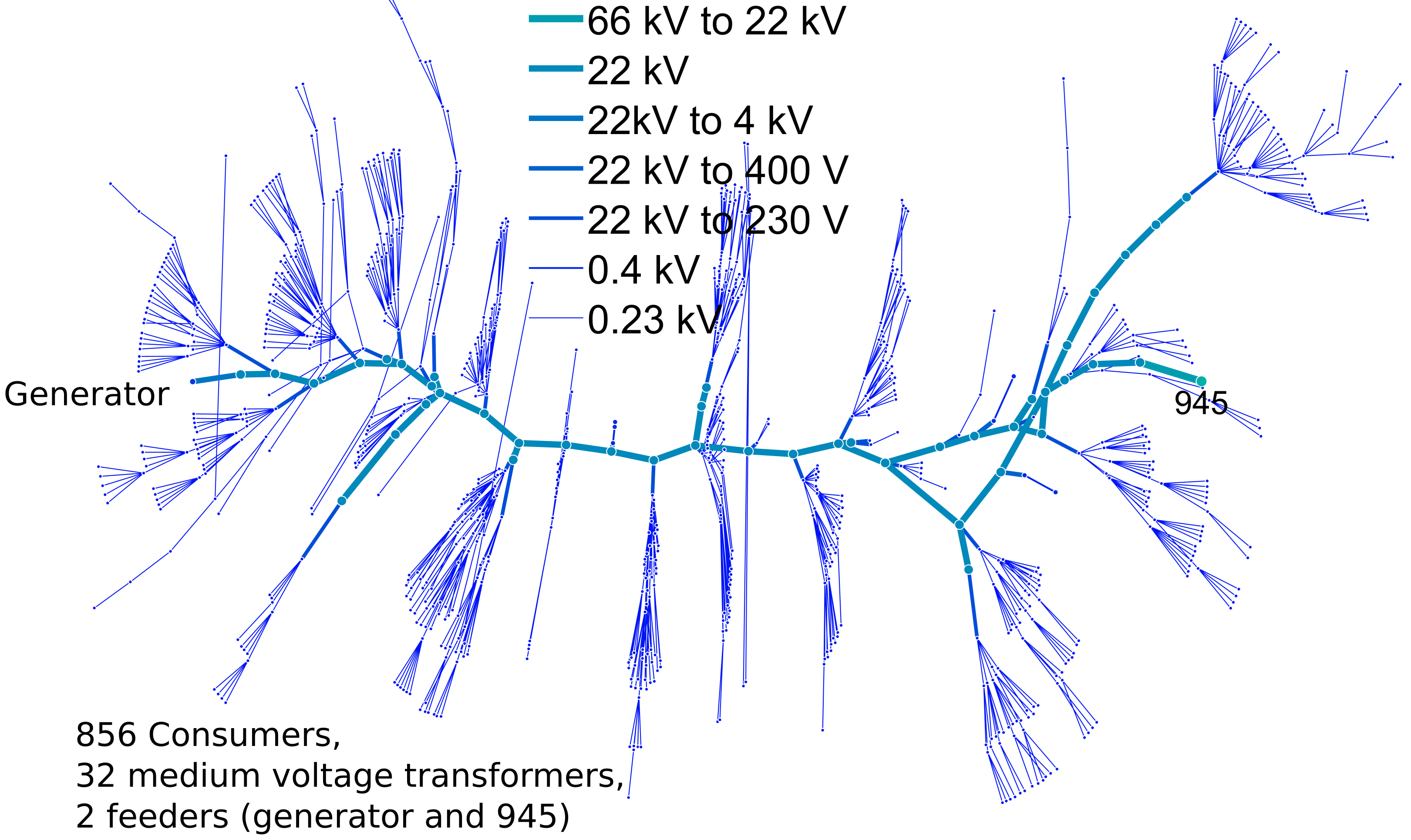}
\caption{Local distribution grid located in Norway with 856 costumers. using the visulaisation technique from \cite{cuffe_visualizing_2017}}
\label{Norwegian1}
\end{figure}

\begin{figure}[!htbp]
\centering
\includegraphics[width=3 in , height=6 in]{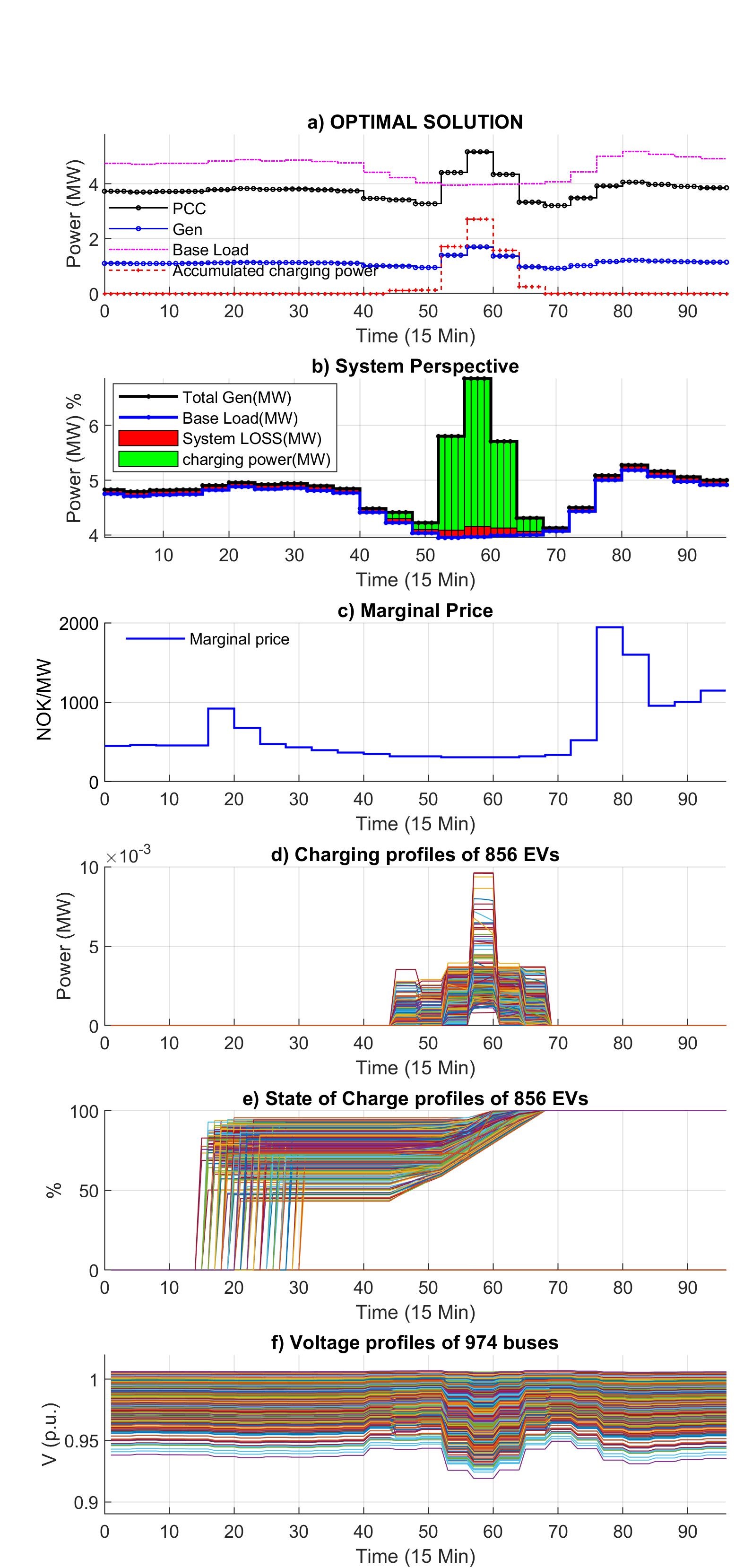}
\caption{Outcome of optimisation of large local distribution grid in mid-Norway, Data for 12:00 PM Feb 1, 2012, until 12:00 PM Feb 2, 2012, with highest pick of electricity price in the year of 2012: a) General perspective of optimisation, total hourly consumption profile, shown as Base Load and optimal production profile of PCC and generator, plus accumulated charging power of 856 EV, 1 EV per costumer. b) accumulation of total generation vs base load and in between two curves, losses in red and charging power in green c) hourly spot price, 8:00 am of Feb 2, 2012, is highest price of the year 2012 d) charging profile of 856 EVs. Outcome of optimisation suggests charging times and power values such that to compromise between total cost and total loss. e) state of charge of 856 EV, f) voltage fluctuations of 974 buses.} 
\label{Norwegian2}
\end{figure}

The real-case Mid-Norway distribution grid is a large distribution case study, shown in Fig. \ref{Norwegian1}, which is a 22 kV medium voltage  to 230 low voltage grid, fed by: 1) a high-voltage 66 kV feeder, Point of Common Coupling (PCC), shown as red circle dot, and 2) a local generator, shown as light blue circle dot. Cost functions of PCC and generator are similar and are a linear function of $f(\boldsymbol{\mathcal{P}}^{\mathrm{g}^{\mathrm{PCC}}},\boldsymbol{\mathcal{P}}^{\mathrm{g}^{\mathrm{gen}}})=\boldsymbol{\mu}^{\top}(\boldsymbol{\mathcal{P}}^{\mathrm{g}^{\mathrm{PCC}}}+ \boldsymbol{\mathcal{P}}^{\mathrm{g}^{\mathrm{gen}}})$ where $\boldsymbol{\mu} \in \mathbb{R}^{T\times 1}$ is the marginal hourly spot price (NOK/MW).  We assumed that the feeder and generator have similar hourly cost functions. The network has 32 MV/LV transformers (shown as dark blue squares in Fig. \ref{Norwegian1}) and feeds 856 registered consumers in the low-voltage area. $\boldsymbol{\mathcal{P}}^\mathrm{d}$ and $\boldsymbol{\mathcal{Q}}^\mathrm{d}$ for Mid-Norway grid are acquired from the local DSO and are hourly real consumption data of 856 consumers.\\ 
{Fig. \ref{Norwegian2} depicts the optimisation outcome for an EV period of 12:00 PM Feb 1, 2012, until 12:00 PM Feb 2, 2012, where optimisation resolution is 15 minutes, and thus $T = 96$}. Fig. \ref{Norwegian2} a) shows the overall picture of a day with the base load of consumers, and optimal production from the generator and PCC as well as optimal charging. Fig. \ref{Norwegian2} b) provides more insight into the optimisation results, where the optimum is a middle ground between: 1) placing all EV chargings with highest capacity of charge for the lowest price (shown in Fig. \ref{Norwegian2} c)), and 2) minimising losses at the same time. The compromise outcome is a pick of charging EV distributed between the time index of 45 until 68 and in a relatively sharp manner. It should be kept in mind that if the same optimisation formulation is run with DCOPF, then the charging pick would be sharp (all in the lowest price timestep), such that the loss would not be seen. Fig. \ref{Norwegian2} d) and e) illustrate the charge and state of charge profile of each EV owner, which total 856. Lastly, Fig. \ref{Norwegian2} f) shows the voltage variation of 974 buses over 24 hours and 96 timesteps. \\     
\subsubsection{Linear Algebra Overhead of EV \texorpdfstring{$\mathbf{A}\mathbf{X}=\mathbf{B}$}{}} \label{AXB2}
\begin{figure}[!htbp]
\centering
\includegraphics[width=3.2 in , height=2.5 in]{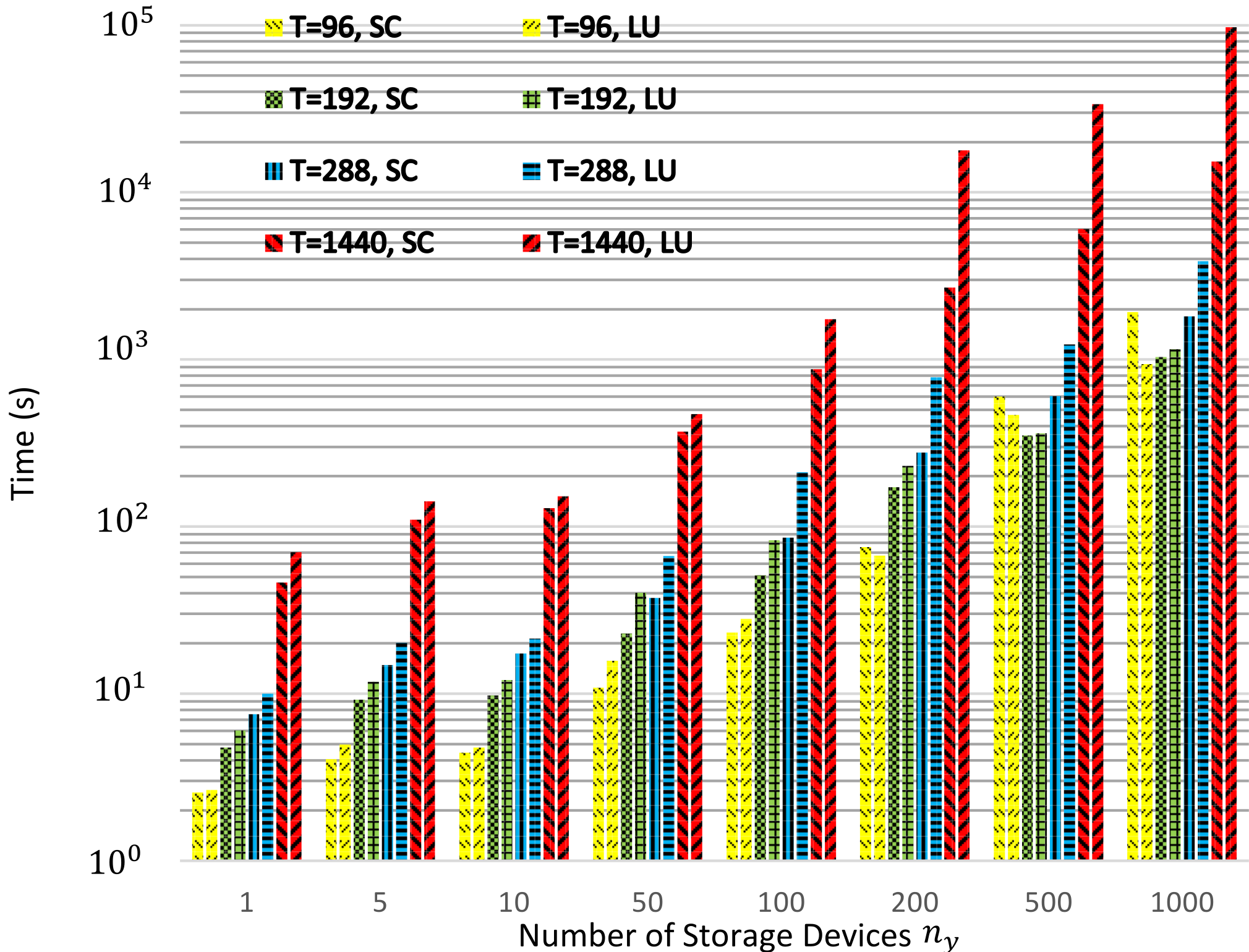}
\caption{Total time ($\mathrm{TotalTime= No_{\cdot} of \ Iter_{\cdot} \times TimePerIter}$) for solution of the linear KKT systems of \eqref{commonStruct} solved by Schur-Complement algorithm vs direct sparse LU solver, applied on CASE85} 
\label{85bus} 
\end{figure}

\begin{figure}[!htbp]
\centering
\includegraphics[width=3.2 in , height=2.5 in]{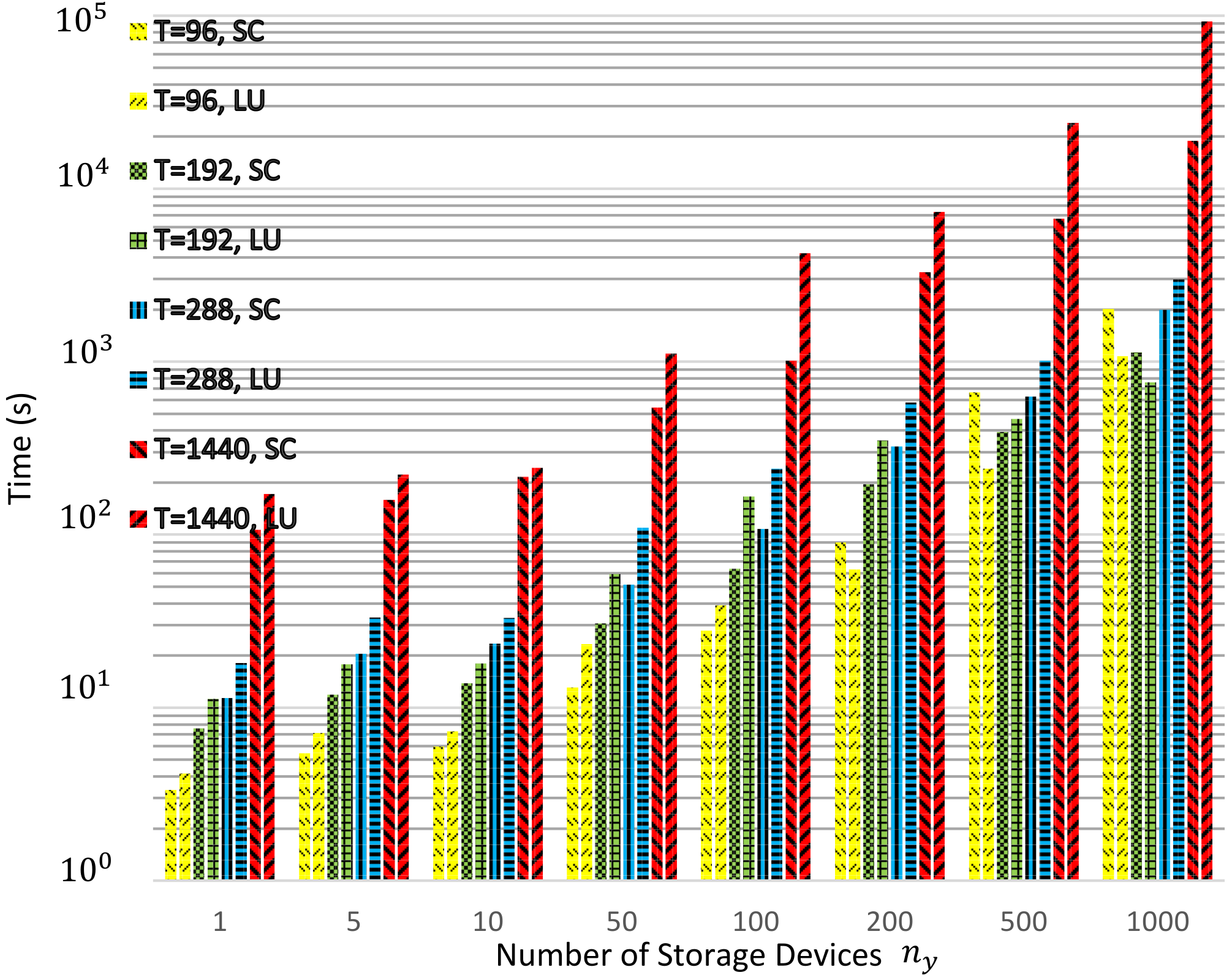}
\caption{Total time ($\mathrm{TotalTime= No_{\cdot} of \ Iter_{\cdot} \times TimePerIter}$) for solution of the linear KKT systems of \eqref{commonStruct} solved by Schur-Complement algorithm vs direct sparse LU solver, applied on CASE141} 
\label{141bus}
\end{figure}

Similar to Section \ref{AXB1}, the computational performance of Schur-Complement algorithms of Algorithm 1 and Algorithm 2 are compared with that of a direct sparse LU solver. Despite the fact that the Schur-Complement structures of SESS and EV are different, the main difference could be that EV have lower number of coupling constraints and some simulation hours could be completely decoupled; it is more efficient to solve them separately from coupled times. Here in this paper, we solve one EV period completely, both with the Schur-Complement solver and the direct sparse LU solver.\\
Figs. \ref{85bus} and \ref{141bus} show the computational time to solve a similar structure of \eqref{commonStruct} with the Schur-Complement algorithm vs the direct sparse LU solver, where number of EVs $n_y$ are increased from 1 to 1000 in the benchmark case study of Case85, with the strategy of First-Last as shown in Table \ref{tab:dist}. The results follow the same pattern as that of optimisation of SESS in IEEE118 shown in fig. \ref{118bus}. \\
In small time horizons when $n_y$ increases the efficiency of computing of Schur-Complement algorithm surpass that of the direct sparse LU solver until a certain point, and then decreases with a slope again. This is more evident when $T=96$ and $n_y>100$. Furthermore, it can be seen $T = 192$ and $n_y>200$. \\ 
The Mid-Norway case study is an interesting case where the direct LU solver mostly performs better. This is due to the fact the $T$ times LU factorisation of the block $\mathbf{\Upsilon}_t$ in Alg. 1 of the Schur-Complement algorithm is more computationally expensive than the solution of the entire coefficient matrix of \eqref{commonStruct}. This is due to the fact that the number of coupled blocks of $\mathbf{\Upsilon}_t$ is reduced, which in turn is because of input matrices, defined by the arrival and departure of EV. In fact the ratio of coupled blocks can be calculated as $\frac{[ \textrm{average of dep time steps}- \textrm{average of arrival time steps}]}{\textrm{total time steps}} =0.52$ which shows that only 52\% of blocks are coupled. One more interesting aspect is the clear observation of the turning point, when $n_y=100$ for two time horizons of $T=1440$ and $T=2880$, such that when $n_y<100$ the direct sparse LU solver performs better, and when $n_y>100$ then the Schur-Complement is supreme. \\ 

\section{Conclusion and Future Work}
A high performance and memory-efficient multi-period ACOPF solver based on a primal-dual IP method is proposed in this paper. In order to boost the computational performance, two mathematical approaches have been investigated. Partial derivatives of linear and non-linear constraints, objective function, and KKT conditions have been extracted analytically and consequently their sparse structures have been explored and exploited. A tailored algorithm has been suggested, using a new re-ordering format, in order to solve the sparse multi-period structure of Newton step in the IP method, with high computational performance. From the numerical results, the performance of the proposed Schur-Complement method is compared with a general sparse LU solver.  Numerical results suggest that a tailored {Schur-Complement} algorithm could be computationally supreme in a problem with certain specifications, such as (1) large networks (large number of bus and branches) (2) {different} optimisation horizon ($T$), and (3) large number of storage devices. In future works, we propose a parallelised Schur-Complement algorithm and benchmark it thoroughly.

\section*{Acknowledgment}
The authors acknowledge the contributions of Iver Bakken Sperstad and Venkatachalam Lakshmanan, researchers at SINTEF Energy Research, Norway, and Jamshid Aghaei, a former postdoctoral researcher at NTNU.

\bibliographystyle{IEEEtran}
\bibliography{IEEEabrv,DOPF}
\begin{IEEEbiography}[{\includegraphics[width=1in,height=1.25in,clip,keepaspectratio]{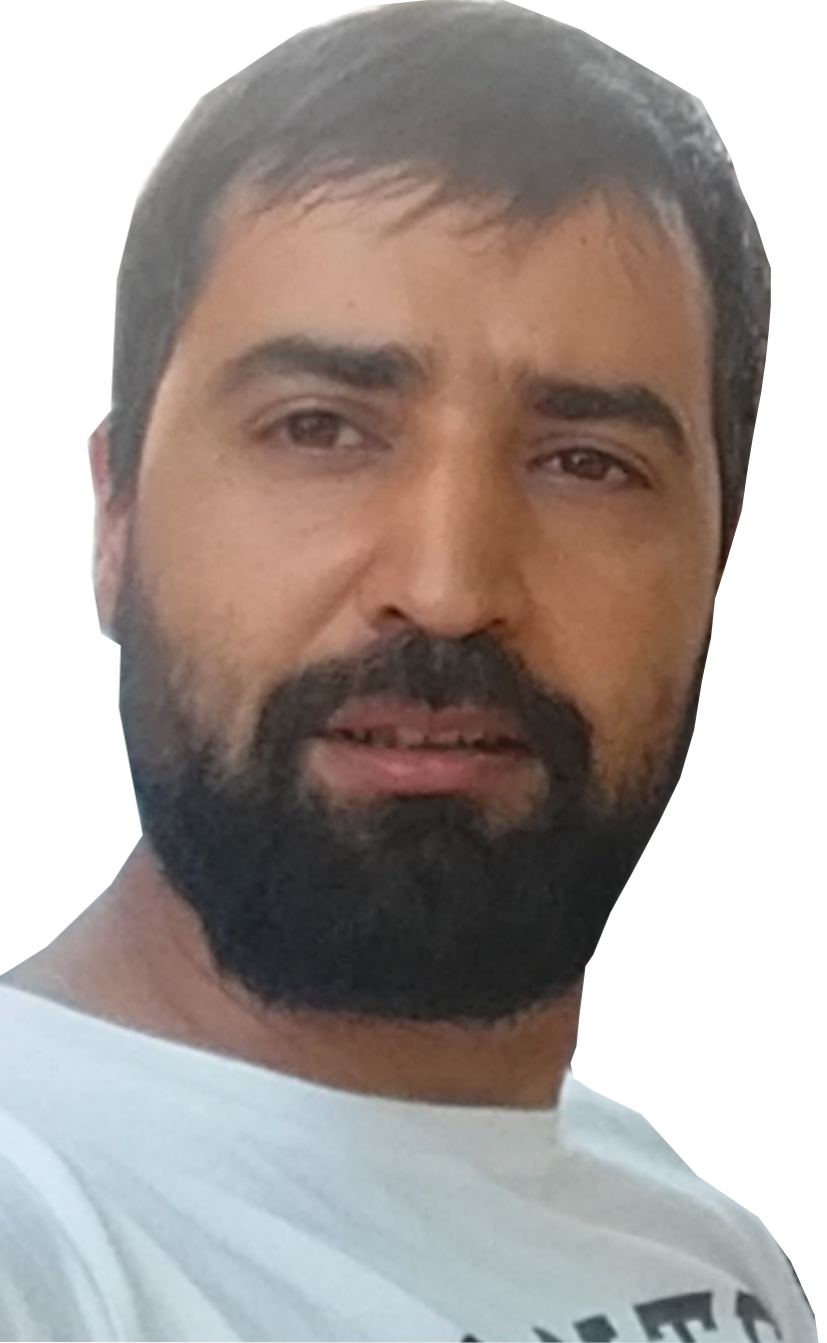}}]%
{Salman Zaferanlouei} (S'15) received his M.Sc. from the department of Energy Engineering and Physics, Amirkabir University of Technology (Tehran Polytechnic), Tehran, Iran, in 2009. He received his Ph.D. from the Department of Electric Power Engineering (IEL), Norwegian University of Science and Technology (NTNU), Trondheim, Norway, in 2020 on the topic of ``Integration of Electric Vehicles into Power Distribution Systems – The Norwegian Case Study; Using High-Performance Multi-Period AC Optimal Power Flow Solver". He is currently a researcher at IEL, NTNU. His research interests are core optimisation problems, high performance computing and power system simulations, economics and modelling. 
\end{IEEEbiography}
\begin{IEEEbiography}[{\includegraphics[width=1in,height=1.25in,clip,keepaspectratio]{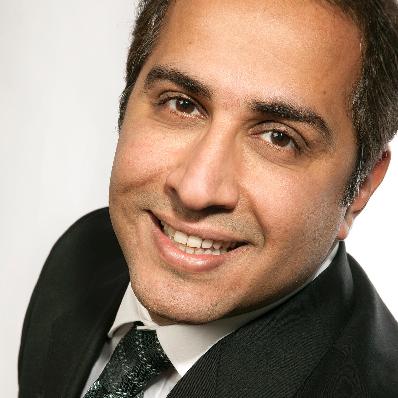}}]%
{Hossein Farahmand}
Farahmand (S'10, M'13, SM'20) received his Ph.D. degree from the Norwegian University of Science and Technology (NTNU), Trondheim, Norway, in 2012. He is currently working as Associate Professor at the Department of Electric Power Engineering, NTNU, and is a member of the Electricity Markets and Energy System Planning research group. His research interests include power system balancing, power market analysis, demand-side management, and flexibility operation in distribution systems. He has been involved in several EU-projects including INVADE H2020, EU FP7 TWENTIES, EU FP7 eHighway2050, and IRPWIND. He is the leader of the research task on the evaluation of the value of flexibility in smart grids in the Centre for Intelligent Electricity Distribution (CINELDI). 

\end{IEEEbiography}
\begin{IEEEbiography}[{\includegraphics[width=1in,height=1.25in,clip,keepaspectratio]{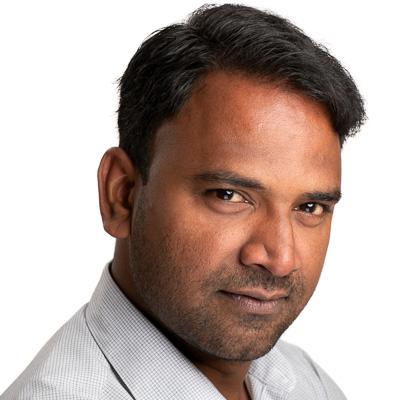}}]%
{Vijay Venu Vadlamudi}
 (S’08, M’11) received his Ph.D. degree from the Indian Institute of Technology Bombay, India, in 2011. He is currently working as Associate Professor at the Department of Electric Power Engineering, Norwegian University of Science and Technology (NTNU), Norway. He has also been the Deputy Head of Department (Education) since 2017, and is a member of the Power System Operation and Analysis research group. His areas of research interest include reliability and risk – based power system operation and planning practices, and probabilistic methods applied to power system analysis. He is on the Editorial Board of the journal IET Generation, Transmission \& Distribution, serving as Subject Editor for the field of power system reliability. 

\end{IEEEbiography}
\begin{IEEEbiography}[{\includegraphics[width=1in,height=1.25in,clip,keepaspectratio]{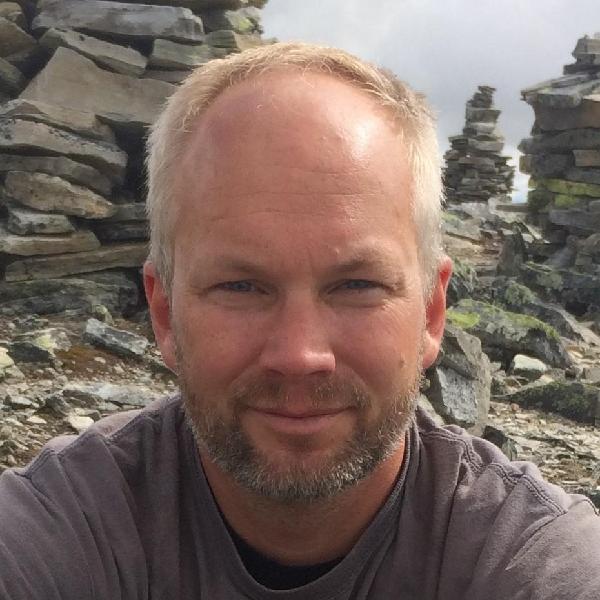}}]%
{Magnus Korpås} 
 (M’15) received his Ph.D. degree from the Norwegian University of Science and Technology (NTNU), Norway, in 2004 on the topic of optimizing the use of energy storage for distributed wind energy in the power market. He is currently working as Professor at the Department of Electric Power Engineering, NTNU, where he also leads the Electricity Markets and Energy System Planning research group. He is a leader and active participant in several large energy research projects at national and European levels. He is a former Research Director of the Department of Energy Systems at SINTEF Energy Research, Norway. He was a visiting researcher in the MIT Laboratory for Information \& Decision Systems (LIDS) in 2018-2019.  He is also the leader of the scientific committee and the leader of the work package on flexible resources in the power system in the Centre for Intelligent Electricity Distribution (CINELDI).
\end{IEEEbiography}
\begin{appendices}
	\section{BATTPOWER Input} \label{Appendix_A}
{BATTPOWER input matrices are introduced and elaborated in this section. The main introduced input is matrix of $\mathbf{BATT}$, which represents a connection matrix of $n_y \in \mathbb{N}$ energy storage devices. It contains charge and discharge rates and efficiencies of storage devices along with their maximum and minimum energy capacities. Moreover, it includes initial points of the charge, discharge, and state of charge variables. Table \ref{tab:inputMatrix} summarises the input matrices fed into the BATTPOWER solver. Note that $T \in \mathbb{N}$ is the time period is optimisation and $t$ is a time in the interval of $t \in \{1,...,T\}$.}
\begin{table*}[!htbp]
 \caption{Definition of Input Matrices }
\label{tab:inputMatrix}
\begin{threeparttable}
\begin{tabularx}{\linewidth}{ l c c X }
\toprule
&   \multicolumn{2}{c}{Size of Matrix}& \\
\cmidrule{2-3} 
Input     & $n$&  $m$ &Description \\ 
\midrule
$\mathbf{BUS}$     &$n_b$ &\tnote{1}& Examples can be found in \cite{zimmerman_matpower:_2011} \\ \midrule
$\mathbf{BRANCH}$  &$n_l$ &\tnote{1}&Examples can be found in \cite{zimmerman_matpower:_2011}\\ \midrule
$\mathbf{GEN}$     &$n_g$ &\tnote{1}&Examples can be found in \cite{zimmerman_matpower:_2011} \\ \midrule
$\mathbf{GENCOST}$ &$n_g$ &\tnote{1}&Examples can be found in \cite{zimmerman_matpower:_2011}\\ \midrule
$\mathbf{BATT}$    &$n_y$  &\tnote{1}&\ \texttt{BATT\_BUS}, \ \texttt{SOC\_OPT}, \ \texttt{PCH\_OPT}, \ \texttt{PDICH\_OPT}, \ \texttt{Q\_INJ\_OPT}, \  $\boldsymbol{\mathcal{SOC}}^\mathrm{max}$, \ $\boldsymbol{\mathcal{SOC}}^\mathrm{min}$, \ $(\boldsymbol{\mathcal{Q}}^\mathrm{s})^\mathrm{max}$, \ $(\boldsymbol{\mathcal{Q}}^\mathrm{s})^\mathrm{min}$, \ \texttt{MBASE}, \ $(\boldsymbol{\mathcal{P}}^\mathrm{ch})^\mathrm{max}$, \ $(\boldsymbol{\mathcal{P}}^\mathrm{dch})^\mathrm{max}$ \ \texttt{EFF\_CH} ($\mathbf{\Psi}^\mathrm{ch}$) \ \texttt{EFF\_DICH} ($\mathbf{\Psi}^\mathrm{dch}$) \\ 
\midrule
$\mathbf{AVBP}$   &$n_y$ &$T$&$\mathbf{AVBP} \in \mathbb{B}^{n_y \times T}$ ($\mathbb{B}$ is a binary set.)\tnote{2} which is the availability matrix of active power provision of storage devices, such that  $\mathbf{AVBP}_{i,t}=1$ if the $i^{th}$ storage at $t^{th}$ time is available and connected to the grid, otherwise $\mathbf{AVBP}_{i,t}=0$, where $T$ is the optimisation horizon. \\ 
\midrule
$\mathbf{CONCH}$  &$n_y$ &$T$&$\mathbf{CONCH} \in \mathbb{B}^{n_y \times T}$ is the charge connectivity matrix in which $\mathbf{CONCH}_{i,t}=1$ if the $i^{th}$ storage at $t^{th}$ time has a charging option, otherwise $\mathbf{CONCH}_{i,t}=0$.  \\ \midrule
$\mathbf{CONDI}$  &$n_y$ &$T$&$\mathbf{CONDI} \in \mathbb{B}^{n_y \times T}$ is the discharge connectivity matrix such that $\mathbf{CONDI}_{i,t}=1$ if the $i^{th}$ storage at $t^{th}$ time has the available discharging option, otherwise $\mathbf{CONDI}_{i,t}=0$\tnote{3}.\\ \midrule
$\mathbf{AVBQ}$   &$n_y$ &$T$&$\mathbf{AVBQ} \in \mathbb{B}^{n_y \times T}$ is the availability matrix of reactive power provision of storage devices such that $\mathbf{AVBQ}_{i,t}=1$ if the $i^{th}$ storage at $t^{th}$ time has the available option for reactive power provision, otherwise $\mathbf{AVBQ}_{i,t}=0$. \\\midrule
$\mathbf{AVG}$    &$n_g$ &$T$&$\mathbf{AVG} \in \mathbb{B}^{n_g \times T}$ which is the availability matrix of generators within the optimisation time horizon and consequently $\mathbf{AVG}_{i,t}=1$ if the $i^{th}$ generator at $t^{th}$ time is available to inject power in the grid. \\\midrule
$\mathbf{SOCi}$   &$n_y$ &$T$&$\mathbf{SOCi} \in \mathbb{R}^{n_y \times T}$ is the matrix consisting of initial state of charge of $n_y$ storage devices over time $t \in \{1,...,T\}$.  A value for initial state of charge $\{0\leq \mathbf{SOCi}_{i,t}\leq 1\}$ is allocated for the $i^{th}$ storage device at time $t$ if and only if one of these conditions is satisfied: 1) $\mathbf{AVBP}_{i,t=1}=1$. 2) $\mathbf{AVBP}_{i,t-1}=0$ and $\mathbf{AVBP}_{i,t}=1$ (arrival definition), otherwise $\mathbf{SOCi}_{i,t}=0$.  \\ \midrule
$\mathbf{SOCMi}$  &$n_y$ &$T$&$\mathbf{SOCMi} \in \mathbb{R}^{n_y \times T}$ matrix which includes the minimum state of charge of $n_y$ storage devices through time $t \in \{1,...,T\}$. The state of charge of the $i^{th}$ storage device at the departure time of $t$ can be settled if one of these two conditions is satisfied: 1) $\mathbf{AVBP}_{i,t}=1$, $\mathbf{AVBP}_{i,t+1}=0$. 2) $\mathbf{AVBP}_{i,t=T}=1$. \\ \midrule
$\mathbf{PD}$ & $n_b$& $T$& Time series of active loads. \\ \midrule
$\mathbf{QD}$ & $n_b$& $T$& Time series of reactive loads.  \\ 
\bottomrule
\end{tabularx}
\begin{tablenotes}
\item[1] {{User Defined.}}
\item[2]{$\mathbb{B}$ is a binary set.}
\item[3] {{Note that  $\mathbf{AVBP}_{i,t}=0$  means that the $i^{th}$ storage\textbackslash EV at time $t$ is not available; therefore, the same element in charge and discharge connectivity matrices must be zero: $\mathbf{CONCH}_{i,t}=0$ and $\mathbf{CONDI}_{i,t}=0$. The converse logic is not valid.}}
\end{tablenotes}
\end{threeparttable}
\end{table*}

	\cleardoublepage
	\section{Extraction of first and second partial derivatives} \label{Appendix_B}
In general, if we assume a complex scalar function $ f:{\rm I\!R}^n\rightarrow\mathbb{C}$ of a real vector such as \eqref{X}, the first derivative can be calculated as:
\begin{subequations}
\begin{alignat}{2}
\label{first-derivative-objective-a}
&f_{\mathbf{X}}=\frac{\partial f}{\partial \mathbf{X}}= \bigg[\frac{\partial f}{\partial \mathbf{x}_1} \quad  \frac{\partial f}{\partial \mathbf{x}_2}\quad ...\quad \frac{\partial f}{\underbrace{\partial \mathbf{x}_t}_{\Downarrow }}
 ...  \frac{\partial f}{\partial \mathbf{x}_T}\bigg]\\ \label{first-derivative-objective-b}
&\big[\frac{\partial f}{\mathbf{\Theta}_t} \frac{\partial f}{\partial \boldsymbol{\mathcal{V}}_{t}}  \frac{\partial f}{\partial \boldsymbol{\mathcal{P}}_{t}}  \frac{\partial f}{\partial \boldsymbol{\mathcal{Q}}_{t}}    \frac{\partial f}{\partial (\boldsymbol{\mathcal{SOC}}_{t})}  \frac{\partial f}{\partial \boldsymbol{\mathcal{P}}_{t}^{\mathrm{ch}}}  \frac{\partial f}{\partial \boldsymbol{\mathcal{P}}_{t}^{\mathrm{dch}}}\frac{\partial f}{\partial \boldsymbol{\mathcal{Q}}_{t}^{\mathrm{s}}}\big]
\end{alignat}
\end{subequations}
\begin{equation}
\begin{multlined}
\label{eqn:second-derivative-objective}
f_{\mathbf{X}\mathbf{X}}=\frac{\partial^2 f}{\partial \mathbf{X}^2}=\frac{\partial}{\partial \mathbf{X}}({\frac{\partial f}{\partial \mathbf{X}}})^\top=
\begin{bmatrix}
    \frac{\partial^2 f}{\partial \mathbf{x}_{1}^2} & \dots  &\frac{\partial^2 f}{\partial \mathbf{x}_{1}\mathbf{x}_{n}}  \\
    \vdots & \ddots & \vdots \\
    \frac{\partial^2 f}{\partial \mathbf{x}_{n}\mathbf{x}_{1}}  & \dots  & \frac{\partial^2 f}{\partial \mathbf{x}_{n}^2} 
\end{bmatrix}
\end{multlined}
\end{equation}
Eqs. \eqref{first-derivative-objective-a}, \eqref{first-derivative-objective-b} and \eqref{eqn:second-derivative-objective} are the basic forms of first and second derivatives of objective function which is $ f:{\rm I\!R}^n\rightarrow\mathbb{C}$. However, constraints $\mathbf{G}(\mathbf{X})$ and $\mathbf{H}(\mathbf{X})$ are vector functions $ f:{\rm I\!R}^n\rightarrow\mathbb{R}^m$ and therefore: 
\begin{equation}
\begin{multlined}
\label{eqn:G(X)}
\mathbf{G}(\mathbf{X})=
{\begin{bmatrix}
    \mathbf{g}_1(\mathbf{X})  & \mathbf{g}_2(\mathbf{X}) &\dots & \mathbf{g}_{N_g}(\mathbf{X})
\end{bmatrix}^\top
 \ \mkern-10mu}_{1 \times N_g}
\end{multlined}
\end{equation}
First derivative of this complex vector function can be written as:
\begin{equation}
\begin{multlined}
\label{eqn:first-derivative-G}
\mathbf{G}_{\mathbf{X}}=\frac{\partial \mathbf{G}}{\partial \mathbf{X}}=
{\begin{bmatrix}
    \frac{\partial \mathbf{g}_1}{\partial \mathbf{x}_{1}} & \dots& \frac{\partial \mathbf{g}_1}{\partial \mathbf{x}_{t}}& \dots &\frac{\partial \mathbf{g}_1}{\partial \mathbf{x}_{T}}  \\
    \vdots & \ddots & \vdots & \ddots & \vdots \\
    \frac{\partial \mathbf{g}_k}{\partial \mathbf{x}_{1}}  & \dots & \frac{\partial \mathbf{g}_k}{\partial \mathbf{x}_{t}}&\dots & \frac{\partial \mathbf{g}_k}{\partial \mathbf{x}_{T}} 
\end{bmatrix}
 \ \mkern-10mu}_{N_g \times N_x}
\end{multlined}
\end{equation}
\begin{equation}
\begin{multlined}
\label{eqn:first-derivative-H}
\mathbf{H}_{\mathbf{X}}=\frac{\partial \mathbf{H}}{\partial \mathbf{X}}=
{\begin{bmatrix}
    \frac{\partial \mathbf{h}_1}{\partial \mathbf{x}_{1}} & \dots& \frac{\partial \mathbf{h}_1}{\partial \mathbf{x}_{t}}& \dots &\frac{\partial \mathbf{h}_1}{\partial \mathbf{x}_{T}}  \\
    \vdots & \ddots & \vdots & \ddots & \vdots \\
    \frac{\partial \mathbf{h}_l}{\partial \mathbf{x}_{1}}  & \dots & \frac{\partial \mathbf{h}_l}{\partial \mathbf{x}_{t}}&\dots & \frac{\partial \mathbf{h}_l}{\partial \mathbf{x}_{T}} 
\end{bmatrix}
 \ \mkern-10mu}_{N_h \times N_x}
\end{multlined}
\end{equation}
Calculation of second derivatives might be somewhat confusing since the three-dimensional set of partial derivatives will not be calculated here \cite{zimmerman2010ac}. The reason is fairly simple and straightforward. In this context, we are using a Newton-Raphson method to find where the partials of a Lagrangian are equal to zero. It is the Hessian of the Lagrangian function in \eqref{eqn:lagrangian} that we need to compute and we always compute it with a known lambda vector. Therefore, it is only the partial derivatives w.r.t. $\mathbf{X}$ of the vector resulting from multiplying the transpose of the Jacobian by lambda that are needed in this context, which means:
\begin{equation}
\begin{multlined}
\label{eqn:second-derivative-HG}
\mathbf{G}_{\mathbf{X}\mathbf{X}}=\frac{\partial}{\partial \mathbf{X}}(\mathbf{G}_\mathbf{X}^\top \boldsymbol{\lambda})
\end{multlined}
\end{equation}
\begin{equation}
\begin{multlined}
\label{eqn:second-derivative-GH}
\mathbf{G}_{\mathbf{X}\mathbf{Y}}=\frac{\partial}{\partial \mathbf{Y}}(\mathbf{G}_{\mathbf{X}}^\top \boldsymbol{\lambda})
\end{multlined}
\end{equation}
The same types of derivatives can also be written for $\mathbf{H}(\mathbf{X})$. More details regarding  the first and second partial differentials of $F(\mathbf{X})$, $\mathbf{G}(\mathbf{X})$ and $\mathbf{H}(\mathbf{X})$, and the arrangement of their matrices can be followed in the Appendix \ref{Appendix_B} of this paper.\\
Partial differentials of equality and inequality constraints are elaborated here in this appendix. According to Eqs. \eqref{first-derivative-objective-a}, \eqref{first-derivative-objective-b}, and \eqref{eqn:second-derivative-objective}, first and second derivatives of $F(\mathbf{X})$, and based on Eqs. \eqref{eqn:first-derivative-G}, \eqref{eqn:first-derivative-H}, \eqref{eqn:second-derivative-HG} and \eqref{eqn:second-derivative-GH} first and second derivatives of $\mathbf{G}(\mathbf{X})$ and $\mathbf{H}(\mathbf{X})$ can be extracted. As we introduced the structure of  $\mathbf{G}(\mathbf{X})$ and $\mathbf{H}(\mathbf{X})$ before reordering as  \eqref{multiperiodG.a} and \eqref{multiperiodH.b}, in the subsections below the first and second partial derivative of $F(\mathbf{X})$, $\mathbf{G}(\mathbf{X})$ and $\mathbf{H}(\mathbf{X})$ are analytically extracted. 
\subsection{First Partial Derivatives of Equality Constraints \texorpdfstring{$\mathbf{G}_\mathbf{X}$}{GX}}
\begin{subequations}
\begin{flalign}
\label{eqn:tildegt}
&\mathbf{G}_{\mathbf{X}}=
    {\begin{bmatrix}
   \widetilde{\mathbf{G}}_{\mathbf{X}}=\frac{\partial \widetilde{\mathbf{G}}}{\partial \mathbf{X}}  & \overline{\mathbf{G}}_{\mathbf{X}}=\frac{\partial \overline{\mathbf{G}}}{\partial \mathbf{X}} &
  \overline{\mathbf{G}}^s_{\mathbf{X}}=\frac{\partial \overline{\mathbf{G}}^s}{\partial \mathbf{X}} 
 \end{bmatrix}}^\top&\\ \label{41.b}
&\widetilde{\mathbf{G}}_\mathbf{X}=
    {\begin{bmatrix}
   \widetilde{\mathbf{g}}_{\mathbf{x}_1} \Scale[0.9]{=}\frac{\partial \widetilde{\mathbf{g}}}{\partial \mathbf{x}_1}  & \widetilde{\mathbf{g}}_{\mathbf{x}_2} \Scale[0.9]{=}\frac{\partial \widetilde{\mathbf{g}}}{\partial \mathbf{x}_2} & \Scale[0.9]{\dots}&
  \widetilde{\mathbf{g}}_{\mathbf{x}_T} \Scale[0.9]{=}\frac{\partial \widetilde{\mathbf{g}}}{\partial \mathbf{x}_T} 
 \end{bmatrix}}^\top\\
 &\overline{\mathbf{G}}_\mathbf{X}=
    {\begin{bmatrix}
  \overline{\mathbf{g}}_{\mathbf{x}_1} \Scale[0.9]{=}\frac{\partial \overline{\mathbf{g}}}{\partial \mathbf{x}_1}  &\overline{\mathbf{g}}_{\mathbf{x}_2} \Scale[0.9]{=}\frac{\partial \overline{\mathbf{g}}}{\partial \mathbf{x}_2} & \Scale[0.9]{\dots}&
  \overline{\mathbf{g}}_{\mathbf{x}_T} \Scale[0.9]{=}\frac{\partial \overline{\mathbf{g}}}{\partial \mathbf{x}_T} 
 \end{bmatrix}}^\top\\
 &\overline{\mathbf{G}}^s_\mathbf{X}=
    {\begin{bmatrix}
   \overline{\mathbf{g}}^s_{\boldsymbol{\tau}_1} \Scale[0.9]{=}\frac{\partial \overline{\mathbf{g}}^s}{\partial \boldsymbol{\tau}_1}  & \overline{\mathbf{g}}^s_{\boldsymbol{\tau}_2} \Scale[0.9]{=}\frac{\partial \overline{\mathbf{g}}^s}{\partial \boldsymbol{\tau}_2} & \Scale[0.9]{\dots}&
  \overline{\mathbf{g}}^s_{\boldsymbol{\tau}_T} \Scale[0.9]{=}\frac{\partial \overline{\mathbf{g}}^s}{\partial \boldsymbol{\tau}_T} 
 \end{bmatrix}}^\top
\end{flalign}
\end{subequations}
Recall from Section \ref{Sec:multiAC} and Eqs. \eqref{active_Power_flow2}-\eqref{Reactive_Power_flow2}; thus, the expression of $\widetilde{\mathbf{g}}_{\mathbf{x}_t}$ in \eqref{41.b} can be expanded as:
\begin{align}
\begin{split}
\label{eqn:tildegxt1}
\widetilde{\mathbf{g}}&_{\mathbf{x}_t}=\frac{\partial \widetilde{\mathbf{g}}}{\partial \mathbf{x}_t}=
\begin{bmatrix}
    \Re{\{\widetilde{\mathbf{g}}_{\mathbf{x}_t}\}}\\
    \Im{\{\widetilde{\mathbf{g}}_{\mathbf{x}_t}\}}\\
\end{bmatrix}\\
=&
\begin{bmatrix}
    \Re{\{\widetilde{\mathbf{g}}_{\boldsymbol{\Theta}_t}\  \widetilde{\mathbf{g}}_{\boldsymbol{\mathcal{V}}_t}\  \widetilde{\mathbf{g}}_{\boldsymbol{\mathcal{P}}_t}\  \widetilde{\mathbf{g}}_{\boldsymbol{\mathcal{Q}}_{t}}\ \widetilde{\mathbf{g}}_{\boldsymbol{\mathcal{SOC}}_{t}}\  \widetilde{\mathbf{g}}_{\boldsymbol{\mathcal{P}}_{t}^{\mathrm{ch}}} \ \widetilde{\mathbf{g}}_{\boldsymbol{\mathcal{P}}_{t}^{\mathrm{dch}}} \ \widetilde{\mathbf{g}}_{\boldsymbol{\mathcal{Q}}_{t}^{\mathrm{s}}}\}}\\
\Im{\{\widetilde{\mathbf{g}}_{\boldsymbol{\Theta}_t}\  \widetilde{\mathbf{g}}_{\boldsymbol{\mathcal{V}}_t}\  \widetilde{\mathbf{g}}_{\boldsymbol{\mathcal{P}}_t}\  \widetilde{\mathbf{g}}_{\boldsymbol{\mathcal{Q}}_{t}}\ \widetilde{\mathbf{g}}_{\boldsymbol{\mathcal{SOC}}_{t}}\  \widetilde{\mathbf{g}}_{\boldsymbol{\mathcal{P}}_{t}^{\mathrm{ch}}} \ \widetilde{\mathbf{g}}_{\boldsymbol{\mathcal{P}}_{t}^{\mathrm{dch}}} \ \widetilde{\mathbf{g}}_{\boldsymbol{\mathcal{Q}}_{t}^{\mathrm{s}}}\}}\\
\end{bmatrix}
\end{split}
\end{align}
In Section \ref{sec:single_OPF}, $\mathbf{\underline{V}}$, $\boldsymbol{\mathcal{V}}$ and $\boldsymbol{\Theta}$ are defined as vectors of complex bus voltages, bus voltage magnitudes and angles respectively. In addition, let $\mathbf{\underline{F}}=\mathbf{diag}(\boldsymbol{\mathcal{V}})^{-1}\mathbf{\underline{V}}$, therefore \eqref{eqn:tildegxt1} can be extended as:
\begin{subequations}\label{eqn:tildegxt2}
\begin{alignat}{2}\label{eqn:tildegxt2a}
    &\widetilde{\mathbf{g}}_{\boldsymbol{\Theta}_t}= j\mathbf{diag}(\mathbf{\underline{V}}_{t})\big(\mathbf{diag}({{\mathbf{\underline{I}}^{\mathrm{bus}}}^*_{t}})-{\mathbf{\underline{Y}}^{\mathrm{bus}}}^*\mathbf{diag}(\mathbf{\underline{V}}_{t}^*)\big)\\
    &\widetilde{\mathbf{g}}_{\boldsymbol{\mathcal{V}}_t}= \mathbf{diag}(\mathbf{\underline{V}}_{t})\big(\mathbf{diag}({{\mathbf{\underline{I}}^{\mathrm{bus}}}^*_{t}})\notag\\
    &\qquad -{\mathbf{\underline{Y}}^{\mathrm{bus}}}^*\mathbf{diag}(\mathbf{\underline{V}}_{t}^*)\big)\mathbf{diag}(\boldsymbol{\mathcal{V}}_{t})^{-1}\\    
    &\widetilde{\mathbf{g}}_{\boldsymbol{\mathcal{P}}^g_t}= -\mathbf{C}_t^\mathrm{g}\\
    &\widetilde{\mathbf{g}}_{\boldsymbol{\mathcal{Q}}^g_{t}} =-j\mathbf{C}_t^\mathrm{g}\\
    &\widetilde{\mathbf{g}}_{\boldsymbol{\mathcal{SOC}}_{t}}=0\\
    &\widetilde{\mathbf{g}}_{\boldsymbol{\mathcal{P}}_{t}^{ch}}=-\mathbf{C}_t^\mathrm{ch}\\
    &\widetilde{\mathbf{g}}_{\boldsymbol{\mathcal{P}}_{t}^{dch}}=\mathbf{C}_t^\mathrm{dch}\\
        &\widetilde{\mathbf{g}}_{\boldsymbol{\mathcal{Q}}_{t}^{s}}=-j\mathbf{C}_t^\mathrm{s}\label{eqn:tildegxt2h}
\end{alignat}
\end{subequations}
It was noted in Section \ref{Sec:multiAC} that linear equality constraints consist of: (a) $\theta^{\mathrm{slack}}=0$, (b) any user-defined custom linear constraint, and (c) binding upper and lower bound variables such that $x_t^\mathrm{min}=x_t^\mathrm{max}$. Thus, it can be written in general format as $\overline{\mathbf{g}}(\mathbf{x}_{t})=
    \mathbf{A}^{\mathrm{grid}}_t \ \mathbf{x}_{t}-\mathbf{B}^{\mathrm{{grid}}}_t$, and subsequently,  $\label{eqn:bargxt}
\overline{\mathbf{g}}_{\mathbf{x}_t}=\frac{\partial \overline{\mathbf{g}}}{\partial \mathbf{x}_t}=
    \mathbf{A}^\mathrm{grid}_t$. The same could be extended to equality constraints regarding the energy storage constraints.
\begin{subequations}
\begin{alignat}{3}
&\overline{\mathbf{g}}^s(\boldsymbol{\tau}_1)=
\begin{bmatrix}
    \mathbf{A}^{\mathrm{s}}_{\boldsymbol{\tau}_1}
\end{bmatrix} \begin{bmatrix}
    \mathbf{x}_{1}
\end{bmatrix} - \begin{bmatrix}
   \mathbf{B}^{\mathrm{{s}}}_{\boldsymbol{\tau}_1}
\end{bmatrix} =0\\
&\mathbf{A}^{\mathrm{s}}_{\boldsymbol{\tau}_1}
 =\Scale[0.9]{\begin{bmatrix}
    0 &\dots &0 & \mathbf{E}^\mathrm{max} &-\mathbf{\Psi}^\mathrm{ch}\Delta t&\frac {\Delta t}{\mathbf{\Psi}^\mathrm{dch}}
\end{bmatrix}}\\
&\mathbf{B}^{\mathrm{s}}_{\boldsymbol{\tau}_1}
 =\Scale[0.9]{\begin{bmatrix}
 \mathbf{SOCi}_{1,t=1} \\
 .\\
 \mathbf{SOCi}_{i,t=1}\\
 .\\
 \mathbf{SOCi}_{n_y,t=1}
\end{bmatrix}} \qquad\in \mathbb{R}^{n_y\times1}\\
&\overline{\mathbf{g}}^s(\boldsymbol{\tau}_t)=
\begin{bmatrix}
    \mathbf{A}^{\mathrm{s}}_{\boldsymbol{\tau}_t}
\end{bmatrix} \begin{bmatrix}
    \mathbf{x}_{t-1}\\
    \mathbf{x}_{t}
\end{bmatrix} - \begin{bmatrix}
   \mathbf{B}^{\mathrm{{s}}}_{\boldsymbol{\tau}_t}
\end{bmatrix} =0\\
&\mathbf{A}^{\mathrm{s}}_{\boldsymbol{\tau}_t}=\notag\\
&\Scale[0.9]{\begin{bmatrix}
    0 \ \dots \ 0 &-\mathbf{E}^\mathrm{max}&0\ \dots\ 0 &  \mathbf{E}^\mathrm{max} &-\mathbf{\Psi}^\mathrm{ch}\Delta t&\frac {\Delta t}{\mathbf{\Psi}^\mathrm{dch}}
\end{bmatrix}}\\
&\mathbf{B}^{\mathrm{s}}_{\boldsymbol{\tau}_t}
 =\Scale[0.9]{\begin{bmatrix}
 \mathbf{SOCi}_{1,t} \\
 .\\
 \mathbf{SOCi}_{i,t}\\
 .\\
 \mathbf{SOCi}_{n_y,t}
\end{bmatrix}} \qquad  \in \mathbb{R}^{n_y\times1}
\end{alignat}
\end{subequations}
Therefore, $\overline{\mathbf{G}}^s_{\boldsymbol{\tau}_t}=\mathbf{A}^{\mathrm{s}}_{\boldsymbol{\tau}_t}$.
\subsection{Second Partial Derivatives of Equality Constraints-\texorpdfstring{$\widetilde{\mathbf{G}}_\mathbf{XX}$}{}}
Second partial derivative of equality constraints w.r.t variables are called in \eqref{LXX}, and can be computed analytically with the format shown in Eqs. \eqref{eqn:second-derivative-HG} and \eqref{eqn:second-derivative-GH}. In detail, they are expanded as shown below. 
\begin{equation}
{\mathbf{G}}_{\mathbf{x}\mathbf{x}}=\frac{\partial}{\partial \mathbf{x}}({\frac{\mathbf{G}_{\mathbf{x}}^\top \boldsymbol{\lambda}}{\partial \mathbf{X}}})
\end{equation}
Therefore 
\begin{subequations}
\begin{alignat}{2}
\label{eqn:tildegxt4}
    &\widetilde{\mathbf{g}}_{\boldsymbol{\Theta}\boldsymbol{\Theta}}(\boldsymbol{\lambda})\bigg\rvert_{t}= \mathbf{diag}(\mathbf{\underline{V}}_{t}^*)\big({{\mathbf{\underline{Y}}^\mathrm{bus}}^*}^\top\mathbf{diag}(\mathbf{\underline{V}}_{t})\mathbf{diag}(\boldsymbol{\lambda})\notag\\
    &-\mathbf{diag}\big({{\mathbf{\underline{Y}}^\mathrm{bus}}^*}^\top\mathbf{diag}(\mathbf{\underline{V}}_{t})\boldsymbol{\lambda}\big)\big)\notag\\
    &\Scale[0.95]{+\mathbf{diag}(\boldsymbol{\lambda})\mathbf{diag}(\mathbf{\underline{V}}_{t})({\mathbf{\underline{Y}}^{\mathrm{bus}}}^*\mathbf{diag}(\mathbf{\underline{V}}_{t}^*)-\mathbf{diag}({\mathbf{\underline{I}}^{\mathrm{bus}}_{t}}^*))}\\
      &\widetilde{\mathbf{g}}_{\boldsymbol{\mathcal{V}}\boldsymbol{\Theta}}(\boldsymbol{\lambda})\bigg\rvert_{t}=\frac{\partial}{\partial \boldsymbol{\mathcal{V}}}(\widetilde{\mathbf{g}}_{\boldsymbol{\Theta}}^\top \boldsymbol{\lambda})
      =j\mathbf{diag}(\boldsymbol{\mathcal{V}}_{t})^{-1}\notag\\
      &\bigg(\big(\mathbf{diag}(\mathbf{\underline{V}}_{t}^*)({{\mathbf{\underline{Y}}^\mathrm{bus}}^*}^\top\mathbf{diag}(\mathbf{\underline{V}}_{t})\mathbf{diag}(\boldsymbol{\lambda})\notag\\
      &-\mathbf{diag}\big({{\mathbf{\underline{Y}}^\mathrm{bus}}^*}^\top\mathbf{diag}(\mathbf{\underline{V}}_{t})\boldsymbol{\lambda}\big)\big)\notag\\
    &\Scale[.95]{-\mathbf{diag}(\boldsymbol{\lambda})\mathbf{diag}(\mathbf{\underline{V}}_{t})({{\mathbf{\underline{Y}}^\mathrm{bus}}^*}\mathbf{diag}(\mathbf{\underline{V}}_{t}^*)-\mathbf{diag}({{\mathbf{\underline{I}}^\mathrm{bus}}^*}))\bigg)}\\
     &\widetilde{\mathbf{g}}_{\boldsymbol{\Theta}\boldsymbol{\mathcal{V}}}(\boldsymbol{\lambda})\bigg\rvert_{t}=\bigg(\widetilde{\mathbf{g}}_{\boldsymbol{\mathcal{V}}\boldsymbol{\Theta}}(\boldsymbol{\lambda})\bigg\rvert_{t}\bigg)^\top\\
      &\widetilde{\mathbf{g}}_{\boldsymbol{\mathcal{V}}\boldsymbol{\mathcal{V}}}(\boldsymbol{\lambda})\bigg\rvert_{t}= \mathbf{diag}(\boldsymbol{\mathcal{V}}_{t})^{-1}\notag\\
      &\bigg(\mathbf{diag}(\boldsymbol{\lambda})\mathbf{diag}(\mathbf{\underline{V}}_{t}){{{\mathbf{\underline{Y}}^\mathrm{bus}}^*}}\mathbf{diag}(\mathbf{\underline{V}}_{t}^*)\notag\\
       &\Scale[0.95]{ +\mathbf{diag}(\mathbf{\underline{V}}_{t}^*){{\mathbf{\underline{Y}}^\mathrm{bus}}^*}^\top\mathbf{diag}(\mathbf{\underline{V}}_{t})\mathbf{diag}(\boldsymbol{\lambda})\bigg)\mathbf{diag}(\boldsymbol{\mathcal{V}}_{t})^{-1}}
\end{alignat}
\end{subequations}
All other partial second derivatives are zero, since the first partial derivatives calculated before were also zero as shown in Eqs. \eqref{eqn:tildegxt2a}-\eqref{eqn:tildegxt2h}. In addition, second derivatives of linear equality constraints are zero, as shown here.  
\begin{subequations}
\begin{flalign}
\label{eqn:second-derivative-BarG}
&\overline{\mathbf{g}}_{\mathbf{x}\mathbf{x}}\bigg\rvert_{t}=\frac{\partial}{\partial \mathbf{x}}({\frac{\overline{\mathbf{g}}_{\mathbf{x}}^\top \boldsymbol{\lambda}}{\partial \mathbf{x}}})\bigg\rvert_{t}=0&\\
&\overline{\mathbf{g}}^s_{\mathbf{x}\mathbf{x}}\bigg\rvert_{t}=\frac{\partial}{\partial \mathbf{x}}({\frac{{\overline{\mathbf{g}}^s_{\mathbf{x}}}^\top \boldsymbol{\lambda}}{\partial \mathbf{x}}})\bigg\rvert_{t}=0&
\end{flalign}
\end{subequations}
\subsection{First Partial Derivatives of Inequality Constraints \texorpdfstring{${\mathbf{H}}_\mathbf{X}$}{}}
\begin{subequations}
\begin{flalign}
\label{eqn:Htildegt}
&\mathbf{H}_{\mathbf{X}}=
    {\begin{bmatrix}
   \widetilde{\mathbf{H}}_{\mathbf{X}}=\frac{\partial \widetilde{\mathbf{H}}}{\partial \mathbf{X}}  & \overline{\mathbf{H}}_{\mathbf{X}}=\frac{\partial \overline{\mathbf{H}}}{\partial \mathbf{X}}  
 \end{bmatrix}}^\top\\ \label{49.b}
&\widetilde{\mathbf{H}}_\mathbf{X}=
    {\begin{bmatrix}
   \widetilde{\mathbf{h}}_{\mathbf{x}_1} \Scale[0.9]{=}\frac{\partial \widetilde{\mathbf{h}}}{\partial \mathbf{x}_1}  & \widetilde{\mathbf{h}}_{\mathbf{x}_2} \Scale[0.9]{=}\frac{\partial \widetilde{\mathbf{h}}}{\partial \mathbf{x}_2} & \Scale[0.9]{\dots}&
  \widetilde{\mathbf{h}}_{\mathbf{x}_T} \Scale[0.9]{=}\frac{\partial \widetilde{\mathbf{h}}}{\partial \mathbf{x}_T} 
 \end{bmatrix}}^\top\\
 &\overline{\mathbf{H}}_\mathbf{X}=
    {\begin{bmatrix}
  \overline{\mathbf{h}}_{\mathbf{x}_1} \Scale[0.9]{=}\frac{\partial \overline{\mathbf{h}}}{\partial \mathbf{x}_1}  &\overline{\mathbf{h}}_{\mathbf{x}_2} \Scale[0.9]{=}\frac{\partial \overline{\mathbf{h}}}{\partial \mathbf{x}_2} & \Scale[0.9]{\dots}&
  \overline{\mathbf{h}}_{\mathbf{x}_T} \Scale[0.9]{=}\frac{\partial \overline{\mathbf{h}}}{\partial \mathbf{x}_T} 
 \end{bmatrix}}^\top
\end{flalign}
\end{subequations}
$\widetilde{\mathbf{H}}(\mathbf{X})$ is derived in \eqref{eqn:h_t}. As elaborated in Section \ref{sec:single_OPF}, $\mathbf{\underline{S}}^{\mathrm{Line}}$ can be extended as $\mathbf{\underline{S}}^{\mathrm{Line}}=
\begin{bmatrix}
    \mathbf{\underline{S}}^{\mathrm{fr}} \\
    \mathbf{\underline{S}}^{\mathrm{to}} \\
\end{bmatrix}$. With the same procedure, $\widetilde{\mathbf{h}}(\mathbf{x}_{t})$ can be extended as:
\begin{align}
\begin{split}
\label{eqn:tildeht}
\widetilde{\mathbf{h}}(\mathbf{x}_{t})=
\begin{bmatrix}
    \widetilde{\mathbf{h}}^{\mathrm{fr}}(\mathbf{x}_t)\\
    \widetilde{\mathbf{h}}^{\mathrm{to}}(\mathbf{x}_t)\\
\end{bmatrix}
= \begin{bmatrix}
   (\mathbf{\underline{S}}^{\mathrm{fr}})^*\mathbf{\underline{S}}^{\mathrm{fr}}-(\abs{\mathbf{\underline{S}}^{\mathrm{Line}}_\mathrm{max}})^{2}\\
   (\mathbf{\underline{S}}^{\mathrm{to}})^*\mathbf{\underline{S}}^{\mathrm{to}}-(\abs{\mathbf{\underline{S}}^{\mathrm{Line}}_\mathrm{max}})^{2}
\end{bmatrix}
\end{split}
\end{align}
Therefore, first derivative of the $\widetilde{\mathbf{h}}(\mathbf{x}_{t})$ can be written as:
\begin{subequations}
\begin{flalign}\label{eqn:tildehxt2}
\widetilde{\mathbf{h}}_{\mathbf{x}_t}&=\frac{\partial \widetilde{\mathbf{h}}}{\partial \mathbf{x}_t}=
\begin{bmatrix}
    \widetilde{\mathbf{h}}_{\mathbf{x}_t}^{\mathrm{fr}}\\
   \widetilde{\mathbf{h}}_{\mathbf{x}_t}^{\mathrm{to}}\\
\end{bmatrix}\\
\widetilde{\mathbf{h}}_{\mathbf{x}_t}^{fr}&=\frac{\partial \widetilde{\mathbf{h}}^{\mathrm{fr}}}{\partial \mathbf{x}_t}=
    2(\Re{\{\mathbf{diag}(\mathbf{\underline{S}}^{\mathrm{fr}})\}}\Re{\{\mathbf{\underline{S}}_{\mathbf{x}_t}^{\mathrm{fr}}\}}&\notag\\
    &+\Im{\{\mathbf{diag}(\mathbf{\underline{S}}^{\mathrm{fr}})\}}\Im{\{\mathbf{\underline{S}}_{\mathbf{x}_t}^{\mathrm{fr}}\}}&\\
    \widetilde{\mathbf{h}}_{\mathbf{x}_t}^{\mathrm{to}}&=\frac{\partial \widetilde{\mathbf{h}}^{\mathrm{to}}}{\partial \mathbf{x}_t} =2(\Re{\{\mathbf{diag}(\mathbf{\underline{S}}^{\mathrm{to}})\}}\Re{\{\mathbf{\underline{S}}_{\mathbf{x}_t}^{\mathrm{to}}\}}&\notag\\
    &+\Im{\{\mathbf{diag}(\mathbf{\underline{S}}^{\mathrm{to}}
    )\}}\Im{\{\mathbf{\underline{S}}_{\mathbf{x}_{t}}^{\mathrm{to}}\}}&
\end{flalign}
\end{subequations}
where $\mathbf{\underline{S}}^{\mathrm{fr}}$ and $\mathbf{\underline{S}}_{\mathbf{x}_t}^{\mathrm{fr}}$ can be written as follows (note that these equations can be extended for $\mathbf{\underline{S}}^{\mathrm{to}}$ and $\mathbf{\underline{S}}_{\mathbf{x}_t}^{\mathrm{to}}$ with the same format):
\begin{align}
\begin{split}
\label{eqn:barhxt1}
\mathbf{\underline{S}}^{\mathrm{fr}}=\mathbf{diag}(\mathbf{\underline{V}}_{t}^{\mathrm{fr}}){\mathbf{\underline{I}}_{t}^{\mathrm{fr}}}^*
\end{split}
\end{align}
\begin{align}
\begin{split}
\label{eqn:barhxt2}
\mathbf{\underline{I}}_{t}^{\mathrm{fr}}=\mathbf{\underline{Y}}^{\mathrm{fr}}\mathbf{\underline{V}}_{t}
\end{split}
\end{align}
\begin{subequations}
\begin{flalign}
\label{eqn:tildegxt3-a}
    &\mathbf{\underline{S}}_{\boldsymbol{\Theta}_t}^{\mathrm{fr}}= j\big(\mathbf{diag}({\mathbf{\underline{I}}_{t}^{\mathrm{fr}}}^*)\mathbf{C}^{\mathrm{fr}}\mathbf{diag}(\mathbf{\underline{V}}_{t})\notag&\\
    &\qquad\qquad-\mathbf{diag}(\mathbf{C}^{\mathrm{fr}}\mathbf{\underline{V}}_{t}){\mathbf{\underline{Y}}^{\mathrm{fr}}}^*\mathbf{diag}(\mathbf{\underline{V}}_{t}^*)\big)&\\
   &\mathbf{\underline{S}}_{\boldsymbol{\mathcal{V}}_t}^{\mathrm{fr}}= \mathbf{diag}({\mathbf{\underline{I}}_{t}^{\mathrm{fr}}}^*)\mathbf{C}^{\mathrm{fr}}\mathbf{diag}(\mathbf{\underline{F}}_{t})\notag&\\
   &\qquad\qquad-\mathbf{diag}(\mathbf{C}^{\mathrm{fr}}\mathbf{\underline{V}}_{t}){\mathbf{\underline{Y}}^{\mathrm{fr}}}^*\mathbf{diag}(\mathbf{\underline{F}}_{t}^*)&\\
    &\mathbf{\underline{S}}_{\boldsymbol{\mathcal{P}}^g_t}^{\mathrm{fr}}=0&\\
    &\mathbf{\underline{S}}_{\boldsymbol{\mathcal{Q}}^g_t}^{\mathrm{fr}}=0&\\
    &\mathbf{\underline{S}}_{\boldsymbol{\mathcal{SOC}}_t}^{\mathrm{fr}}=0&\\
    &\mathbf{\underline{S}}_{\boldsymbol{\mathcal{P}}_{t}^{\mathrm{ch}}}^{\mathrm{fr}}=0&\\
    &\mathbf{\underline{S}}_{\boldsymbol{\mathcal{P}}_{t}^{\mathrm{dch}}}^{\mathrm{fr}}=0&\\
    &\mathbf{\underline{S}}_{\boldsymbol{\mathcal{Q}}_{t}^{\mathrm{s}}}^{\mathrm{fr}}=0&\label{eqn:tildegxt3-h}
\end{flalign}
\end{subequations}
\begin{align}
\begin{split}
\label{eqn:barht}
\overline{\mathbf{h}}(\mathbf{x}_{t})=
    \begin{bmatrix}
       \mathbf{A}^{\mathrm{BOX}} \ \mathbf{x}_{t}-\mathbf{B}^{\mathrm{BOX}}
    \end{bmatrix}\leq 0
\end{split}
\end{align}
Thus, $\overline{\mathbf{h}}_{\mathbf{x}_{t}} = \mathbf{A}^{\mathrm{BOX}}$.
\subsection{Second Partial Derivatives of Inequality Constraints \texorpdfstring{${\mathbf{H}}_\mathbf{XX}$}{}}
For $\widetilde{\mathbf{h}}_{\mathbf{x}\mathbf{x}}^{\mathrm{fr}}$, in general it can be written (note that these equations can be extended to $\widetilde{\mathbf{h}}_{\mathbf{x}\mathbf{x}}^{\mathrm{to}}$ with the same logic):

\begin{align}
\begin{split}
\label{eqn:tildegxt5}
&\widetilde{\mathbf{h}}_{\mathbf{x}\mathbf{x}}^{\mathrm{fr}}\bigg\rvert_{t}=\frac{\partial}{\partial \mathbf{x}}({\frac{\widetilde{\mathbf{h}}_{\mathbf{x}}^{\mathrm{fr}}{}^\top \boldsymbol{\mu}}{\partial \mathbf{x}}})\bigg\rvert_{t}\\
&=\frac{\partial}{\partial \mathbf{x}}\bigg(\mathbf{\underline{S}}_{\mathbf{x}_t}^{\mathrm{fr}}{}^\top\mathbf{diag}({\mathbf{\underline{S}}^{\mathrm{fr}}}^*)\boldsymbol{\mu}+{{\mathbf{\underline{S}}_{\mathbf{x}_t}^{\mathrm{fr}}}^*}^\top\mathbf{diag}(\mathbf{\underline{S}}^{\mathrm{fr}})\boldsymbol{\mu}\bigg)\bigg\rvert_{t}\\
&=2.\Re{\bigg\{\mathbf{\underline{S}}_{\mathbf{x}\mathbf{x}}^{\mathrm{fr}}\bigg\rvert_{t}\mathbf{diag}({\mathbf{\underline{S}}^{\mathrm{fr}}}^*)\boldsymbol{\mu}+{\mathbf{\underline{S}}_{\mathbf{x}_t}^{\mathrm{fr}}}^\top\mathbf{diag}(\boldsymbol{\mu}){{\mathbf{\underline{S}}_{\mathbf{x}_t}^{\mathrm{fr}}}^*}}\bigg\}
\end{split}
\end{align}

\begin{subequations}
\begin{flalign}
\label{eqn:tildegxt6}
    &\widetilde{\mathbf{h}}_{\boldsymbol{\Theta}\boldsymbol{\Theta}}^{\mathrm{fr}}(\boldsymbol{\mu})\bigg\rvert_{t}\notag\\
    &= 2\Re{\bigg\{\mathbf{\underline{S}}_{\boldsymbol{\Theta}\boldsymbol{\Theta}}^{\mathrm{fr}}\bigg\rvert_{t}\mathbf{diag}({\mathbf{\underline{S}}^{\mathrm{fr}}}^*)\boldsymbol{\mu}+{\mathbf{\underline{S}}_{\boldsymbol{\Theta}}^{\mathrm{fr}}}^\top\mathbf{diag}(\boldsymbol{\mu}){{\mathbf{\underline{S}}_{\boldsymbol{\Theta}}^{\mathrm{fr}}}^*}}\bigg\}\\
      &\widetilde{\mathbf{h}}_{\boldsymbol{\mathcal{V}}\boldsymbol{\Theta}}^{\mathrm{fr}}(\boldsymbol{\mu})\bigg\rvert_{t}\notag\\
      &=2\Re{\bigg\{\mathbf{\underline{S}}^{\mathrm{fr}}_{\boldsymbol{\mathcal{V}}\boldsymbol{\Theta}}\bigg\rvert_{t}\mathbf{diag}({\mathbf{\underline{S}}^{\mathrm{fr}}}^*)\boldsymbol{\mu}+
      {\mathbf{\underline{S}}_{\boldsymbol{\mathcal{V}}}^\mathrm{fr}}^\top\mathbf{diag}(\boldsymbol{\mu}){{\mathbf{\underline{S}}^{\mathrm{fr}}_{\boldsymbol{\Theta}}}^*}}\bigg\}\\
     &{\widetilde{\mathbf{h}}_{\boldsymbol{\Theta}\boldsymbol{\mathcal{V}}}}^{\mathrm{fr}}(\boldsymbol{\mu})\bigg\rvert_{t}\notag\\&=
     2\Re{\bigg\{\mathbf{\underline{S}}^{\mathrm{fr}}_{\boldsymbol{\Theta}\boldsymbol{\mathcal{V}}}\bigg\rvert_{t}\mathbf{diag}({\mathbf{\underline{S}}^{\mathrm{fr}}}^*)\boldsymbol{\mu}+{\mathbf{\underline{S}}_{\boldsymbol{\Theta}}^{\mathrm{fr}}}^\top\mathbf{diag}(\boldsymbol{\mu}){{\mathbf{\underline{S}}_{\boldsymbol{\mathcal{V}}}^{\mathrm{fr}}}^*}}\bigg\}\\
      &\widetilde{\mathbf{h}}_{\boldsymbol{\mathcal{V}}\boldsymbol{\mathcal{V}}}^{\mathrm{fr}}(\boldsymbol{\mu})\bigg\rvert_{t}\notag\\&=
      2\Re{\bigg\{\mathbf{\underline{S}}_{\boldsymbol{\mathcal{V}}\boldsymbol{\mathcal{V}}}^{\mathrm{fr}}\bigg\rvert_{t}\mathbf{diag}({\mathbf{\underline{S}}^{\mathrm{fr}}}^*)\boldsymbol{\mu}+{\mathbf{\underline{S}}_{\boldsymbol{\mathcal{V}}}^{\mathrm{fr}}}^\top\mathbf{diag}(\boldsymbol{\mu}){{\mathbf{\underline{S}}_{\boldsymbol{\mathcal{V}}}^{\mathrm{fr}}}^*}}\bigg\}
\end{flalign}
\end{subequations}
All other second partial derivatives are zero, since the first partial derivatives calculated before were also zero as shown in Eqs. \eqref{eqn:tildegxt2a}-\eqref{eqn:tildegxt2h}.
 $\mathbf{\underline{S}}_{\boldsymbol{\Theta}\boldsymbol{\Theta}}^{\mathrm{fr}}$, $\mathbf{\underline{S}}_{\boldsymbol{\mathcal{V}}\boldsymbol{\Theta}}^{\mathrm{fr}}$, $\mathbf{\underline{S}}_{\boldsymbol{\Theta}\boldsymbol{\mathcal{V}}}^{\mathrm{fr}}$, and $\mathbf{\underline{S}}_{\boldsymbol{\mathcal{V}}\boldsymbol{\mathcal{V}}}^{\mathrm{fr}}$ can be extracted as formulated in Eqs. \eqref{eqn:second_der_S-a}-\eqref{eqn:second_der_S-d}. Note that these equations can be extended to $\mathbf{\underline{S}}_{\boldsymbol{\Theta}\boldsymbol{\Theta}}^{\mathrm{to}}$, $\mathbf{\underline{S}}_{\boldsymbol{\mathcal{V}}\boldsymbol{\Theta}}^{\mathrm{to}}$, $\mathbf{\underline{S}}_{\boldsymbol{\Theta}\boldsymbol{\mathcal{V}}}^{\mathrm{to}}$ and $\mathbf{\underline{S}}_{\boldsymbol{\mathcal{V}}\boldsymbol{\mathcal{V}}}^{\mathrm{to}}$ with the same logic.\\
Second derivatives of linear inequalities are zero
\begin{equation}
\begin{multlined}
\label{eqn:second-derivative-barH}
\overline{h}_{\mathbf{x}\mathbf{x}}\bigg\rvert_{t}=\frac{\partial}{\partial \mathbf{x}}({\frac{\overline{\mathbf{h}}_{\mathbf{x}}^\top \boldsymbol{\lambda}}{\partial \mathbf{x}}})\bigg\rvert_{t}=0
\end{multlined}
\end{equation}
\begin{subequations}
\begin{flalign}
\label{eqn:second_der_S-a}
    \mathbf{\underline{S}}_{\boldsymbol{\Theta}\boldsymbol{\Theta}}^{\mathrm{fr}}(\boldsymbol{\mu})\bigg\rvert_{t}&=\mathbf{diag}(\underline{\mathbf{V}}_{t}^*){{\underline{\mathbf{Y}}^{\mathrm{fr}}}^*}^\top\mathbf{diag}(\boldsymbol{\mu})\mathbf{C}^{\mathrm{fr}}\mathbf{diag}(\underline{\mathbf{V}}_{t})\notag\\
    +&\mathbf{diag}(\underline{\mathbf{V}}_{t}){\mathbf{C}^{\mathrm{fr}}}^\top\mathbf{diag}(\boldsymbol{\mu}){\underline{\mathbf{Y}}^{\mathrm{fr}}}^*\mathbf{diag}(\underline{\mathbf{V}}_{t}^*)\notag
    \\-&\mathbf{diag}\bigg({{\underline{\mathbf{Y}}^{\mathrm{fr}}}^*}^\top\mathbf{diag}(\boldsymbol{\mu})\mathbf{C}^{\mathrm{fr}}\underline{\mathbf{V}}_{t}\bigg)\mathbf{diag}(\underline{\mathbf{V}}_{t}^*)\notag\\
    -&\mathbf{diag}\bigg({\mathbf{C}^{\mathrm{fr}}}^\top\mathbf{diag}(\boldsymbol{\mu}){\underline{\mathbf{Y}}^{\mathrm{fr}}}^*\underline{\mathbf{V}}_{t}^*\bigg)\mathbf{diag}(\underline{\mathbf{V}}_{t})\\
      \mathbf{\underline{S}}_{\boldsymbol{\mathcal{V}}\boldsymbol{\Theta}}^{\mathrm{fr}}(\boldsymbol{\mu})\bigg\rvert_{t}&=j\mathbf{diag}(\boldsymbol{\mathcal{V}}_{t})^{-1}\notag\\
      &\bigg(\mathbf{diag}(\underline{\mathbf{V}}_{t}^*){{\underline{\mathbf{Y}}^{\mathrm{fr}}}^*}^\top\mathbf{diag}(\boldsymbol{\mu})\mathbf{C}^{\mathrm{fr}}\mathbf{diag}(\underline{\mathbf{V}}_{t})\notag\\
      -&\mathbf{diag}(\underline{\mathbf{V}}_{t}){\mathbf{C}^{\mathrm{fr}}}^\top\mathbf{diag}(\boldsymbol{\mu}){\underline{\mathbf{Y}}^{\mathrm{fr}}}^*\mathbf{diag}(\underline{\mathbf{V}}_{t}^*)\notag
   \\-&\mathbf{diag}\bigg({{\underline{\mathbf{Y}}^{\mathrm{fr}}}^*}^\top\mathbf{diag}(\boldsymbol{\mu})\mathbf{C}^{\mathrm{fr}}\underline{\mathbf{V}}_{t}\bigg)\mathbf{diag}(\underline{\mathbf{V}}_{t}^*)\notag\\
   +&\mathbf{diag}\bigg({\mathbf{C}^{\mathrm{fr}}}^\top\mathbf{diag}(\boldsymbol{\mu}){\underline{\mathbf{Y}}^{\mathrm{fr}}}^*\underline{\mathbf{V}}_{t}^*\bigg)\mathbf{diag}(\underline{\mathbf{V}}_{t})\bigg)\\
     \mathbf{\underline{S}}_{\boldsymbol{\Theta}\boldsymbol{\mathcal{V}}}^{\mathrm{fr}}(\boldsymbol{\mu})\bigg\rvert_{t}&={\mathbf{\underline{S}}_{\boldsymbol{\mathcal{V}}\boldsymbol{\Theta}}^{\mathrm{fr}}}^\top(\boldsymbol{\mu})\bigg\rvert_{t}\\
      \mathbf{\underline{S}}_{\boldsymbol{\mathcal{V}}\boldsymbol{\mathcal{V}}}^{\mathrm{fr}}(\boldsymbol{\mu})\bigg\rvert_{t}&=\mathbf{diag}(\boldsymbol{\mathcal{V}}_{t})^{-1}\notag\\
      &\bigg(\mathbf{diag}(\underline{\mathbf{V}}_{t}^*){{\underline{\mathbf{Y}}^{\mathrm{fr}}}^*}^\top\mathbf{diag}(\boldsymbol{\mu})\mathbf{C}^{\mathrm{fr}}\mathbf{diag}(\underline{\mathbf{V}}_{t})\notag\\
      +&\mathbf{diag}(\underline{\mathbf{V}}_{t}){\mathbf{C}^{\mathrm{fr}}}^\top\mathbf{diag}(\boldsymbol{\mu}){\underline{\mathbf{Y}}^{\mathrm{fr}}}^*\mathbf{diag}(\underline{\mathbf{V}}_{t}^*)\bigg)\notag\\
      &\mathbf{diag}(\boldsymbol{\mathcal{V}}_{t})^{-1}\label{eqn:second_der_S-d}
\end{flalign}
\end{subequations}
\subsection{Partial Derivatives of Objective Function \texorpdfstring{$F(\mathbf{X})$}{}}
Since MATPOWER case files are used for the sake of benchmarking the solution proposal, we introduce only quadratic cost functions and their first and second partial derivatives here. Moreover, no operational costs for  storage devices are considered. Thus, we can assume the following function as the total operational cost.
\begin{equation}
    F(\mathbf{X}) = F_1+F_2+\dots+F_t+\dots+F_T
\end{equation}
where $F_1 = F_2= F_t$ and $F_t= f^g_t(\mathcal{P}^g_t)+f^q_t(\mathcal{Q}^g_t)$. Therefore first partial derivatives of $F(\mathbf{X})$ w.r.t. $\mathbf{x}_t$ can be extended as
\begin{subequations}
\begin{flalign}\label{60a}
&F_{\boldsymbol{\theta}_t}=0&\\
&F_{\boldsymbol{\mathcal{V}}_t}=0\\
&F_{\boldsymbol{\mathcal{P}}_t^\mathrm{g}}=\frac{\partial f^\mathrm{g}_t}{\partial \boldsymbol{\mathcal{P}}^\mathrm{g}_t}={f^\mathrm{g}_t}'\\
&F_{\boldsymbol{\mathcal{Q}}_t^g}=\frac{\partial f^\mathrm{q}_t}{\partial \boldsymbol{\mathcal{Q}}^\mathrm{g}_t}={f^\mathrm{q}_t}'\\
&F_{\boldsymbol{\boldsymbol{\mathcal{SOC}}}_t}=0\\
&F_{\boldsymbol{\boldsymbol{\mathcal{P}}}^\mathrm{ch}_t}=0\\
&F_{\boldsymbol{\boldsymbol{\mathcal{P}}}^\mathrm{dch}_t}=0\label{60g}\\
&F_{\boldsymbol{\boldsymbol{\mathcal{Q}}}^\mathrm{s}_t}=0
\end{flalign}
\end{subequations}
and subsequently, the second partial derivatives of $F(\mathbf{X})$ w.r.t. $\mathbf{x}_t$ can be extended using \eqref{60a}-\eqref{60g}. 
\begin{subequations}
\begin{flalign}
&F_{\boldsymbol{\mathcal{P}}^\mathrm{g}_t\boldsymbol{\mathcal{P}}^\mathrm{g}_t}=\frac{\partial {f^\mathrm{g}_t}'}{\partial \boldsymbol{\mathcal{P}}^\mathrm{g}_t}={f^\mathrm{g}_t}{''}&\\
&F_{\boldsymbol{\mathcal{Q}}^\mathrm{g}_t\boldsymbol{\mathcal{Q}}^\mathrm{g}_t}=\frac{\partial {f^\mathrm{q}_t}'}{\partial \boldsymbol{\mathcal{Q}}^\mathrm{g}_t}={f^\mathrm{q}_t}{''}&
\end{flalign}
\end{subequations}
and the rest of the partial derivatives w.r.t. the other variables are zero.
	\section{Sparsity Structure of Partial derivatives} \label{Appendix_C}
In this section, the sparsity structure of computed partial derivatives is illustrated. These structures have a great importance since they contribute to the efficient computational operations when it comes to reordering steps as described in Section \ref{reordering}. The number of non-zero elements in matrices are accurately estimated and the specific amount of memory is pre-allocated to these structures and blocks for further computational purposes. 
\subsection{Sparsity Structure of \texorpdfstring{$\mathbf{G}_\mathbf{X}$}{} and \texorpdfstring{$\mathbf{G}_\mathbf{XX}$}{}}
Overall structure of $\widetilde{\mathbf{G}}_{\mathbf{x}_t}$ is
\begin{subequations}
\begin{flalign}
\setcounter{MaxMatrixCols}{12}
\label{eqn:Firstex-derivative-G}
&{\widetilde{\mathbf{G}}}_{\mathbf{X}}=
\begin{bmatrix}
   \widetilde{\mathbf{G}}_{\mathbf{x}_1}& & &\\
   & \ddots& & \\
   & & \widetilde{\mathbf{G}}_{\mathbf{x}_t}\\
   & & & \ddots&\\
   & & & &\widetilde{\mathbf{G}}_{\mathbf{x}_T}
\end{bmatrix}\\
 &\widetilde{\mathbf{G}}_{\mathbf{x}_t}=\notag\\
&\Scale[0.91]{\begin{bmatrix}
   \Scale[0.8]{\Re{\{\widetilde{\mathbf{g}}_{\boldsymbol{\Theta}}\}}\big\rvert_{t} }&  \Scale[0.8]{\Re{\{\widetilde{\mathbf{g}}_{\boldsymbol{\mathcal{V}}}\}}\big\rvert_{t} } &
    \Scale[0.7]{-}\mathbf{C}_t^{\mathrm{g}}&0&0&\mathbf{C}_t^{\mathrm{ch}}&\Scale[0.7]{-}\mathbf{C}_t^{\mathrm{dch}}&0
    \\
  \Scale[0.8]{ \Im{\{\widetilde{\mathbf{g}}_{\boldsymbol{\Theta}}\}}\big\rvert_{t}}&  \Scale[0.8]{ \Im{\{\widetilde{\mathbf{g}}_{\boldsymbol{\mathcal{V}}}\}}\big\rvert_{t}} & 0 & \Scale[0.7]{-}\mathbf{C}_t^{\mathrm{g}}& 0 & 0&0&\Scale[0.7]{-}\mathbf{C}_t^{\mathrm{s}}\\
\end{bmatrix}}
\end{flalign}
 \end{subequations}
 
In matrix format, $\overline{\mathbf{G}}_\mathbf{X}$ and $\overline{\mathbf{G}}^s_\mathbf{X}$  can be simply written as:
  \begin{flalign}
      &\overline{\mathbf{G}}_\mathbf{X}= \begin{bmatrix}
         \mathbf{A}^{\mathrm{grid}}_{\mathbf{x}_1}& 0 &\dots &0\\
         0& \mathbf{A}^{\mathrm{grid}}_{\mathbf{x}_2}&0&\vdots\\
         \vdots&0&\ddots &   0 \\
         0&\dots&0&\mathbf{A}^{\mathrm{grid}}_{\mathbf{x}_T}
      \end{bmatrix}&
\end{flalign}
 \begin{flalign}
      &\overline{\mathbf{G}}^s_\mathbf{X}= \begin{bmatrix}
         \mathbf{A}^{\mathrm{s}}_{\boldsymbol{\tau}_1}& 0 &\dots &0\\
         0& \mathbf{A}^{\mathrm{s}}_{\boldsymbol{\tau}_2}&0&\vdots\\
         \vdots&0&\ddots &   0 \\
         0&\dots&0&\mathbf{A}^{\mathrm{s}}_{\boldsymbol{\tau}_T}
      \end{bmatrix}&
\end{flalign}
According to \eqref{eqn:second-derivative-HG}, the second derivatives of G and H can be written as:
\begin{subequations}
\begin{flalign}
\label{eqn:second-derivative-G}
&\Scale[0.93]{\widetilde{\mathbf{G}}_{\mathbf{X}\mathbf{X}}\bigg\rvert_{t}=\frac{\partial}{\partial \mathbf{X}}({\frac{\widetilde{\mathbf{g}}_{\mathbf{X}}^\top\boldsymbol{\lambda}}{\partial \mathbf{X}}})\bigg\rvert_{t}=
\begin{bmatrix}
   \widetilde{\mathbf{g}}_{\mathbf{x}\mathbf{x}}\bigg\rvert_{1}  \\
        &\widetilde{\mathbf{g}}_{\mathbf{x}\mathbf{x}}\bigg\rvert_{2}& &      \\
        & &\ddots& \\
    && &\widetilde{\mathbf{g}}_{\mathbf{x}\mathbf{x}}\bigg\rvert_{T} 
\end{bmatrix}}&\\
&\widetilde{\mathbf{g}}_{\mathbf{x}\mathbf{x}}\bigg\rvert_{t} =\begin{bmatrix}
   \widetilde{\mathbf{g}}_{\boldsymbol{\Theta}\boldsymbol{\Theta}}(\boldsymbol{\lambda})\bigg\rvert_{t}&  \widetilde{\mathbf{g}}_{\boldsymbol{\Theta}\boldsymbol{\mathcal{V}}}(\boldsymbol{\lambda})\bigg\rvert_{t}& 0  &0  \\
    \widetilde{\mathbf{g}}_{\boldsymbol{\mathcal{V}}\boldsymbol{\Theta}}(\boldsymbol{\lambda})\bigg\rvert_{t}&  \widetilde{\mathbf{g}}_{\boldsymbol{\mathcal{V}}\boldsymbol{\mathcal{V}}}(\boldsymbol{\lambda})\bigg\rvert_{t}& 0  &0  \\
    0 & 0  & 0&0\\ 
     0 & 0  & 0&0
\end{bmatrix}&
\end{flalign}
\end{subequations}
\subsection{Sparsity Structure of \texorpdfstring{$\mathbf{H}_\mathbf{X}$}{} and \texorpdfstring{$\mathbf{H}_\mathbf{XX}$}{}}

\begin{subequations}
\begin{flalign}
\label{eqn:firstex-derivative-H}
&{\widetilde{\mathbf{H}}}_{\mathbf{X}}=
\begin{bmatrix}
   \widetilde{\mathbf{h}}_{\mathbf{x}_1}& & &\\
   & \ddots& & \\
   & & \widetilde{\mathbf{h}}_{\mathbf{x}_t}\\
   & & & \ddots&\\
   & & & &\widetilde{\mathbf{h}}_{\mathbf{x}_T}
\end{bmatrix}\\
&\widetilde{\mathbf{h}}_{\mathbf{X}_t}={\frac{ \partial \widetilde{\mathbf{h}} }{\partial \mathbf{x}_t}}
=\begin{bmatrix}
        \widetilde{\mathbf{h}}_{\boldsymbol{\Theta}}^{\mathrm{fr}}\big\rvert_{t}
   &  \widetilde{\mathbf{h}}_{\boldsymbol{\mathcal{V}}}^{\mathrm{fr}}\big\rvert_{t}&0 &0\\
    \widetilde{\mathbf{h}}_{\boldsymbol{\Theta}}^{\mathrm{to}}\big\rvert_{t} &   \widetilde{\mathbf{h}}_{\boldsymbol{\mathcal{V}}}^{\mathrm{to}}\big\rvert_{t}    &0&0 \\
    0&0&0&0\\
    0&0&0&0
\end{bmatrix}&
\end{flalign}
\end{subequations}
Finally, we can extend  $H_{\mathbf{X}\mathbf{X}}$ as in \eqref{eqn:second-derivative-HXX}.
\begin{subequations}
\begin{flalign}
\label{eqn:second-derivative-HXX}
&{\mathbf{H}}_{\mathbf{X}\mathbf{X}}=\frac{\partial}{\partial \mathbf{X}}({\frac{{\mathbf{H}}_{\mathbf{X}}^\top \boldsymbol{\lambda}}{\partial \mathbf{X}}})
=\notag\\
&\Scale[0.83]{\begin{bmatrix}
   (\widetilde{\mathbf{h}}_{\mathbf{X}\mathbf{X}}^{\mathrm{fr}}+\widetilde{\mathbf{h}}_{\mathbf{X}\mathbf{X}}^{\mathrm{to}})\bigg\rvert_{1}     \\
                            &(\widetilde{\mathbf{h}}_{\mathbf{X}\mathbf{X}}^{\mathrm{fr}}+\widetilde{\mathbf{h}}_{\mathbf{X}\mathbf{X}}^{\mathrm{to}})\bigg\rvert_{2}&       \\
       & & \ddots& \\
    && &(\widetilde{\mathbf{h}}_{\mathbf{X}\mathbf{X}}^{\mathrm{fr}}+\widetilde{\mathbf{h}}_{\mathbf{X}\mathbf{X}}^{\mathrm{to}})\bigg\rvert_{T}  \\
\end{bmatrix}}\\
\label{eqn:second-derivative-tildeH}
&\Scale[0.91]{\widetilde{\mathbf{h}}_{\mathbf{x}\mathbf{x}}^{\mathrm{fr}}\bigg\rvert_{t}=\frac{\partial}{\partial \mathbf{x}}({\frac{\widetilde{\mathbf{h}}_{\mathbf{x}}^{\mathrm{fr}}{}^\top \boldsymbol{\mu}}{\partial \mathbf{x}}})\bigg\rvert_{t}=}
\begin{bmatrix}
   \Scale[0.9]{\widetilde{\mathbf{h}}_{\boldsymbol{\Theta}\boldsymbol{\Theta}}^{\mathrm{fr}}(\boldsymbol{\mu})\bigg\rvert_{t}}&  \Scale[0.9]{\widetilde{\mathbf{h}}_{\boldsymbol{\Theta}\boldsymbol{\mathcal{V}}}^{\mathrm{fr}}(\boldsymbol{\mu})\bigg\rvert_{t}}& 0  &0  \\
    \Scale[0.9]{\widetilde{\mathbf{h}}_{\boldsymbol{\mathcal{V}}\boldsymbol{\Theta}}^{\mathrm{fr}}(\boldsymbol{\mu})\bigg\rvert_{t}}&  \Scale[0.9]{\widetilde{\mathbf{h}}_{\boldsymbol{\mathcal{V}}\boldsymbol{\mathcal{V}}}^{\mathrm{fr}}(\boldsymbol{\mu})\bigg\rvert_{t}}& 0  &0  \\
    0 & 0  & 0&0\\ 
     0 & 0  & 0&0
\end{bmatrix}\\
&(\widetilde{\mathbf{h}}_{\mathbf{X}\mathbf{X}}^{\mathrm{fr}}+\widetilde{\mathbf{h}}_{\mathbf{X}\mathbf{X}}^{\mathrm{to}})\bigg\rvert_{t}=\notag\\ &\Scale[0.9]{\begin{bmatrix}
   \widetilde{\mathbf{h}}_{\boldsymbol{\Theta}\boldsymbol{\Theta}}^{\mathrm{fr}}(\boldsymbol{\mu})\big\rvert_{t}+\widetilde{\mathbf{h}}_{\boldsymbol{\Theta}\boldsymbol{\Theta}}^{\mathrm{to}}(\boldsymbol{\mu})\big\rvert_{t}&  \widetilde{\mathbf{h}}_{\boldsymbol{\Theta}\boldsymbol{\mathcal{V}}}^{\mathrm{fr}}(\boldsymbol{\mu})\big\rvert_{t}+\widetilde{\mathbf{h}}_{\boldsymbol{\Theta}\boldsymbol{\mathcal{V}}}^{\mathrm{to}}(\boldsymbol{\mu})\big\rvert_{t}&{0} \\
    {\widetilde{\mathbf{h}}_{\boldsymbol{\mathcal{V}}\boldsymbol{\Theta}}^{\mathrm{fr}}(\boldsymbol{\mu})\big\rvert_{t}}+{\widetilde{\mathbf{h}}_{\boldsymbol{\mathcal{V}}\boldsymbol{\Theta}}^{\mathrm{to}}(\boldsymbol{\mu})\big\rvert_{t}}&  {\widetilde{\mathbf{h}}_{\boldsymbol{\mathcal{V}}\boldsymbol{\mathcal{V}}}^{\mathrm{fr}}(\boldsymbol{\mu})\big\rvert_{t}}   +{\widetilde{\mathbf{h}}_{\boldsymbol{\mathcal{V}}\boldsymbol{\mathcal{V}}}^{\mathrm{to}}(\boldsymbol{\mu})\big\rvert_{t}}&{0}  \\
   {0}&  {0}&0 \\
\end{bmatrix}}
\end{flalign}
\end{subequations}
\subsection{Sparsity Structure of \texorpdfstring{${F}_\mathbf{X}$}{} and \texorpdfstring{${F}_\mathbf{XX}$}{}}
First derivative of the objective function has the following structure
\begin{equation}
    {F}_\mathbf{X}= \begin{bmatrix}
           0&0&{f^\mathrm{g}_t}'&{f^\mathrm{q}_t}'&0&0&0&0
    \end{bmatrix}
\end{equation}

Second derivative of $f_{\mathbf{X}\mathbf{X}}$ also can be written as \eqref{eqn:second-derivative-f} according to \eqref{eqn:second-derivative-objective}
\begin{subequations}
\begin{flalign}
\label{eqn:second-derivative-f}
&{f}_{\mathbf{X}\mathbf{X}}
=\begin{bmatrix}
   f_{\mathbf{x}\mathbf{x}}\bigg\rvert_{1}  \\
        &f_{\mathbf{x}\mathbf{x}}\bigg\rvert_{2}& &      \\
        & &\ddots& \\
    && &f_{\mathbf{x}\mathbf{x}}\bigg\rvert_{T} 
\end{bmatrix}\\
&f_{\mathbf{x}\mathbf{x}}\bigg\rvert_{t}=
\begin{bmatrix}
       0
   &  \Scale[1]{\dots}& \Scale[1]{\dots}&\Scale[1]{\dots}&\Scale[1]{\dots}\\
    \Scale[1]{\vdots}&   \Scale[1]{\ddots}    &\Scale[1]{\hdots}  &\Scale[1]{\dots}& \Scale[1]{\dots}\\
     \Scale[1]{\vdots} & \Scale[1]{\vdots} &  \Scale[1]{\begin{bmatrix}f_{PP}\big\rvert_{t}\end{bmatrix}} & \Scale[1]{\dots}&\Scale[1]{\dots}  \\
   \Scale[1]{\vdots}&\Scale[1]{\vdots}&\Scale[1]{\vdots}  &\ \Scale[1]{\begin{bmatrix}f_{QQ}\big\rvert_{t}\end{bmatrix}}&
   \Scale[1]{\dots}\\
    \Scale[1]{\vdots}&\Scale[1]{\vdots}& \Scale[1]{\vdots}  &\Scale[1]{\vdots}& 0  \\
\end{bmatrix}
\end{flalign}
\end{subequations}
	\section{Function Evaluation} \label{Appendix_D}
The efficiency of calculating analytical derivatives and their structures is illustrated here. Table \ref{tab:numericVSanalytic} shows the total computational time in order to calculate $F_\mathbf{X}$, $\mathbf{G}_\mathbf{X}$, $\mathbf{H}_\mathbf{X}$ and $\boldsymbol{\mathcal{L}}_{\mathbf{X}\mathbf{X}}^{\gamma}$ until the algorithm converges with the corresponding iterations, where description of each term can be seen in Table \ref{tab:derivative}. As can be seen, hand-coded derivatives outperform significantly faster than numerical methods. For large networks, numerical derivatives are intractable. The numerical derivatives are computed using central finite differences. Since analytical derivatives are the accurate model of partial derivatives of functions, further accuracy comparison with the finite numerical method applied here is neglected. 
\begin{table}[htbp!]
\caption{First and Second Partial Derivatives}
\label{tab:derivative}
\begin{center}
\begin{tabular}{c c} 
\toprule
 Term & Description  \\  
\midrule \midrule
$F_\mathbf{X}$ & $F_{\mathbf{X}}=\frac{\partial \mathbf{G}}{\partial \mathbf{X}}$ \\
\hline
$\mathbf{G}_\mathbf{X}$ &  $\mathbf{G}_{\mathbf{X}}=\frac{\partial \mathbf{G}}{\partial \mathbf{X}}$\\
\hline
$\mathbf{H}_\mathbf{X}$ & $\mathbf{H}_{\mathbf{X}}=\frac{\partial \mathbf{H}}{\partial \mathbf{X}}$\\
\hline
$\boldsymbol{\mathcal{L}}_{\mathbf{X}\mathbf{X}}^{\gamma}$& $\boldsymbol{\mathcal{L}}_{\mathbf{X}\mathbf{X}}^{\gamma} = F_{\mathbf{XX}}+\mathbf{G}_{\mathbf{XX}}(\boldsymbol{\lambda})+\mathbf{H}_{\mathbf{XX}}(\boldsymbol{\mu})$\\
\bottomrule
\end{tabular}
\end{center}
\end{table}
\begin{table*}
 \caption{Total time ($\mathrm{TotalTime= No_{\cdot} of \ Iter_{\cdot} \times TimePerIter}$) elapsed to calculate: 1) Analytical (hand-coded) derivatives, and 2) Numerical derivatives }
\label{tab:numericVSanalytic}
\begin{threeparttable}
\begin{tabularx}{\textwidth}{m s s s s m s s s m m }
\toprule
&&&&   \multicolumn{3}{c}{Analytical}&&\multicolumn{3}{c}{Numerical} \\
\cmidrule{5-7}  \cmidrule{9-11} 
Case     & $T$&$n_y$ &iter & ${F}_\mathbf{X}$(s)&$\mathbf{G}_\mathbf{X}$+ $\mathbf{H}_\mathbf{X}$(s)&$\boldsymbol{\mathcal{L}}_{\mathbf{X}\mathbf{X}}^{\gamma}$(s)& & ${F}_\mathbf{X}$(s)&$\mathbf{G}_\mathbf{X}$+ $\mathbf{H}_\mathbf{X}$(s)&$\boldsymbol{\mathcal{L}}_{\mathbf{X}\mathbf{X}}^{\gamma}$(s) \\ 
\midrule
Case9     &2 &5&13& 0.03 & 0.13 & 0.14 & &0.43 & 0.98 &140.07\\ 
Case9     &10&5&23& 0.08  & 0.36 & 0.37 & &11.32 & 30.62  & 22815.29\\ 
IEEE30    &2&5&12 & 0.04 & 0.25 & 0.18    &    & 1.01  & 2.16  & 682.70\\ 
IEEE30    &10&5&16& 0.05 & 0.24 & 0.25 && 16.73 & 49.12  & 79712.78\\ 
IEEE118   &2&5&22 & 0.04 & 0.19 & 0.20 &   &7.41 & 18.07 & 24557.09\\ 
IEEE118   &10&5&37&  0.09 & 0.62 & 0.82&    & 158.21\tnote{1}  & 572.09 \tnote{1} & 4599735\tnote{1} \\ 
PEGASE1354&2 &5&23& 0.05 & 0.61 & 0.78  &   &85.54\tnote{1}  & 496.18 \tnote{1} & 7185888 \tnote{1}  \\ 
PEGASE1354&10&5&33& 0.10 & 3.77 & 5.15   && 588.06\tnote{1}  & 3530  \tnote{1} & 51550941\tnote{1}  \\ 
\bottomrule
\end{tabularx}
\begin{tablenotes}
\item[1] {Estimated total time: The time elapsed for one iteration multiplied to the iteration that would take to converge}
\end{tablenotes}
\end{threeparttable}
\end{table*}
\end{appendices}
\end{document}